\documentclass[a4paper,11pt]{article}

\usepackage{amsmath}
\usepackage{amsfonts}
\usepackage{amssymb}
\usepackage{latexsym}
\usepackage{epsfig}
\usepackage{graphicx}
\usepackage{oldgerm}
\usepackage{theorem}
\usepackage{bm}
\usepackage{mathrsfs}

\def\C{\mathbb{C}}
\def\R{\mathbb{R}}
\def\N{\mathbb{N}}
\def\Z{\mathbb{Z}}
\def\HH{\mathbb{H}}
\def\S{\mathbb{S}}

\def\I{\mathbb{I}}
\def\W{\mathbb{W}}
\def\D{\mathbb{D}}
\def\O{\mathbb{O}}
\def\P{\mathbb{P}}
\def\E{\mathbb{E}}
\def\mbK{\mathbb{K}}
\def\J{\mathbb{J}}
\def\A{\mathbb{A}}

\def\Conf{{\rm Conf}}
\def\supp{{\rm supp}\,}
\def\hcap{{\rm hcap}}

\def\rN{{\rm N}}
\def\rP{{\rm P}}
\def\rO{{\rm O}}
\def\rE{{\rm E}}

\def\rC{{\rm C}}
\def\rK{{\rm K}}
\def\rB{{\rm B}}
\def\bE{{\bf E}}
\def\bP{{\bf P}}

\def\bK{{\bf K}}
\def\bW{{\bf W}}
\def\x{\bm{x}}
\def\y{\bm{y}}
\def\z{\bm{z}}

\def\t{\bm{t}}
\def\f{\bm{f}}
\def\bv{\bm{v}}

\def\X{\bm{X}}
\def\Y{\bm{Y}}
\def\1{{\bf 1}}
\def\0{{\bf 0}}
\def\V{\bm{V}}
\def\bxi{\bm{\xi}}
\def\B{\bm{B}}

\def\cC{{\cal C}}

\def\cH{{\cal H}}

\def\cF{{\cal F}}
\def\cI{{\cal I}} 

\def\cB{{\cal B}}

\def\cM{{\cal M}}
\def\cV{{\cal V}}
\def\cD{{\cal D}}
\def\cE{{\cal E}}

\def\cG{{\cal G}}
\def\cZ{{\cal Z}}
\def\cF{{\cal F}}

\def\sS{{\sf S}}

\def\sD{{\sf D}}
\def\sm{{\sf m}}

\def\sL{{\sf L}}
\def\sX{{\sf X}}
\def\sY{{\sf Y}}
\def\sA{{\sf A}}
\def\sB{{\sf B}}
\def\sU{{\sf U}}
\def\sM{{\sf M}}
\def\law={\stackrel{\rm (law)}{=}}
\def\dis={\stackrel{\rm d}{=}}

\def\zbar{\overline{z}}
\def\wbar{\overline{w}}

\def\Re{{\rm Re} \,}
\def\Im{{\rm Im} \,}
\def\arg{{\rm arg}\,}
\def\GFF{\mathrm{GFF}}
\def\SLE{\mathrm{SLE}}
\def\Cov{\mathrm{Cov}}
\def\Var{\mathrm{Var}}

\def\bra{\langle}
\def\ket{\rangle}

\def\mfh{\mathfrak{h}}

\def\alphabar{\overline{\alpha}}

\def\zbar{\overline{z}}
\def\wbar{\overline{w}}

\def\alphahat{\widehat{\alpha}}

\def\perdet{\mathop{ \mathrm{perdet} }}
\def\per{\mathop{ \mathrm{per} }}

\theorembodyfont{\upshape}

\newtheorem{thm}{Theorem}[section]
\newtheorem{lem}[thm]{Lemma}
\newtheorem{cor}[thm]{Corollary}
\newtheorem{prop}[thm]{Proposition}

\newtheorem{df}[thm]{Definition}

\newtheorem{example}[thm]{Example}
\newtheorem{rem}[thm]{Remark}

\newcommand{\SSC}[1]{\section{#1}\setcounter{equation}{0}}
\newcommand{\qed}{\hbox{\rule[-2pt]{3pt}{6pt}}}

\begin{document}

\title{\bf
Point Processes and Multiple SLE/GFF Coupling
\footnote{
This manuscript was prepared for the 
lectures given in
the 4th ZiF Summer School
`Randomness in Physics and Mathematics: 
From Integrable Probability to Disordered Systems'
held at 
ZiF--Center for Interdisciplinary Research, Bielefeld University,
Germany, from 1st to 13th August 2022,
which has been organized by Gernot Akemann
and Friedrich G\"{o}tze.
}
}
\author{
Makoto Katori
\footnote{
Department of Physics,
Faculty of Science and Engineering,
Chuo University,
Kasuga, Bunkyo-ku, Tokyo 112-8551, Japan;
e-mail: katori@phys.chuo-u.ac.jp}
}
\date{12 August 2022}
\pagestyle{plain}
\maketitle

\begin{abstract}
In the series of lectures, we will discuss probability
laws of random points, curves, and surfaces.

Starting from a brief review of the notion of
martingales, one-dimensional Brownian motion (BM),
and the $D$-dimensional Bessel processes, 
BES$_{D}$, $D \geq 1$, 
first we study Dyson's Brownian motion model with parameter
$\beta >0$, DYS$_{\beta}$,
which is a one-parameter family
of repulsively interacting $N$ Brownian motions
on ${\mathbb{R}}$, 
$N \in {\mathbb{N}}:=\{1,2, \dots\}$, 
and is regarded as multivariate extensions of BES$_D$
with the relation $\beta=D-1$. 
In particular, the determinantal structures 
are proved for DYS$_{2}$, which is realized as
the eigenvalue process of Hermitian-matrix-valued BM 
studied in random matrix theory. 

Next, using the reproducing kernels of Hilbert function spaces,
the Gaussian analytic functions (GAFs) are defined on
a unit disk ${\mathbb{D}}$ and an annulus 
${\mathbb{A}}_q :=\{z \in {\mathbb{C}} : 
q < |z| < 1 \}$, $q \in (0, 1)$, 
which provide models of random surfaces.
As zeros of the GAFs, 
determinantal point processes and 
permanental-determinantal point processes 
are obtained, which have symmetry and invariance
associated with conformal transformations. 

Then, the Schramm--Loewner evolution 
with parameter $\kappa >0$, SLE$_{\kappa}$, is
introduced, which is driven by
a BM on ${\mathbb{R}}$ and 
generates a family of
conformally invariant probability laws of 
random curves on
the upper half complex plane ${\mathbb{H}}$. 
We regard SLE$_{\kappa}$ as a complexification of
BES$_D$ with the relation
$\kappa=4/(D-1)$.

The last topic of lectures is 
the construction of the multiple SLE$_{\kappa}$, 
which is driven by the $N$-particle DYS$_{\beta}$ 
on ${\mathbb{R}}$
and generates $N$ interacting random curves
in ${\mathbb{H}}$.
There, we define the Gaussian free field (GFF)
and its generalization called the
imaginary surface with parameter $\chi$,
which are considered as the
distribution-valued random fields on ${\mathbb{H}}$. 
Under the relation 
$\chi=2/\sqrt{\kappa}-\kappa/\sqrt{2}$, 
we characterize the SLE/GFF coupling
studied by Dub\'edat, Sheffield, and Miller 
by its temporal stationarity, and extend 
it to multiple cases. We prove that
the multiple SLE/GFF coupling is established,  
if and only if the driving $N$-particle process
on ${\mathbb{R}}$ is identified with DYS$_{\beta}$ with
the relation $\beta=8/\kappa$.
This relation between parameters is very simple, but
highly nontrivial, since
if we simply combine the relations
$\beta=D-1$ and $\kappa=4/(D-1)$ mentioned above,
we will have a different relation $\beta=4/\kappa$.
Under the present multiple SLE/GFF coupling with
$\beta=8/\kappa$, we can prove that
the multiple SLE driven by DYS$_{\beta}$ exhibits
the phase transitions at 
$\kappa=4$ and $\kappa=8$ in the similar way
to the original SLE$_{\kappa}$ with a single SLE curve. 

\vskip 0.2cm


\end{abstract}
\footnotesize
\tableofcontents
\vspace{3mm}
\normalsize

\SSC
{Determinantal Martingales and Determinantal Stochastic Processes}
\label{sec:stochastic}
\subsection{One-dimensional standard Brownian motion (BM)}
\label{sec:1dimBM}
First of all, we consider the motion of a 
Brownian particle in one-dimensional space $\R$ 
starting from the origin 0 at time $t=0$.
At each time $t>0$, the particle position is randomly distributed, 
and each realization of its path (trajectory) 
is denoted by $\omega$ and called a 
\textbf{sample path}
or simply a 
\textbf{path}.
Let $\Omega$ be the collection of all sample paths
and we call it the 
\textbf{sample path space}.
The position of the Brownian particle at time $t \geq 0$ 
in a path $\omega \in \Omega$
is written as $B(t, \omega)$. Usually we omit $\omega$
and simply write it as $B(t), \ t \geq 0$.

We represent each event associated with the process by a subset of $\Omega$,
and the collection of all events is denoted by $\cF$.
The whole sample path space $\Omega$ and the empty set $\emptyset$
are in $\cF$.
For any two sets $\sA, \sB \in \cF$, we assume that
$\sA \cup \sB \in \cF$ and $\sA \cap \sB \in \cF$.
If $\sA \in \cF$, then its complement $\sA^{\rm c}$ is also in $\cF$.
It is closed for any infinite sum of events in the sense that, 
if $\sA_n \in \cF, n=1,2, \dots$, then
$\cup_{n \geq 1} \sA_n \in \cF$.
Such a collection of events $\cF$ is said to be a 
\textbf{$\sigma$-field} (sigma-field). 
The symbol $\sigma$ is for `sum'.

A \textbf{probability measure} $\rP$ is a nonnegative function
defined on the $\sigma$-field $\cF$.
Since any element of $\cF$ is given by a set as above,
any input of $\rP$ is a set; $\rP$ is a 
\textbf{set function}.
It satisfies the following properties: 
$\rP[\sA] \geq 0$ for all $\sA \in \cF$, 
$\rP[\Omega]=1$, $\rP[\emptyset]=0$, 
and if $\sA, \sB \in \cF$ are disjoint, $\sA \cap \sB = \emptyset$, 
then $\rP[\sA \cup \sB]=\rP[\sA]+\rP[\sB]$.
In particular,
$\rP[\sA^{\rm c}]=1-\rP[\sA]$ for all $\sA \in \cF$.
The triplet $(\Omega, \cF, \rP)$ is called
the \textbf{probability space}.

The smallest $\sigma$-field containing all intervals on $\R$
is called the \textbf{Borel $\sigma$-field} and denoted by $\cB(\R)$.
A \textbf{random variable} or 
\textbf{measurable function}
is a real-valued function $f(\omega)$ on $\Omega$
such that, for every Borel set $A \in \cB(\R)$,
$f^{-1}(A) \in \cF$.
Two events $\sA$ and $\sB$ are said to be 
\textbf{independent}
if $\rP[\sA \cap \sB]=\rP[\sA] \rP[\sB]$.
Two random variables $X$ and $Y$ are 
\textbf{independent}
if the events $\sA=\{X: X \in A\}$ and $\sB=\{Y: Y \in B\}$ are independent
for any $A, B \in \cB(\R)$.

The 
\textbf{one-dimensional standard Brownian motion},
$\{B(t, \omega) : t \geq 0 \}$,
has the following three properties: 
\begin{description}
\item{{\bf (BM1)}} \quad
$B(0, \omega)=0$ with probability one.

\item{{\bf (BM2)}} \quad
There is a subset of the sample path space
$\widetilde{\Omega} \subset \Omega$, such that
$\rP[\widetilde{\Omega}]=1$ and 
for any $\omega \in \widetilde{\Omega}$,
$B(t, \omega)$ is a real continuous function of $t$.
We say that $B(t)$ has 
a 
\textbf{continuous path} almost surely
(a.s., for short).

\item{{\bf (BM3)}} \quad
For an arbitrary $M \in \N := \{1,2,3, \dots\}$,
and for any sequence of times, 
$t_0 := 0 < t_1< \cdots < t_M$, 
the increments $B(t_m)-B(t_{m-1}), m=1,2, \dots, M,$
are independent, 
and each increment is 
in \textbf{the normal distribution} 
(\textbf{the Gaussian distribution}) 
with mean 0 and variance $\sigma^2=t_{m}-t_{m-1}$.
It means that for any $1 \leq m \leq M$ and $a < b$,
$$
\rP[B(t_m)-B(t_{m-1}) \in [a, b]]
=\int_{a}^{b} p(t_m-t_{m-1}, z |0) d z,
$$
where we define for $x,y \in \R$
\begin{equation}
p(t,y|x)= \left\{
\begin{array}{ll}
\displaystyle{\frac{1}{\sqrt{2 \pi t}} e^{-(x-y)^2/2t}},
& \quad \mbox{for $t > 0$},
\cr
\delta(x-y),
& \quad \mbox{for $t=0$}.
\end{array} \right.
\label{eqn:pt1}
\end{equation}
\end{description}

Unless otherwise noted,
the one-dimensional standard Brownian motion is
simply abbreviated to BM in this lecture note.
The probability measure $\rP$ for the BM
in particular is called the 
\textbf{Wiener measure}.
The expectation with respect to the probability measure $\rP$
is denoted by $\rE$.
We write the conditional probability as
$\rP[\cdot | C]$,
where $C$ denotes the condition.
The conditional expectation is similarly
written as $\rE[\cdot | C]$.

The third property {\bf (BM3)} implies that
for any $0 \leq s \leq t < \infty$
\begin{equation}
\rP[B(t) \in A | B(s)=x]
=\int_{A} p(t-s, y|x) dy
\label{eqn:pt1b}
\end{equation}
holds, $^{\forall} A \in \cB(\R), ^{\forall} x \in \R$.
Therefore the integral kernel $p(t, y|x)$ given by (\ref{eqn:pt1}) 
is called the 
\textbf{transition probability density function}
of Brownian motion starting from $x$.
The probability that the BM is observed
in a region $A_m \in \cB(\R)$ at time $t_m$
for each $m=1, 2, \dots, M$ is then given by
\begin{equation}
\rP[B(t_m) \in A_m, m=1,2, \dots, M]
=\int_{A_1}dx_1 \cdots \int_{A_M} dx_M \,
\prod_{m=1}^{M} p(t_m-t_{m-1}, x_m|x_{m-1}),
\label{eqn:BM3}
\end{equation}
where $x_0 := 0$.

By {\bf (BM3)}, we can see that, for any $c>0$,
the probability distribution of $B(c^2t)/c$ is 
equivalent with that of $B(t)$ at arbitrary time $t \geq 0$.
It is written as
$$
\frac{1}{c} B(c^2t) \dis= B(t), \quad ^{\forall} c > 0,
$$
where the symbol $\dis=$ is for 
\textbf{equivalence in distribution}.
Moreover, (\ref{eqn:BM3}) implies that,
for any $c>0$, $B(t), t \geq 0$ and
its time changed process with $t \mapsto c^2t$ multiplied by
a factor $1/c$ (dilatation) follow the same probability law.
This 
\textbf{equivalence in probability law} 
of stochastic processes is expressed as
\begin{equation}
(B(t))_{t \geq 0} \law= \left( \frac{1}{c} B(c^2 t) \right)_{t \geq 0},
\quad ^{\forall} c \geq 0,
\label{eqn:BMscaling}
\end{equation}
and called the 
\textbf{scaling property of Brownian motion}.

For $a > 0$, let
$T_a=\inf\{t \geq 0: B(t)=a \}$. Then for any $t \geq 0$,
\begin{equation}
\rP[T_a <t, B(t)<a]
=\rP[T_a<t, B(t) > a],
\label{eqn:reflection_principle}
\end{equation}
since the transition probability density (\ref{eqn:pt1}) is a symmetric function
of the increment $y-x$.
This property is called the 
\textbf{reflection principle} of BM.
For $\{\omega: B(t) >a\} \subset \{\omega: T_a < t\}$, $a > 0$, 
the above is equal to $\rP[B(t) >a]$.

The formula (\ref{eqn:BM3}) also means that for any fixed $s \geq 0$,
under the condition that $B(s)$ is given,
$\{B(t): t \leq s\}$ and $\{B(t) : t > s\}$ 
are independent.
This independence of the events in the future
and those in the past is called the
\textbf{Markov property}.
A positive random variable $\tau$ is called 
\textbf{stopping time}
(or \textbf{Markov time}), 
if the event $\{\omega: \tau \leq t\}$ is determined 
by the behavior of the process
until time $t$ and independent of that after $t$.
For any stopping time $\tau$, 
$\{B(t): t \leq \tau \}$ and $\{B(t) : t > \tau \}$ 
are independent.
It is called the
\textbf{strong Markov property}.
A stochastic process which has the strong Markov property
and has a continuous path almost surely 
is generally called a 
\textbf{diffusion process}.

For each time $t \in [0, \infty)$, 
we write the smallest $\sigma$-field generated by
the BM up to time $t \geq 0$ as
$\sigma(B(s): 0 \leq s \leq t)$ and define
\begin{equation}
\cF_t := \sigma(B(s) : 0 \leq s \leq t), \quad t \geq 0.
\label{eqn:fil1}
\end{equation}
By definition, with respect to any event in $\cF_t$,
$B(s)$ is measurable at every $s \in [0, t]$.
Then we have a nondecreasing family
$\{\cF_t : t \geq 0\}$ of sub-$\sigma$-fields
of the original $\sigma$-field $\cF$ in the probability space
$(\Omega, \cF, \rP)$ such that
$\cF_s \subset \cF_t \subset \cF$ for
$0 \leq s < t < \infty$.
We call this family of $\sigma$-fields 
a \textbf{filtration}.

Since the probability density of increment 
in any time-interval $t-s>0$, $p(t-s, z |0)$, 
has mean zero,
BM satisfies the equality
\begin{equation}
\rE[B(t) | \cF_s]=B(s),
\quad 0 \leq s < t < \infty,
\quad \mbox{a.s.}
\label{eqn:martingale}
\end{equation}
That is, the mean is constant in time,
even though the variance 
increases in time as $\sigma^2 = t$.
Processes with such a property are called 
\textbf{martingales}.
We note that for $0 \leq s < t < \infty$
\begin{eqnarray}
\rE[B(t)^2|\cF_s]
&=& \rE[(B(t)-B(s))^2+2(B(t)-B(s))B(s)+B(s)^2|\cF_s]
\nonumber\\
&=& \rE[(B(t)-B(s))^2|\cF_s]
+2 \rE[(B(t)-B(s))B(s)|\cF_s]
+ \rE[B(s)^2|\cF_s].
\nonumber
\end{eqnarray}
By the property {\bf (BM3)} and the definition of $\cF_s$, 
\begin{eqnarray}
&& \rE[(B(t)-B(s))^2|\cF_s]=t-s,
\nonumber\\
&& \rE[(B(t)-B(s))B(s)|\cF_s]
=\rE[B(t)-B(s)|\cF_s] \, B(s)=0,
\nonumber\\
&& \rE[B(s)^2|\cF_s]=B(s)^2.
\nonumber
\end{eqnarray}
Then we have the equality
\begin{equation}
\rE[B(t)^2-t |\cF_s]=B(s)^2-s, \quad
0 \leq s < t < \infty, \quad \mbox{a.s.}
\label{eqn:martingale2}
\end{equation}
It means that $B(t)^2-t$ is a martingale.

For the transition probability density of BM (\ref{eqn:pt1}), 
it should be noted that $p(\cdot, y|x)=p(\cdot, x|y)$
for any $x, y \in \R$, 
and $u_{t}(x) := p(t, y|x)$
is a unique solution of the 
\textbf{heat equation} 
(\textbf{diffusion equation})
\begin{equation}
\frac{\partial}{\partial t} u_t(x)
=\frac{1}{2} \frac{\partial^2}{\partial x^2} u_t(x),
\quad x \in \R, \quad t \geq 0
\label{eqn:heat_eq1}
\end{equation}
with the initial condition $u_0(x)=\delta(x-y)$.
The solution of (\ref{eqn:heat_eq1})  with the initial condition
$u^f_0(x)=f(x), x \in \R$ is then given by
\begin{equation}
u^f_t(x)= \rE^x[f(B(t))]= \int_{-\infty}^{\infty} f(y) p(t, y|x) dy,
\label{eqn:heat_eq2}
\end{equation}
if $f$ is a measurable function
satisfying the condition
$\int_{-\infty}^{\infty} e^{-a x^2}|f(x)| dx < \infty$
for some $a>0$.
Since $p(t,y|x)$ plays as an integral kernel in (\ref{eqn:heat_eq2}),
it is also called the 
\textbf{heat kernel}.

For $0 \leq s < t < \infty$, $\xi \in \R$, consider 
$\rE[e^{\sqrt{-1} \xi (B(t)-B(s))}| \cF_s]$.
It is calculated as 
the \textbf{Fourier transform of the
transition probability density} $p(t-s, \cdot|0)$;
\begin{eqnarray}
\int_{-\infty}^{\infty} e^{\sqrt{-1} \xi z}
p(t-s, z|0) d z
&=& \int_{-\infty}^{\infty} e^{\sqrt{-1} \xi z}
\frac{e^{-z^2/2(t-s)}}{\sqrt{2 \pi (t-s)}} d z
\nonumber\\
&=& e^{-\xi^2 (t-s)/2}.
\label{eqn:charaBM0}
\end{eqnarray}
The obtained function of $\xi \in \R$,
\begin{equation}
\rE[e^{\sqrt{-1} \xi (B(t)-B(s))} | \cF_s]
=e^{-\xi^2(t-s)/2}, \quad 0 \leq s < t < \infty,
\label{eqn:charaBM1}
\end{equation}
is called the 
\textbf{characteristic function} of BM.

\subsection{It\^o's formula and Kolmogorov equation}
\label{sec:Ito}

Let $X(t), t \geq 0$ be a one-dimensional diffusion process
on the probability space $(\Omega_X, \cF_X, \rP_X)$,
where the expectation is written as $\rE_X$
and the filtration is given by the 
\textbf{natural filtration} of $X$; 
$(\cF_X)_t=\sigma(X(s): 0 \leq s \leq t), t \geq 0$.
For each time interval $[0, t], t > 0$, put $n \in \N$ and
let $\Delta_n=\Delta_n([0, t])$ be a subdivision of $[0, t]$ with 
$t_0:=0 < t_1 < \cdots < t_{n-1} < t_n := t$.
Then we define
$Q^{\Delta_n}(t) = \sum_{m=1}^n (X(t_m)-X(t_{m-1}))^2$.
If there is a process $Q(t), t \geq 0$ such that
\begin{equation}
\lim_{n \to \infty} \rP_X
[|Q^{\Delta_n}(t)-Q(t)| > \varepsilon]=0,
\quad ^{\forall}\varepsilon >0
\label{eqn:qv0}
\end{equation}
holds provided 
$\max_{1 \leq m \leq n} |t_m-t_{m-1}| \to 0$
as $n \to \infty$, 
then we call $Q(t)$, $t \geq 0$,  
the 
\textbf{quadratic variation} of $X(t), t \geq 0$
and express it by
$\langle X, X \rangle_t, t \geq 0$.

For BM, $B(t), t \geq 0$, set $n \in \N$, 
$t_0 :=0 < t_1 < \cdots < t_{n-1} < t_n := t$
and put
$
Q^{\Delta_n}_{\rm BM}(t)= \sum_{m=1}^n (B(t_m)-B(t_{m-1}))^2.
$
By the property {\bf (BM3)}, the mean is given by
$
\rE[Q^{\Delta_n}_{\rm BM}(t)]=\sum_{m=1}^n (t_m-t_{m-1})=t.
$
The variance of $Q^{\Delta_n}_{\rm BM}(t)$ 
\begin{eqnarray}
\sigma^{\Delta_n}_{\rm BM}(t)^2 &:=&
\rE[(Q^{\Delta_n}_{\rm BM}(t)-t)^2]
\nonumber\\
&=& \rE \left[ \left\{
\sum_{m=1}^n \Big\{ (B(t_m)-B(t_{m-1}))^2-(t_m-t_{m-1}) \Big\}
\right\}^2 \right]
\nonumber
\end{eqnarray}
is calculated as
\begin{eqnarray}
&& \sum_{m=1}^n \Big\{
\rE[(B(t_m)-B(t_{m-1}))^4]-2 (t_m-t_{m-1}) \rE[(B(t_m)-B(t_{m-1}))^2]
\nonumber\\
&& \qquad \qquad 
+(t_m-t_{m-1})^2 \Big\}
\nonumber\\
&=& \sum_{m=1}^{n} \Big\{
3(t_m-t_{m-1})^2-2(t_m-t_{m-1})^2+(t_m-t_{m-1})^2 \Big\}
\nonumber\\
&=& 2 \sum_{m=1}^n (t_m-t_{m-1})^2
\leq 2 t \max_{1 \leq m \leq n} |t_m-t_{m-1}|,
\label{eqn:sigma1}
\end{eqnarray}
where
independence of increments of BM mentioned in {\bf (BM3)} and  
\begin{equation}
\rE[(B(t)-B(s))^4]= 3 (t-s)^2, \quad 0 \leq s \leq t
\label{eqn:dB4}
\end{equation}
were used.
By 
\textbf{Chebyshev's inequality}, we have
$$
\rP[|Q^{\Delta_n}_{\rm BM}(t)-t| > \varepsilon] \leq \frac{\sigma^{\Delta_n}_{\rm BM}(t)^2}{\varepsilon^2},
\quad ^{\forall} \varepsilon >0.
$$
Provided that $\max_{1 \leq m \leq n}|t_m-t_{m-1}| \to 0$
as $n \to \infty$, 
(\ref{eqn:sigma1}) gives $\lim_{n \to \infty} \sigma^{\Delta_n}_{\rm BM}(t)=0$
and it proves
\begin{equation}
\langle B, B \rangle_t=t, \quad
\quad t \geq 0.
\label{eqn:Bqv1}
\end{equation}

For a stopping time $\tau$, we put
$X^{\tau}(t) := X(t \wedge \tau), t \geq 0$,
where $t \wedge \tau := \min\{t, \tau\}$. 
We define a diffusion process $X(t), t \geq 0$
as a \textbf{local martingale}, 
if there exists stopping times $\tau_n, n \in \N$ such that
(i) the sequence $\{\tau_n\}_{n \in \N}$ is nondecreasing
and $\lim_{n \to \infty} \tau_n=\infty$ a.s., and
(ii) for every $n$, the process 
$X^{\tau_n}(t), t \geq 0$
is a martingale. 
When $X(t)$ is a local martingale,
we can prove that a unique increasing continuous process
is given by $\langle X, X \rangle_t, t \geq 0$
such that $\langle X, X \rangle_0=0$ and
$X(t)^2-\langle X, X \rangle_t, t \geq 0$
provides a local martingale. 

Assume that $X(t)$ and $Y(t), t \geq 0$
are both local martingales.
Then $(X(t)+Y(t))^2-\langle X+Y, X+Y \rangle_t$
and $(X(t)-Y(t))^2-\langle X-Y, X-Y \rangle_t$, $t \geq 0$
are local martingales.
Therefore, their difference
$4 X(t) Y(t)-\{ \langle X+Y, X+Y \rangle_t-\langle X-Y, X-Y \rangle_t\}, t \geq 0$
is also a local martingale.
For any pair of local martingales $X(t)$ and $Y(t), t \geq 0$, we define
the 
\textbf{mutual quadratic variation} (\textbf{cross variation}) \textbf{process} as
\begin{equation}
\langle X, Y \rangle_t
:= \frac{1}{4}
\{ \langle X+Y, X+Y \rangle_t
- \langle X-Y, X-Y \rangle_t \},
\quad t \geq 0.
\label{eqn:qv1}
\end{equation}
We can prove that, 
if $B_i(t), t \geq 0, i=1,2, \dots, D$
are independent BMs, then
\begin{equation}
\langle B_i, B_j \rangle_t
=\delta_{ij} t,
\quad i, j =1, \dots, D, \quad t \geq 0.
\label{eqn:Bqv2}
\end{equation}

For a continuous process $A(t), t \geq 0$, here
we consider the quantity 
$$
S^{\Delta_n}(t)=\sum_{m=1}^{n}|A(t_m)-A(t_{m-1})|
$$
instead of $Q^{\Delta_n}(t)$,
where $\Delta_n, n \in \N$ is a subdivision
of the time interval $[0, t]$.
We can see that if $\Delta_n \subset \Delta_{n+1}$,
then $S^{\Delta_n}(t) \leq S^{\Delta_{n+1}}(t), ^{\forall} t \geq 0$.
Assume that $\max_{1 \leq m \leq n} |t_m-t_{m-1}| \to 0$
as $n \to \infty$.
Let $\lim_{n \to \infty} \sup_{\Delta_n} S^{\Delta_n}(t)=S(t) 
\leq \infty$
and call it the 
\textbf{variation} of $A$ on $[0, t]$.
If $S(t) < \infty$ for every $t$, then
the process $A(t), t \geq 0$ is of 
\textbf{finite variation}.
Let $S_{\rm BM}(t), t \geq 0$ be the variation of BM.
We have
$$
Q^{\Delta_n}_{\rm BM}(t) \leq 
\sup_{1 \leq m \leq n} |B(t_m)-B(t_{m-1})| S_{\rm BM}(t),
\quad t \geq 0, \quad n \in \N. 
$$
By the property {\bf (BM2)},
RHS becomes 0 as $n \to \infty$ if $S_{\rm BM}(t)$ is
finite.
On the other hand, we have proved 
$Q^{\Delta_n}_{\rm BM}(t) \to \langle B, B \rangle_t$ 
as $n \to \infty$ in probability and the fact (\ref{eqn:Bqv1}).
Hence $S_{\rm BM}(t)= \infty$ a.s. $^{\forall} t \geq 0$.

Let $N \in \N$ and $\{X_1(t), \dots, X_N(t)\}$, $t \geq 0$
be a set of diffusion processes.
Put $\X(t)=(X_1(t), \dots, X_N(t))$. 
Let $F$ be a real function 
of $(t, \x) \in [0, \infty) \times \R^N$,
which is bounded and has a bounded first-order derivative with respect to $t$ and
bounded first- and
second-order derivatives with respect to $x_j, 1 \leq j \leq N$, 
and we denote this by $F \in \rC^{1,2}_{\rm b}$.
We know that, 
in order to describe the statistics of a function 
of several random variables, we have to take into account
the `propagation of error'.
For the process $F(t, \X(t))$
that is defined as a function of $t$ as well as a functional of
processes $X_1(t), \dots, X_N(t)$, $t \geq 0$, 
\textbf{It\^o's formula} gives an equation which
governs the difference of $F(t, \X(\cdot))$ as
\begin{eqnarray}
d F(t, \X(t)) &=&
\sum_{i=1}^N \frac{\partial F}{\partial x_i}(t, \X(t)) d X_i(t)
+ \frac{\partial F}{\partial t}(t, \X(t))dt
\nonumber\\
&+& 
\frac{1}{2} \sum_{1 \leq i, j \leq N}
\frac{\partial^2 F}{\partial x_i \partial x_j}(t, \X(t)) 
d \langle X_i, X_j \rangle_t,
\quad t \geq 0.
\label{eqn:Ito1}
\end{eqnarray}
The first term gives a local martingale
and the second and third terms give processes of finite variations.
A continuous process $X$
given by the sum of a local martingale $M$
and a finite-variation process $A$, 
\[
X(t)=M(t)+A(t), \quad t \geq 0,
\]
is called a 
\textbf{semi-martingale}.
We will use the fact that if the continuous process
$F(t, \X(t)), \ t \geq 0$ is a local martingale, then
its finite-variation part should vanish, 
and vice versa.

Let $X(t), t \geq 0$ be a one-dimensional diffusion 
in $(\Omega_X, \cF_X, \rP_X)$,
which satisfies the following SDE
\begin{equation}
dX(t)=\sigma(X(t)) dB(t)+b(X(t))dt,
\quad X(0)=x.
\label{eqn:SDEb1}
\end{equation}
Here $B(t),t \geq 0$ is BM and the functions 
$\sigma, b: \R \mapsto \R$ satisfy 
the condition that
$^{\exists} K \geq 0$, s.t.
$|\sigma(x)-\sigma(y)| \leq K|x-y|$,
$|b(x)-b(y)| \leq K|x-y|, x , y \in \R$.
(This is called the 
\textbf{Lipschitz continuity}.)
Put
\begin{equation}
u(s,x)=\rE^x[f(X(T-s))],
\quad 0 \leq s < T < \infty,
\label{eqn:Kol1}
\end{equation}
with an $(\cF_X)_T$-measurable bounded function $f$.
By the Markov property, 
for $0 \leq s < t < T < \infty$,
\begin{eqnarray}
u(s,x) &=& \rE^x \Big[ \rE^{X(t-s)} [f(X(T-t))] \Big]
\nonumber\\
&=& \rE^x[u(t, X(t-s))].
\label{eqn:Kol2}
\end{eqnarray}
Assume that $u(s,x) \in \rC^{1,2}_{\rm b}$.
Then by It\^o's formula (\ref{eqn:Ito1}),
\begin{eqnarray}
u(t, X(t-s))-u(s,x)
&=& \int_{0}^{t-s} 
\left( \frac{\partial u}{\partial s}+ \sL u \right)
(s+r, X(r)) dr
\nonumber\\
&+&
\int_0^{t-s} \sigma(X(r))
\frac{\partial u}{\partial x}(s+r, X(r)) dB(r),
\label{eqn:K3}
\end{eqnarray}
where
\begin{equation}
(\sL f)(x)=\frac{1}{2} a(x) \frac{d^2}{dx^2} f(x)
+ b(x) \frac{d}{dx} f(x)
\label{eqn:K4}
\end{equation}
with
\begin{equation}
a(x)=\sigma(x)^2.
\label{eqn:K5}
\end{equation}
Since (\ref{eqn:Kol2}) holds, taking expectation of (\ref{eqn:K3}) gives
$$
\rE^{x} \left[
\int_{0}^{t-s} 
\left( \frac{\partial u}{\partial s}+ \sL u \right)
(s+r, X(r)) dr
\right] =0,
$$
where the expectation of the second term 
in RHS of (\ref{eqn:K3}) vanished
for it is a martingale.
Set $h=t-s$, divide the both side of the above equation
by $h$ and take the limit $h \to 0$.
Then we have
\begin{equation}
\frac{\partial u(s,x)}{\partial s}
+ \sL u(s,x)=0.
\label{eqn:K6}
\end{equation}
Since $\sL$ acts as a differential operator
with respect to the initial value $x$, 
it is called the
\textbf{backward Kolmogorov equation}.
The differential operator (\ref{eqn:K4})
is called the 
\textbf{generator} of the diffusion process.

\subsection{Martingale polynomials of BM}
\label{sec:poly_mar}

Let
\begin{equation}
G_{\alpha}(t, x)=e^{\alpha x-\alpha^2 t/2}.
\label{eqn:G1}
\end{equation}
For BM, we perform the transformation
with the parameter 
$\alpha \in \C := \{z=x+\sqrt{-1} y : x, y \in \R\}$, 
$B \mapsto \check{B}_{\alpha}$ defined by
\begin{equation}
\check{B}_{\alpha}(t)= G_{\alpha}(t, B(t)), 
\quad t \geq 0,
\label{eqn:Esscher1}
\end{equation}
which is called the
\textbf{Esscher transformation}.
Since
$$
\rE[e^{\alpha B(t)}]
=\int_{-\infty}^{\infty} e^{\alpha x} p(t, x|0) dx
=e^{\alpha^2 t/2}, \quad t \geq 0.
$$
we can see that
\begin{equation}
G_{\alpha}(t, B(t))= \frac{e^{\alpha B(t)}}{\rE[e^{\alpha B(t)}]},
\quad t \geq 0,
\label{eqn:Esscher1b}
\end{equation}
and hence
\begin{eqnarray}
\rE[G_{\alpha}(t, B(t)) | \cF_s ]
&=& \frac{\rE[e^{\alpha B(t)} | \cF_s]}{\rE[ e^{\alpha B(t)}]}
\nonumber\\
&=& \frac{\rE[e^{\alpha B(s)} e^{\alpha(B(t)-B(s))} | \cF_s]}
{\rE[e^{\alpha B(s)} e^{\alpha(B(t)-B(s))} ]},
\quad 0 < s < t.
\nonumber
\end{eqnarray}
By the definition of $\cF_s$ and independence of increment
of BM (the property {\bf (BM3)}),
the numerator is equal to $e^{\alpha B(s)} \rE[e^{\alpha(B(t)-B(s))}]$,
and the denominator is equal to $\rE[e^{\alpha B(s)}] \rE[e^{\alpha (B(t)-B(s))}]$.
Hence the above equals 
$e^{\alpha B(s)}/\rE[e^{\alpha B(s)}]
= G_{\alpha}(s, B(s))$.
Therefore, $G_{\alpha}(t, B(t))$ is a martingale: 
\begin{equation}
\rE[G_{\alpha}(t, B(t)) | \cF_s ]=G_{\alpha}(s, B(s)), \quad
0 \leq s \leq t.
\label{eqn:G1b}
\end{equation}

The function (\ref{eqn:G1}) is expanded as
\begin{equation}
G_{\alpha}(t, x) = \sum_{n=0}^{\infty} m_n(t, x) \frac{\alpha^n}{n!}
\label{eqn:G2}
\end{equation}
with
\begin{equation}
m_n(t, x)=\left( \frac{t}{2} \right)^{n/2}
H_n \left( \frac{x}{\sqrt{2t}} \right),
\quad n \in \N_0 := \{0,1,2, \dots\}. 
\label{eqn:mn1}
\end{equation}
Here $\{H_n(x) \}_{n \in \N_0}$ are the 
\textbf{Hermite polynomials} of degrees
$n \in \N_0$,
\begin{eqnarray}
\label{eqn:Hn1}
H_n(x) 
:=(-1)^n e^{x^2} \frac{d^n e^{-x^2}}{d x^n}
=
\sum_{k=0}^{[n/2]} (-1)^k \frac{n!}{k! (n-2k)!} (2x)^{n-2k},
\end{eqnarray}
where for $a \geq 0$, $[a]$ denotes the largest integer
which is not larger than $a$. 
We can show the following contour integral representations
of the Hermite polynomials,
\begin{equation}
H_n(z)= \frac{n!}{2 \pi \sqrt{-1}} \oint_{C(\delta_0)} d \eta \,
\frac{e^{2 \eta z-\eta^2}}{\eta^{n+1}}, \quad n \in \N_0,
\label{eqn:Hncontour}
\end{equation}
where $C(\delta_0)$ is a closed contour on the complex plane $\C$
encircling the origin 0 once in the positive direction.

\begin{lem}
\label{thm:mn}
The functions $\{m_n(t, x) \}_{n \in \N_0}$ satisfy the following.\\
(i) They are monic polynomials of degrees $n \in \N_0$
with time-dependent coefficients: 
$$
m_n(t, x)=x^n+\sum_{k=0}^{n-1} c_n^{(k)}(t) x^k, \quad t \geq 0.
$$
(ii) For $0 \leq k \leq n-1$, $c_n^{(k)}(0)=0$. That is,
$$
m_n(0, x)=x^n, \quad n \in \N_0.
$$
(iii) If we set $x=B(t)$, they provide martingales: 
\begin{equation}
\rE[m_n(t, B(t)) |\cF_s]=m_n(s, B(s)), \quad 0 \leq s \leq t,
\quad n \in \N_0.
\label{eqn:mn_martingale}
\end{equation}
\end{lem}

We call $\{m_n(t, x)\}_{n \in \N_0}$
the \textbf{fundamental martingale polynomials}
associated with BM \cite{Kat14}. 
For $n=2$, (\ref{eqn:Hn1}) gives
$H_2(x)=4x^2-2$, and then
$m_2(t,x)=x^2-t$ by (\ref{eqn:mn1}).
We already proved in (\ref{eqn:martingale2})
that $m_2(t, B(t))=B(t)^2-t$ is a martingale.

The Fourier transformation of $G_{\alpha}(t, x)$ with respect to
the parameter $\alpha \in \R$ is calculated as
\begin{equation}
\widehat{G}_{w}(t, x) :=
\int_{-\infty}^{\infty} \frac{e^{-\sqrt{-1} \alpha w}}{2 \pi}
G_{\alpha}(t, x) d \alpha
= \frac{e^{-(\sqrt{-1} x+ w)^2/2t}}{\sqrt{2 \pi t}}.
\label{eqn:Ghat1}
\end{equation}
Owing to the factor $e^{-w^2/2t}$ in $\widehat{G}_{w}(t, x)$,
the following calculations are justified,
\begin{eqnarray}
G_{\alpha}(t, x) &=& \int_{-\infty}^{\infty} e^{\sqrt{-1} \alpha w} 
\widehat{G}_{w}(t, x) d w
=\sum_{n=0}^{\infty} \frac{\alpha^n}{n!} \int_{-\infty}^{\infty} 
(\sqrt{-1} w)^n \widehat{G}_w(t,x) dw,
\label{eqn:Ghat1b}
\end{eqnarray}
which proves the integral representation of 
$m_n(t, x)$,
\begin{eqnarray}
m_n(t, x) &=& \int_{-\infty}^{\infty} (\sqrt{-1} w)^n
\widehat{G}_w(t, x) dw
\nonumber\\
&=& \int_{-\infty}^{\infty} (\sqrt{-1} w)^n 
\frac{e^{-(\sqrt{-1} x+w)^2/2t}}{\sqrt{2 \pi t}} dw,
\quad t \geq 0, \quad n \in \N_0.
\label{eqn:mn2}
\end{eqnarray}

We define this type of integral transformation of a function $f$ as
\begin{equation}
\cI[f(W)|(t, x)]
= \int_{-\infty}^{\infty} f(\sqrt{-1} w) \widehat{G}_w(t,x) dw.
\label{eqn:I1}
\end{equation}
Then the above results are written as
\begin{equation}
m_n(t,x)= \cI[W^n|(t,x)],  \quad t \geq 0, \quad n \in \N_0.
\label{eqn:mn3}
\end{equation}


\subsection{Bessel processes (BES$_D$)}
\label{sec:BES}
\subsubsection{Radial part of $D$-dimensional BM}
\label{sec:radialBES}

Let $D \in \N$ denote the spatial dimension.
For $D \geq 2$, the $D$-dimensional BM in $\R^D$ starting
from the position $\x=(x_1, \dots, x_D) \in \R^D$
is defined by the $D$-dimensional
vector-valued diffusion process,
\begin{equation}
\B^{\x}(t)=(B^{x_1}_1(t), B^{x_2}_2(t), \dots,
B^{x_D}_D(t)),
\quad t \geq 0,
\label{eqn:DBM}
\end{equation}
where 
\[
B^{x_i}_i(t) := x_i + B_i(t), \quad t \geq 0, \quad i=1, \dots, D
\]
are independent one-dimensional BMs
starting from $x_i \in \R$.

The $D$-dimensional 
\textbf{Bessel process} 
is defined as the absolute value
({\it i.e.}, the radial coordinate) of the
$D$-dimensional Brownian motion,
\begin{eqnarray}
R^{x}(t) &:=& |\B^{\x}(t)|
\nonumber\\
&=& \sqrt{ B^{x_1}_1(t)^2+ \cdots
+ B^{x_D}_D(t)^2},
\quad t \geq 0,
\label{eqn:BES1}
\end{eqnarray}
where the initial value is given by
$R^{x}(0)=x=|\x|=\sqrt{x_1^2+ \dots+ x_D^2} \geq 0$.
By definition $R^{x}(t)$ is nonnegative,
$R^{x}(t) \in \R_+ \cup \{0\}$,
where $\R_{+} := \{x \in \R: x > 0\}$.
We will abbreviate the $D$-dimensional Bessel process
to BES$_{D}$.

By this definition, $R^{x}(t)$ is a functional of
$D$-tuples of diffusion processes 
$\{B^{x_i}_i(t) \}_{i=1}^D$, $t \geq 0$. 
Now we apply It\^o's formula (\ref{eqn:Ito1}) 
to BES$_{D}$ (\ref{eqn:BES1}).
Assume $x=|\x|>0$.
In this case
\begin{equation}
R^x(t)=F(\B^{\x}(t)), \quad t \geq 0
\quad \mbox{with} \quad
F(\y)=\sqrt{\sum_{i=1}^D y_i^2}.
\label{eqn:BES2}
\end{equation}
We see that 
$$
\frac{\partial F}{\partial t}=0, \quad
\frac{\partial F}{\partial y_i}=\frac{y_i}{F}, \quad
\frac{\partial^2 F}{\partial y_i \partial y_j}
=\frac{\delta_{ij}}{F}-\frac{y_i y_j}{F^3},
\quad 1 \leq i, j \leq D.
$$
From (\ref{eqn:Bqv2}), we have
$d \langle B_i^{x_i}, B_j^{x_j} \rangle_t=\langle d B_i^{x_i}, d B_j^{x_j} \rangle_t = \delta_{ij} dt$, $1 \leq i, j \leq D, t \geq 0$.
Then, the second term of (\ref{eqn:Ito1}) for BES$_{D}$
becomes
$$
\frac{1}{2} \sum_{1 \leq i, j \leq D}
\left\{ \frac{\delta_{ij}}{F(\B^{\x}(t))}-\frac{B_i^{x_i}(t) B_j^{x_j}(t) }{F(\B^{\x}(t))^3} \right\}\delta_{ij} dt
= \frac{D-1}{2} \frac{1}{F(\B^{\x}(t))} dt
=\frac{D-1}{2} \frac{dt}{R^x(t)}.
$$
On the other hand, the first term of (\ref{eqn:Ito1}) for BES$_{D}$ is
\begin{equation}
\frac{1}{R^x(t)} \sum_{i=1}^D B_i^{x_i}(t) dB_i^{x_i}(t).
\label{eqn:BESb1}
\end{equation}
It seems to be complicated, but the quadratic variation is calculated as
\begin{eqnarray}
\left\langle \frac{1}{R^x} \sum_{i=1}^D B_i^{x_i} dB_i^{x_i}, \frac{1}{R^x} \sum_{j=1}^D B_j^{x_j} dB_j^{x_j}  \right\rangle_t
&=& \frac{1}{R^x(t)^2}
\sum_{i=1}^D \sum_{j=1}^D  
B_i^{x_i}(t) B_j^{x_j}(t) \langle dB_i^{x_i},  dB_j^{x_j} \rangle_t
\nonumber\\
&=& \frac{1}{R^x(t)^2} \sum_{i=1}^D \sum_{j=1}^D B_i^{x_i}(t) B_j^{x_j}(t) \delta_{ij} dt
= dt,
\nonumber
\end{eqnarray}
where the independence of BMs (\ref{eqn:Bqv2})
and the definition, $R^x(t)^2=\sum_{i=1}^D B_i^{x_i}(t)^2$, 
have been used.
That is, (\ref{eqn:BESb1}) is equivalent in probability law with
an infinitesimal increment of a diffusion process with
quadratic variation $dt$.
Then, by introducing a BM, $B^x(t), t \geq 0$,
which is different from $B_i^{x_i}(t), t \geq 0, x_i \in \R, i=1,2, \dots, D$
and is started at $x=|\x|>0$, 
(\ref{eqn:BESb1}) is identified with $dB^{x}(t), t \geq 0$.
Hence we have obtained the following equation for BES$_{D}$,
\begin{equation}
dR^x(t)= dB^x(t)
+\frac{D-1}{2} \frac{dt}{R^x(t)}, \quad x > 0, \quad 0 \leq t < T^x,
\label{eqn:BESeq1}
\end{equation}
where $T^x=\inf \{t > 0: R^x(t)=0\}$.

The first term of RHS, $dB^x(t)$, denotes the infinitesimal
increment of BM starting
from $x >0$ at time $t=0$.
This martingale term gives randomness to the motion.
On the other hand, if $D >1$, 
for $dt > 0$, the second term in RHS of (\ref{eqn:BESeq1})
is positive definite.
It means that there is a drift to increase the value of $R^{x}(t)$.
This drift term is increasing in $D$ and decreasing in $R^{x}(t)$.
Since as $R^{x}(t) \downarrow 0$, the drift term
$\uparrow \infty$, it seems that 
a `repulsive force' is acting to the $D$-dimensional BM, 
$\B^{\x}(t), |\x| > 0$ to keep the distance
from the origin be positive, $R^{x}(t)=|\B^{\x}(t)| > 0$
and avoid a collision of the Brownian particle
with the origin.
A differential equation such as (\ref{eqn:BESeq1}), which involves
a random fluctuation term and a drift term is called a 
\textbf{stochastic differential equation} (SDE).

The following equivalence in probability law is established
for arbitrary $x>0$,
\begin{equation}
\left( \frac{1}{x} R^x(x^2 t) \right)_{t \geq 0} \law=
(R^1(t))_{t \geq 0}.
\label{eqn:BESscaling}
\end{equation}
It inherits (\ref{eqn:BMscaling}) and is called
the \textbf{scaling property of the Bessel process}.

The SDE for BES$_{D}$ is given by the equation (\ref{eqn:SDEb1})
with $\sigma(x) := 1$ and $b(x)=(D-1)/(2x)$.
Then the generator of BES$_{D}$ is obtained as
\begin{equation}
\sL^{D}=\frac{1}{2} \frac{\partial^2}{\partial x^2}
+\frac{D-1}{2x} \frac{\partial}{\partial x}.
\label{eqn:K7}
\end{equation}

Let $p_{D}(t-s, y|x)$ be the transition probability density of
BES$_{D}$ from $x$ at times $s \geq 0$
to $y$ at times $t \geq s$.
For any $t \geq s, y \in \R_{+}$,
$p_{D}(t-s, y|x), x>0$ solves 
(\ref{eqn:K6}) with (\ref{eqn:K7})
under the condition 
$\lim_{s \uparrow t}p_{D}(t-s,y|x)
=\delta(x-y)$.
In other words, 
$p_{D}(t,y|x)$ solves
\begin{equation}
\frac{\partial}{\partial t} p_{D}(t,y|x)
=\sL_{D} p_{D}(t,y|x)
\label{eqn:BES3}
\end{equation}
under the initial condition $p_{D}(0,y|x)=\delta(x-y)$,
which is called the backward Kolmogorov equation
for BES$_{D}$.

Let $I_{\nu}(z)$ be the 
\textbf{modified Bessel function}
of the first kind defined by 
\begin{equation}
I_{\nu}(z) = \sum_{n=0}^{\infty} 
\frac{1}{\Gamma(n+1) \Gamma(n+1+\nu)}
\left( \frac{z}{2} \right)^{2n+\nu}
\label{eqn:I}
\end{equation}
with the gamma function 
\begin{equation}
\Gamma(z)=\int_0^{\infty} e^{-u} u^{z-1} du, 
\quad \Re z > 0.
\label{eqn:Gamma_def}
\end{equation}
The function $I_{\nu}(z)$ solves the 
\textbf{Bessel differential equation}
\begin{equation}
\frac{d^2 w}{dz^2}
+\frac{1}{z} \frac{dw}{dz} 
- \left( 1+ \frac{\nu^2}{z^2} \right) w=0.
\label{eqn:Besseleq}
\end{equation}
Then we can show that 
\begin{equation}
p_{D}(t,y|x)
=  \left\{ \begin{array}{ll}
\displaystyle{
\frac{1}{t} \frac{y^{\nu+1}}{x^{\nu}}
e^{-(x^2+y^2)/2t}
I_{\nu} \left( \frac{xy}{t} \right)},
& \quad t>0, \, x >0, \, y \geq 0, \cr
\displaystyle{
\frac{y^{2\nu+1}}{2^{\nu} t^{\nu+1} \Gamma(\nu+1)} e^{-y^2/2t}},
& \quad t>0, \, x=0, \, y \geq 0 \cr
& \cr
\delta(y-x),
& \quad t=0, \, x, y \geq 0,
\end{array} \right.
\label{eqn:pD}
\end{equation}
where the index $\nu$ is specified by the dimension $D$ as
\begin{equation}
\nu=\frac{D-2}{2}
\quad \Longleftrightarrow \quad
D=2(\nu+1).
\label{eqn:nuD}
\end{equation}
This fact that $p_{D}(t,y|x)$ is expressed by
using $I_{\nu}(z)$ gives the reason why
the process $R^{x}(t)$ is called the Bessel process.

\subsubsection{BES$_{3}$ and absorbing BM}
\label{sec:BES3}

When $D=3$, $\nu=1/2$ by (\ref{eqn:nuD}),
and we can use the equality 
$I_{1/2}(z)=\sqrt{2/\pi z} \sinh z$
$=(e^{z}-e^{-z})/\sqrt{2 \pi z}$.
Then (\ref{eqn:pD}) gives
\begin{equation}
p_{3}(t,y|x)=\frac{y}{x} \Big\{
p(t,y|x)-p(t,y|-x) \Big\}
\label{eqn:p3}
\end{equation}
for $t >0, x >0, y \geq 0$, where $p(t,y|x)$ is 
the transition probability density of BM started at $x$ given by (\ref{eqn:pt1}). 
If we put 
\begin{equation}
q_{\rm abs}(t,y|x)=p(t,y|x)-p(t,y|-x), 
\label{eqn:p_abs}
\end{equation}
we see that
$q_{\rm abs}(t,0|x)=0$ for any $x > 0$,
since the transition probability density of BM, $p(t, y|x)$, 
is an even function of $y-x$.

\begin{figure}[t]
\begin{center}
\includegraphics[scale=0.6]{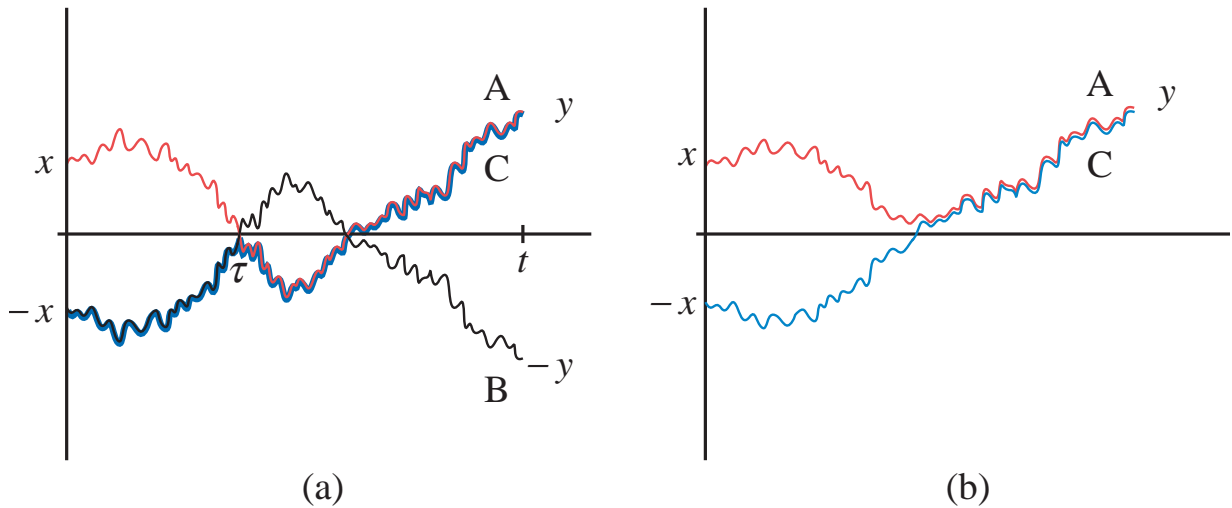}
\caption{
One realization of a Brownian path
from $x>0$ to $y>0$ is represented by a red curve and denoted by path A,
which visits the nonpositive region $\R_{-} \cup \{0\}$, where 
$\R_{-} := \{x \in \R: x < 0 \}$.
The first time the path A hits the origin is
denoted by $\tau$. 
Path B, which is represented by a black curve, is a mirror image of path A with respect
to the origin $x=0$, which is running from $-x<0$ to $-y < 0$. 
The path C (blue curve) is then defined as 
the concatenation of the part of path B from time zero
up to time $\tau$ 
and the part of path A after $\tau$
such that it runs from $-x < 0$ to $y>0$.
By the reflection principle of BM (\ref{eqn:reflection_principle})
applied at time $\tau$, we see that 
there exists a bijection between path A and path C,
which have the same probability weight as Brownian paths.
Since the Brownian path A contributes
to $p(t, y|x)$ and the Brownian path C contributes to $p(t, y|-x)$,
such a path from $x>0$ to $y>0$ that visits $\R_-$ is cancelled
in $q_{\rm abs}(t, y|x)$ given by (\ref{eqn:p_abs}). 
In Fig.\ref{fig:absBM} (b), a path from $x>0$ to $y>0$ 
which stays in the positive region $\R_+$
is considered (path A).
In this case there is no such a path
that perfectly cancels the contribution of path A to (\ref{eqn:p_abs}) 
as path B in the case (a).
In summary, $q_{\rm abs}(t, y|x)=p(t, y|x)-p(t, y|-x)$
gives the total weight of Brownian paths
which do not hit the origin.
}
\end{center}
\label{fig:absBM}       
\end{figure}

We consider the situation where an absorbing wall is put at the origin and,
if the Brownian particle starting from $x>0$ arrives at the origin,
it is absorbed there and the motion is stopped.
Such a process is called
the 
\textbf{absorbing Brownian motion in $\R_+$}.
Its transition probability density is given by
$q_{\rm abs}$.

By absorption, the total mass of paths from $x >0$ to $y>0$
is then reduced, if we compare the
original BM and the absorbing Brownian motion in $\R_+$.
The factor $y/x$ appearing in the transition probability density
(\ref{eqn:p3}) of BES$_{3}$ is for 
renormalization so that
\begin{equation}
\int_{\R_+} p_{3}(t,y|x)dy=1, \quad ^{\forall} t >0, 
\quad ^{\forall} x >0. 
\label{eqn:BES3_normalization}
\end{equation}
We regard this renormalization procedure from $q_{\rm abs}$
to $p_{3}$ as a transformation.
Since 
\begin{equation}
h_1(x)=x
\label{eqn:h1}
\end{equation}
is a one-dimensional harmonic function 
in a rather trivial sense
\[
\Delta^{(1)} h_1(x) := \frac{d^2 x}{d x^2}=0,
\]
we say that the BES$_{3}$ is
a \textbf{harmonic transformation}
($h$-transformation) of the one-dimensional absorbing BM
in the sense of Doob.
This implies the following equivalence.
\begin{eqnarray}
\mbox{BES$_{3}$} \,
&\Longleftrightarrow& \,
\mbox{one-dimensional Brownian motion
conditioned to stay positive}
\nonumber\\
&\Longleftrightarrow& \,
\mbox{$h$-transform of absorbing BM in $\R_+$}
\nonumber
\end{eqnarray}
Let $\E_{{\rm BES}_{3}}^x$ denote the expectation
with respect to BES$_{3}$, $(R(t))_{t \geq 0}$, started at $x \in \R_+$.
For an independent BM, $(B(t))_{t \geq 0}$ started at the same point $x \in \R_+$,
let $\tau=\inf \{t > 0: B(t)=0\}$.
Then the above equivalence is written as follows;
for any $\cF_t$-measurable bounded function
$F$, $t \geq 0$,
\begin{align}
\E_{{\rm BES}_{3}}^x[F(R(t))]
&=\rE^x \left[
F(B(t)) \1_{(\tau >t)} \frac{B(t)}{x} \right]
\nonumber\\
&=\rE^x \left[
F(B(t)) \1_{(\tau >t)} \frac{h_1(B(t)}{h_1(x)} \right],
\quad t \geq 0,
\label{eqn:BES3_b1}
\end{align}
where $\rE^x$ is an expectation with respect to
BM started at $x \in \R_+$.
If $F$ is an even function; $F(-x)=F(x)$, then
the above gives
\begin{equation}
\E_{{\rm BES}_{3}}^x[F(R(t))]
=\rE^x \left[
F(B(t)) \frac{h_1(B(t))}{h_1(x)} \right],
\quad t \geq 0,
\label{eqn:BES3_b2}
\end{equation}
since by the reflection principle of BM
(\ref{eqn:reflection_principle}),
all contribution from paths $\{\omega: \tau \leq t\}$
should be canceled out. 

Here we emphasize the obvious fact that $p_{3}(t, 0|x)=0,
^{\forall} x >0$.
It implies that BES$_{3}$ does not
visit the origin.
When $D=3$, the outward drift is strong enough to avoid
any visit to the origin.
Moreover, we can prove that for any $x>0$, $R^{x}(t) \to \infty$
as $t \to \infty$ with probability 1
and we say the process is 
\textbf{transient}
(see Theorem \ref{thm:BESD1} (ii) below).

\subsubsection{Critical dimension $D_{\rm c}=2$}
\label{sec:Dc}

Originally, the Bessel process was defined by (\ref{eqn:BES1})
for $D \in \N$.
We find that, however, the modified Bessel function 
(\ref{eqn:I}) is an analytic function of 
$\nu$ for all values of $\nu$.
So we will be able to define the Bessel process
for any value of dimension $D \geq 1$
as a diffusion process in $\R_+$
such that the transition probability density function
is given by (\ref{eqn:pD}), where the index 
$\nu \geq -1/2$ is determined by (\ref{eqn:nuD})
for each value of $D \geq 1$.

For BES$_{D}$ starting from $x>0$, denote
its first visiting time at the origin by
\begin{equation}
T^x=\inf \{t > 0: R^x(t)=0 \}.
\label{eqn:Tx}
\end{equation}
The following theorem is proved 
(see, for instance, \cite{Kat16_Springer}). 

\begin{thm}
\label{thm:BESD1}
\begin{description}
\item{\rm (i)} \quad
$D \geq 2 \quad \Longrightarrow \quad
T^{x}=\infty, ^{\forall}x >0$, with probability 1.
\item{\rm (ii)} \quad
$D > 2 \quad \Longrightarrow \quad
\displaystyle{\lim_{t \to \infty}} R^{x}(t)=\infty,
\, ^{\forall}x > 0$, with probability 1,
{\it i.e.} the process is transient.
\item{\rm (iii)} \quad
$D=2 \quad \Longrightarrow \quad
\displaystyle{\inf_{t >0}} \, R^{x}(t)=0, \, 
^{\forall}x >0$, with probability 1D \\
That is, BES$_{(2)}$ starting from $x>0$ does not
visit the origin, but it can visit any neighborhood of
the origin.
\item{\rm (iv)} \quad
$1 \leq D < 2 \quad \Longrightarrow \quad
T^x < \infty, \, ^{\forall} x > 0$, with probability 1,
{\it i.e.} the process is recurrent.
\end{description}
\end{thm}
\subsubsection{Bessel flow and 
another critical dimension $\overline{D}_{\rm c}=3/2$}
\label{sec:BES_Flow}

Theorem \ref{thm:BESD1} states that we have a
\textbf{critical dimension},
$$
 D_{\rm c}=2,
$$
for competition between the two effects acting 
the Bessel process,
the `random force' (the martingale term)
and the `entropy force' (the outward drift term)
in (\ref{eqn:BESeq1}): 
when $D > D_{\rm c}$, the latter dominates the former
and the process becomes transient, and when 
$D < D_{\rm c}$, the former is relevant and recurrence
to the origin of the process is realized frequently.

Here we show that there is another critical dimension \cite{Law05},
$$
\overline{D}_{\rm c}=\frac{3}{2}.
$$
In order to characterize the transition at $\overline{D}_{\rm c}$,
we have to investigate the dependence of the behavior
of $R^{x}(t)$ on its initial value, $x >0$.
We call the one-parameter family 
$\{R^{x}(t) : t \geq 0 \}_{x >0}$ the 
\textbf{Bessel flow} 
for each fixed $D > 1$.

For $0 < x < y$, we trace the motions of two BES$_{D}$'s
starting from $x$ and $y$ by solving
(\ref{eqn:BESeq1}) using the {\it common} BM, $B(t), t \geq 0$,
\begin{eqnarray}
&& R^{x}(t)= x+ B(t)+\frac{D-1}{2}
\int_0^{t} \frac{ds}{R^{x}(s)},
\nonumber\\
&& R^{y}(t)= y+ B(t)+\frac{D-1}{2}
\int_0^{t} \frac{ds}{R^{y}(s)}, \quad 0 \leq t < T^x.
\label{eqn:Bessel_flow}
\end{eqnarray}
We will see that 
\begin{eqnarray}
 x<y \quad &\Longrightarrow& \quad
R^x(t) < R^{y}(t), \quad 0 \leq t < T^{x}
\quad \mbox{with probability 1}
\nonumber\\
&\Longrightarrow& \quad
T^x \leq T^{y} \quad \mbox{with probability 1}.
\nonumber
\end{eqnarray}

The interesting fact is that in the 
intermediate fractional dimensions,
$\overline{D}_{\rm c} < D < D_{\rm c}$,
it is possible to have a situation where
$T^{x}=T^{y}$ even for $x < y$.

\begin{figure}[t]
\begin{center}
\includegraphics[scale=.50]{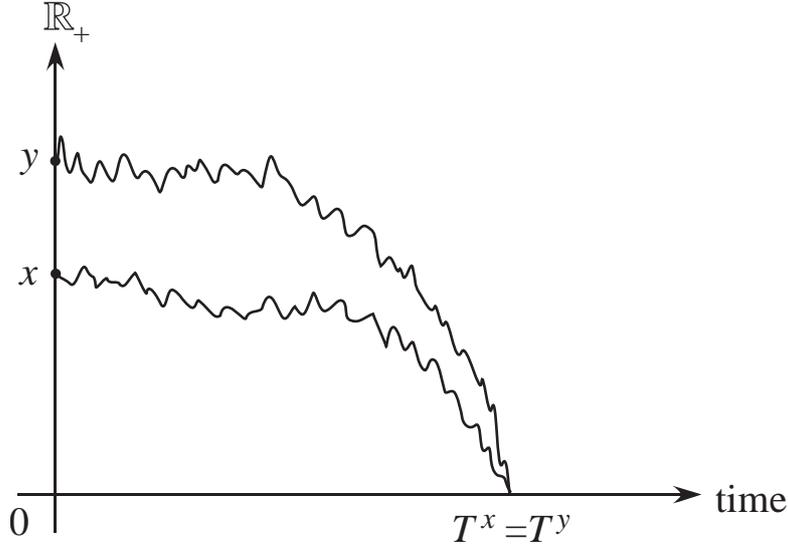}
\caption{In the intermediate fractional dimensions,
$3/2 < D < 2$, there is a positive probability
that two Bessel processes starting from different
initial positions, $0<x<y< \infty$, return
to the origin simultaneously, $T^x=T^y$.
}
\end{center}
\label{fig:BESflow}       
\end{figure}

\begin{thm}
\label{thm:BESD2}
For $0 < x < y < \infty$, 
\begin{description}
\item{\rm (i)} \quad
$1 \leq D \leq 3/2 \quad \Longrightarrow \quad
T^{x} < T^{y} \quad \mbox{with probability 1}$.
\item{\rm (ii)} \quad
$3/2 < D < 2 \quad \Longrightarrow \quad
\rP[T^x=T^y] > 0$.
\end{description}
\end{thm}
\vskip 0.5cm

\subsection{Hermitian-matrix-valued BM and 
the Dyson model with parameter $\beta$ (DYS$_{\beta}$)}
\label{sec:Dyson}
\subsubsection{Multivariate extensions of Bessel processes}
\label{sec:multi}

Here we consider the stochastic motion of two particles
$(X_1(t), X_2(t))$ in one dimension $\R$
satisfying the following SDEs,
\begin{eqnarray}
&& dX_1(t) = dB_1(t)
+\frac{\beta}{2} \frac{dt}{X_1(t)-X_2(t)}, 
\nonumber\\
&& dX_2(t) = dB_2(t)
+\frac{\beta}{2} \frac{dt}{X_2(t)-X_1(t)}, 
\label{eqn:Dyson1}
\end{eqnarray}
with the initial condition $x_1=X_1(0) <  x_2=X_2(0)$
for $0 \leq t < \inf\{t>0: X_1(t)=X_2(t)\}$, 
where $B_1(t)$ and $B_2(t)$, $t \geq 0$ are independent BMs 
and $\beta >0$ is the `coupling constant' of the two particles.
The second terms in (\ref{eqn:Dyson1})
represent the repulsive force acting between two 
particles, which is proportional to the inverse
of the distance between them, $X_2(t)-X_1(t)$.
Since it is a central force ({\it i.e.}, 
depending only on distance, and thus symmetric
for two particles), the `center of mass'
$X_{\rm c}(t) := (X_2(t)+X_1(t))/2$
is just a time change of BM; 
we can calculate the quadratic variation as
$d \langle X_{\rm c}, X_{\rm c} \rangle_t
=\langle dX_{\rm c}, dX_{\rm c} \rangle_t=
\langle (dB_1+dB_2)/2, (dB_1+dB_2)/2 \rangle_t=dt/2$,
since $d\langle B_1, B_1 \rangle_t=\langle dB_1, dB_1 \rangle_t= dt$,
$d \langle B_2, B_2 \rangle_t= \langle dB_2, dB_2 \rangle_t=dt$ and
$d \langle B_1, B_2 \rangle_t=\langle dB_1, dB_2 \rangle_t=0$.
Then
$$
(X_{\rm c}(t))_{t \geq 0}
\law= \left( \frac{1}{\sqrt{2}} B(t) \right)_{t \geq 0}
\law= (B(t/2))_{t \geq 0},
$$
where $B(t)$ is a BM independent from $B_1(t)$ and $B_2(t)$. 
On the other hand, if we define the relative
coordinate by
$X_{\rm r}(t) := (X_2(t)-X_1(t))/\sqrt{2}$, 
it satisfies the SDE
\begin{equation}
dX_{\rm r}(t)=d \widetilde{B}(t)
+\frac{\beta}{2} \frac{dt}{X_{\rm r}(t)}
\label{eqn:dXr}
\end{equation}
for $0 < t < \inf\{t > 0: X_{\rm r}(t)=0\}$, 
where $\widetilde{B}(t), \ t \geq 0$ is a BM independent from
$B_1(t), B_2(t), B(t)$, $t \geq 0$.
It is nothing but the SDE for BES$_{D}$ with
$D=\beta+1$: 
\begin{equation}
(X_{\rm r}(t))_{t \geq 0}
\law= (R(t))_{t \geq 0}, \quad D=\beta+1.
\label{eqn:DbetaA}
\end{equation}

Now we consider the following $N$-particle systems
of \textbf{interacting Brownian motions} 
in $\R$ as a solution 
$\X(t)=(X_1(t), X_2(t), \dots, X_N(t))$
of the following system of SDEs: 
with $\beta >0$ and the condition
$x_1 < x_2 < \cdots < x_N$
for initial positions $x_i=X_i(0), i=1, \dots, N$, 
\begin{equation}
dX_i(t)=dB_i(t)+ \frac{\beta}{2}
\sum_{\substack{1 \leq j \leq N, \cr j \not=i}}
\frac{dt}{X_i(t)-X_j(t)},
\quad t \in [0, T^{\x}),
\quad i=1, \dots, N,
\label{eqn:Dyson3}
\end{equation}
where $\{B_i(t) \}_{i=1}^N$, $t \geq 0$ are independent BMs
and 
\begin{eqnarray}
T^{\x}_{ij} &=& \inf \{ t > 0 :
X_i(t)=X_j(t)\}, \quad 1 \leq i < j \leq N,
\nonumber\\
T^{\x} &=& \min_{1 \leq i < j \leq N} T^{\x}_{ij}.
\nonumber
\end{eqnarray}
It is usually called 
\textbf{Dyson's Brownian motion model with parameter} $\beta$.
But in this lecture we simply call it the
\textbf{Dyson model with parameter} $\beta$
and denoted by DYS$_{\beta}$.

As shown above, in the case where $N=2$, DYS$_{\beta}$
is a composition of a BM (the center of mass) 
and a BES$_{(\beta+1)}$ (the relative coordinate).
In this sense, DYS$_{\beta}$ can be regarded
as a multivariate (multidimensional) extension
of BES$_{(\beta+1)}, \ \beta>0$. 

We can prove that,
for any $\x \in \R^N$ with $x_1 < x_2 < \cdots < x_N$, 
$T^{\x} < \infty$ if $\beta <1$,
and $T^{\x} = \infty$ if $\beta \geq 1$ \cite{RS93,GM13}.
The critical value $\beta_{\rm c}=1$ seems to
correspond to $D_{\rm c}=\beta_{\rm c}+1=2$ of
BES$_{D}$.
Together with (\ref{eqn:DbetaA}), this observation suggests the 
following relation between parameters, 
\begin{equation}
D=\beta+1 \quad \iff \quad
\beta=D-1.
\label{eqn:DbetaB}
\end{equation}

In particular, we study the special case of
Dyson's BM model with parameter $\beta=2$, DYS$_{2}$
As shown above, the case where $\beta=2$ 
corresponds to a BES$_{D}$ with $D=3$.
In Section \ref{sec:BES}, we have shown that
BES$_{3}$ has two aspects: 
{\bf [Aspect 1]} as a radial coordinate of three-dimensional
Brownian motion, which was used to define the Bessel process
in Section \ref{sec:radialBES}, 
and
{\bf [Aspect 2]} as a one-dimensional Brownian motion conditioned
to stay positive as explained in Section \ref{sec:BES3}.
We show that the Dyson model inherits these two aspects
from BES$_{3}$ \cite{Kat16_Springer}.

\subsubsection{DYS$_2$ realized as 
eigenvalue process of Hermitian-matrix-valued BM}
\label{sec:eigenvalue}

Dyson \cite{Dys62} introduced the processes (\ref{eqn:Dyson3})
with $\beta=1, 2$, and 4 as the eigenvalue processes
of matrix-valued stochastic processes
in order to realize the point processes in equilibrium called
the Gaussian orthogonal ensemble (GOE), 
the Gaussian unitary ensemble (GUE),
and the Gaussian symplectic ensemble (GSE),
which are extensively studied in
\textbf{random matrix theory}
\cite{Meh04,For10,AGZ10,ABD11}. 

For $\beta=2$ with given $N \in \N$, we prepare
$N$-tuples of BMs
$\{B_{ii}^{x_i}(t) \}_{i=1}^{N}$, $t \geq 0$, 
each of which starts from $x_i \in \R$, and
$N(N-1)/2$-tuples of pairs of BMs
$\{B_{ij}(t), \widetilde{B}_{ij}(t) \}_{1 \leq i< j \leq N}$, $t \geq 0$, 
starting from the origin.
Here, there is a total of $N+2 \times N(N-1)/2=N^2$ BMs, 
each of them independent from the rest.
Then consider an $N \times N$ 
\textbf{Hermitian-matrix-valued Brownian motion}, 
\begin{equation}
\sB^{\x}(t) = \left( \begin{array}{cccc}
B_{11}^{x_1}(t) 
& \displaystyle{\frac{B_{12}(t)+\sqrt{-1} \widetilde{B}_{12}(t)}{\sqrt{2}}}
& \cdots & 
\displaystyle{\frac{B_{1N}(t)+\sqrt{-1} \widetilde{B}_{1N}(t)}{\sqrt{2}}}
\cr
\displaystyle{\frac{B_{12}(t)-\sqrt{-1} \widetilde{B}_{12}(t)}{\sqrt{2}}}
& B_{22}^{x_2}(t) & \cdots &
\displaystyle{\frac{B_{2N}(t)+\sqrt{-1} \widetilde{B}_{2N}(t)}{\sqrt{2}}}
\cr
\cdots & \cdots & \cdots & \cdots \cr
\displaystyle{\frac{B_{1N}(t)-\sqrt{-1} \widetilde{B}_{1N}(t)}{\sqrt{2}}}
& 
\displaystyle{\frac{B_{2N}(t)-\sqrt{-1} \widetilde{B}_{2N}(t)}{\sqrt{2}}}
 & \cdots & B_{NN}^{x_N}(t) 
\end{array} \right).
\label{eqn:Mxt}
\end{equation}

Remember that when we introduced BES$_{D}$ in Section \ref{sec:BES},
we considered the $D$-dimensional vector-valued Brownian motion in $\R^N$,
(\ref{eqn:DBM}), by preparing $D$-tuples of
independent BMs for its elements.
Here we consider the space of $N \times N$ Hermitian matrices
denoted by $\cH(N)$.
Since the dimension of this space is
${\rm dim} \, \cH(N)=N^2$,
we need $N^2$ independent BMs for elements
to describe a Brownian motion in this space $\cH(N)$.
Hence we can regard 
the process $\sB^{\x}(t), \ t \geq 0$ defined by (\ref{eqn:Mxt})
as a `Brownian motion in $\cH(N)$'.
By definition, the initial state of this Brownian motion is
the diagonal matrix 
\begin{equation}
\sB^{\x}(0)={\rm diag} (x_1, x_2, \dots, x_N).
\label{eqn:Mx0}
\end{equation}
We assume $x_1 \leq x_2 \leq \cdots \leq x_N$. 

Corresponding to calculating the absolute value (\ref{eqn:BES1}) of 
$\B^{\x}(t)$, by which BES$_{D}$ was introduced, 
here we calculate the eigenvalues of $\sB^{\x}(t)$.
For any $t \geq 0$, there is a family of $N \times N$
unitary matrices $\{\sU(t)\}$ which diagonalizes $\sB^{\x}(t)$,
$$
\sU^*(t) \sB^{\x}(t) \sU(t) =
{\rm diag} (\lambda_1(t), \dots, \lambda_N(t)) =: \Lambda(t),
\quad t \geq 0.
$$
Here for a matrix $\sM=(M_{ij})_{1 \leq i, j \leq N}$,
we define its Hermitian conjugate by 
$\sM^*=(\overline{M_{ji}})_{1 \leq i, j \leq N}$,
where $\overline{z}$ denotes the complex conjugate of $z \in \C$.
Consider a subspace of $\R^N$ defined by
\begin{equation}
\W_N
:= \{\x=(x_1, x_2, \dots, x_N)
\in \R^N: x_1 < x_2 < \cdots < x_N\},
\label{eqn:Weyl1}
\end{equation}
which is called the 
\textbf{Weyl chamber}
in representation theory.
If we impose the condition
$(\lambda_i(t))_{i=1}^N \in \W_N$,
$\sU(t)$ is uniquely determined at each time $t \geq 0$.

The following theorem is established.
\begin{thm}
\label{thm:Dyson_main}
The
\textbf{eigenvalue process}
$(\lambda_i(t))_{i=1}^N, \, t \geq 0$
of the Hermitian-matrix-valued Brownian motion (\ref{eqn:Mxt})
started at (\ref{eqn:Mx0}) satisfies the SDEs, 
\begin{equation}
d \lambda_i(t)=dB_i^{x_i}(t)
+ \sum_{\substack{1 \leq j \leq N, \cr j \not=i}}
\frac{dt}{\lambda_i(t)-\lambda_j(t)},
\quad t \geq 0, \quad i=1, \dots, N,
\label{eqn:lambda}
\end{equation}
where $(B_i^{x_i}(t))_{i=1}^{N}, \ t \geq 0$ are independent
BMs different from the $N^2$-tuples of BMs used to define $\sB^{\x}(t)$
in (\ref{eqn:Mxt}).
That is, this process realizes DYS$_{2}$. 
\end{thm}
The correspondence 
between BES$_{3}$ and DYS$_{2}$
in terms of equivalent processes 
is summarized as follows.
$$
\begin{array}{lccl}
{\mbox{\bf [Aspect 1]}} &  &  &  \cr
  &  &  & \quad
  \mbox{radial coordinate of} \cr
  & \mbox{BES$_{3}$}  & \Longleftrightarrow &
  \quad \mbox{$D=3$ vector-valued} \cr
  &  & & 
  \quad \mbox{Brownian motion} \cr
  &  &  & \cr
  & \mbox{} \quad & & \quad \mbox{eigenvalue process of} \cr
  & \mbox{$N$-particle DYS$_{2}$} \quad & \Longleftrightarrow & \quad \mbox{$N \times N$ Hermitian-matrix-valued} \cr
  & & & \quad \mbox{Brownian motion} \cr
\end{array}
$$

Dyson derived (\ref{eqn:lambda}) by applying the perturbation theory in 
quantum mechanics \cite{Dys62}.
Since $(\lambda_i(t))_{i=1}^N, \ t \geq 0$
are functionals of $\{B_{ii}^{x_i}(t),
B_{ij}(t), \widetilde{B}_{ij}(t) \}_{1 \leq i < j \leq N}$, $t \geq 0$, 
we can use It\^o's formula
to prove Theorem \ref{thm:Dyson_main}.
A key point to prove the theorem is applying
\textbf{It\^o's rule for} differentiating the product of
\textbf{matrix-valued semi-martingales} : 
If $\sX(t)=(X_{ij}(t))$ and $\sY(t)=(Y_{ij}(t))$ are $N \times N$ matrices
with semi-martingale elements, 
then
\begin{equation}
d( \sX^*(t) \sY(t))= d \sX^*(t) \sY(t) + \sX^{*}(t) d \sY(t) +
\langle d \sX^*, d \sY \rangle_t, \quad t \geq 0,
\label{eqn:ItoM1}
\end{equation}
where $\langle d \sX^*, d \sY \rangle_t$ 
denotes an $N \times N$ matrix-valued process,
whose $(i,j)$-th element is given by the finite-variation process
$\sum_k \langle d \overline{X_{k i}}, d Y_{k j} \rangle_t$, 
$1 \leq i, j \leq N$.

\subsection{Determinantal martingale representation (DMR)
of DYS$_2$}
\label{sec:determinantal}
\subsubsection{General formula}
\label{sec:general_DMR}

Denote the $N$-dimensional Laplacian 
with respect to the variables
$\x=(x_1, \dots, x_N)$ by
\[
\Delta^{(N)} := \sum_{i=1}^{N} \frac{\partial^2}{\partial x_i^2}.
\]
We set
\begin{equation}
h_N(\x):=\prod_{1 \leq i < j \leq N} (x_j-x_i)
=\det_{1 \leq i, j \leq N} [x_j^{i-1}], 
\label{eqn:hN0}
\end{equation}
where the determinant appearing here is
called the
\textbf{Vandermonde determinant}. 
As we have seen in \ref{sec:BES3}, 
$h_1(x) := x$ is a positive harmonic function in $\R_+=(0, \infty)$ 
which satisfies the boundary condition $h_1(0)=0$.
Similarly, we can see that
\begin{equation}
\Delta^{(N)} h_N(\x) =0. 
\label{eqn:hN3}
\end{equation}
and
\begin{equation}
h_N(\x) >0, \, \mbox{if $\x \in \W_N$}, \quad \mbox{and} \quad
h_N(\x)=0, \, \mbox{if $\x \in \partial \W_N$}.
\label{eqn:hN2}
\end{equation}

We have another correspondence between BES$_{3}$ and
DYS$_{2}$. 
$$
\begin{array}{lccl}
{\mbox{\bf [Aspect 2]}} &  &  &  \cr
  & \mbox{BES$_{3}$} \quad & \Longleftrightarrow & 
  \quad \mbox{$h$-transformation by $h_1$ of absorbing BM in $\R_+$} \cr
  &                    & \Longleftrightarrow &
  \quad \mbox{BM conditioned to stay positive} \cr
  &  &  & \cr
  & \mbox{$N$-particle DYS$_2$} \quad & \Longleftrightarrow &
  \quad \mbox{$h$-transformation by $h_N$ of absorbing BM 
  in $\W_N$} \cr
  & & \Longleftrightarrow &
  \quad \mbox{noncolliding BM}
\end{array}
$$
\vskip 0.5cm

As {\bf [Aspect 2]}, the Dyson model is constructed as the
$h$-transformation of the absorbing Brownian motion in $\W_N$.
Therefore, at any positive time $t > 0$ the configuration
is an element of $\W_N$,
\begin{equation}
\X(t)=(X_1(t), X_2(t), \dots, X_N(t)) \in
\W_N \quad t > 0,
\label{eqn:Xt}
\end{equation}
and hence there are no multiple points at which
coincidence of particle positions, 
$X_i(t)=X_j(t), i \not= j$, occurs.
We can consider, however, the Dyson model
starting from initial configurations with multiple points.
In order to describe configurations with multiple points, 
we represent each particle configuration by a sum 
of delta measures in the form
\begin{equation}
\xi(\cdot)=\sum_{i \in \I} \delta_{x_i}(\cdot)
\label{eqn:xi}
\end{equation}
with a sequence of points in $\R,
\x=(x_i)_{i \in \I}$,
where $\I$ is a countable index set.
Here for $y \in \R$,
$\delta_y(\cdot)$ denotes the delta measure such that
$\delta_y(\{x\})=1$ for $x=y$ and
$\delta_y(\{x\})=0$ otherwise.
Then, for (\ref{eqn:xi}) and $A \subset \R$,
$\xi(A) = \int_{A} \xi(dx)
=\sum_{i \in \I: x_i \in A} 1
=\sharp\{ x_i, x_i \in A\}$.

The measures of the form (\ref{eqn:xi})
satisfying the condition
$\xi(K) < \infty$ for any compact subset $K \subset \R$ are called the
\textbf{nonnegative integer-valued Radon measures} on $\R$
and we denote the space they form by $\Conf(\R)$;
\begin{equation}
\Conf(\R):=\left\{
\xi=\sum_{i \in \I} \delta_{x_i} : x_i \in \R, 
\xi(K) < \infty \, \, 
\mbox{for all bounded set $K \subset \R$} \right\}.
\label{eqn:conf}
\end{equation}
The set of configurations without multiple points
is denoted by 
\[
\Conf_0(\R):=\{\xi \in \Conf(\R) : \xi(\{x\}) \leq 1, ^{\forall} x \in \R\}.
\]
There is a trivial correspondence between $\W_N$
and $\Conf_0(\R)$.
We call $\x \in \R^N$ a
\textbf{labeled configuration}
and $\xi \in \Conf(\R)$ an 
\textbf{unlabeled configuration}.

We introduce a sequence of independent BMs, 
$\B^{\x}(t)=(B_i^{x_i}(t))_{i \in \I}, \ t \geq 0$, 
in $(\Omega, \cF, \rP^{\x})$ with expectation
written as $\rE^{\x}$.

In this section we assume that 
$\xi=\sum_{i \in \I} \delta_{x_i} \in \Conf_0(\R)$,
$\xi(\R)=N \in \N$
and consider DYS$_2$ as an $\Conf_0(\R)$-valued 
diffusion process,
\begin{equation}
\Xi(t, \cdot)=\sum_{i=1}^N \delta_{X_i(t)} (\cdot),
\quad t \geq 0, 
\label{eqn:Xit}
\end{equation}
starting from the initial configuration $\xi=\sum_{i=1}^N \delta_{x_i}$,
where $\X(t)=(X_1(t), \cdots, X_N(t))$ is the
solution of (\ref{eqn:Dyson3}) with $\beta=2$ under the initial 
configuration $\x=(x_1, \dots, x_N) \in \W_N$.
We write the process as $(\Xi, \P^{\xi})$
and express the expectation with respect to
the probability law $\P^{\xi}$ of the Dyson model by
$\E^{\xi}[\, \cdot \,]$.
We introduce a filtration $\{(\cF_{\Xi})_t\}_{t \in [0, \infty)}$
on the space of continuous paths
$\rC([0, \infty) \to \Conf(\R))$ defined by
$(\cF_{\Xi})_t=\sigma(\Xi(s), s \in [0, t])$, where
$\sigma$ denotes the smallest $\sigma$-field.

{\bf [Aspect 2]} of the Dyson model is expressed by
the following equality: 
for any $(\cF_{\Xi})_t$-measurable bounded function $F$, 
$0 \leq t \leq T < \infty$,
\begin{equation}
\E^{\xi}[F(\Xi(\cdot))]
=\rE^{\x} \left[ F\left( \sum_{i=1}^N \delta_{B_i(\cdot)} \right)
\1_{(\tau > T)}
\frac{h_N(\B(T))}{h_N(\x)} \right],
\label{eqn:Asp2_1}
\end{equation}
where 
\begin{equation}
\tau := \inf \{t > 0: \B^{\x}(t) \notin \W_N \}
\label{eqn:tau_collision}
\end{equation} 
and 
we have assumed the relations
$\xi=\sum_{i=1}^N \delta_{x_i} \in \Conf_0(\R),
\x=(x_1, \dots, x_N) \in \W_N$
and (\ref{eqn:Xit}).

In the following lemma, we claim that 
even if we delete the indicator $\1_{(\tau > T)}$
in RHS of (\ref{eqn:Asp2_1}),
still the equality holds.
It is a multivariate extension of the claim by which
we replaced (\ref{eqn:BES3_b1}) by (\ref{eqn:BES3_b2})
in Section \ref{sec:BES3}.

\begin{lem}
\label{thm:h_trans_iden}
Assume that $\xi=\sum_{i=1}^N \delta_{x_i} \in \Conf_0(\R)$,
$\x=(x_1, \dots, x_N) \in \W_N$.
For any $(\cF_{\Xi})_t$-measurable bounded function $F$, 
$0 \leq t \leq T < \infty$,
\begin{equation}
\E^{\xi}[F(\Xi(\cdot))]
=\rE^{\x} \left[ F\left( \sum_{i=1}^N \delta_{B_i(\cdot)} \right)
\frac{h_N(\B(T))}{h_N(\x)} \right].
\label{eqn:Asp2_2}
\end{equation}
\end{lem}

In Section \ref{sec:poly_mar}, we introduced 
the fundamental martingale polynomials
associated with BM, $\{m_n(t, x)\}_{n \in \N_0}$.
Since they are monic polynomials, we can prove the equalities
\begin{eqnarray}
\frac{h_N(\y)}{h_N(\x)}
&=& \frac{1}{h_N(\x)} \det_{1 \leq i, j \leq N} [y_j^{i-1}]
\nonumber\\
&=& \frac{1}{h_N(\x)} \det_{1 \leq i, j \leq N}[ m_{i-1}(t, y_j)]
\label{eqn:hB1}
\end{eqnarray}
for an arbitrary $t \in [0, \infty)$.
This implies that $(h_N(\B(t))/h_N(\x))_{t \geq 0}$ is a martingale.

Here we extend the integral transformation
defined by (\ref{eqn:I1}) with (\ref{eqn:Ghat1}) in
Section \ref{sec:poly_mar} to a linear integral transformation
of multivariate functions as follows.
When $F^{(i)}(\x)=\prod_{j=1}^N f_j^{(i)}(x_j)$, $i=1,2$
are given for $\x=(x_1, \dots, x_N) \in \R^N$, then we define
$$
\cI \left[ F^{(i)}(\bW) \left| \{(t_{\ell}, x_{\ell})\}_{\ell=1}^N  \right. \right]
:= \prod_{j=1}^N \cI \left[ \left. f^{(i)}_j(W_j) \right| (t_{j}, x_{j}) \right],
\quad i=1,2,
$$
and
\begin{eqnarray}
&& \cI \Big[c_1 F^{(1)}(\bW) +c_2 F^{(2)}(\bW) \left| \{(t_{\ell}, x_{\ell})\}_{\ell=1}^N 
\right. \Big]
\nonumber\\
&& \quad
:= c_1 \cI \Big[F^{(1)}(\bW) \left| \{(t_{\ell}, 
x_{\ell})\}_{\ell=1}^N \right. \Big]
+ c_2 \cI \Big[F^{(2)}(\bW) \left| \{(t_{\ell}, 
x_{\ell})\}_{\ell=1}^N \right. \Big],
\nonumber
\end{eqnarray}
$c_1, c_2 \in \C$, 
for $0 < t_i < \infty, i=1, \dots, N$, 
where $\bW=(W_1, \dots, W_N) \in \R^N$.
In particular, if $t_{\ell}=t, 1 \leq {^{\forall}\ell} \leq N$, we write
$\cI[\, \cdot \,| \{(t_{\ell}, x_{\ell})\}_{\ell=1}^N]$ simply as
$\cI[\, \cdot \,|(t, \x)]$ with $\x=(x_1, \dots, x_N)$.
Then (\ref{eqn:hB1}) is further rewritten as
\begin{eqnarray}
\frac{h_N(\y)}{h_N(\x)}
&=& \frac{1}{h_N(\x)}
\det_{1 \leq i, j \leq N}
\Big[ \cI [(W_j)^{i-1} | (t, y_j) ] \Big]
\nonumber\\
&=& \cI \left[ \left. \frac{1}{h_N(\x)}
\det_{1 \leq i, j \leq N} [(W_j)^{i-1}] \right| (t, \y) \right]
\nonumber\\
&=& \cI \left[ \left.
\frac{h_N(\bW)}{h_N(\x)} \right| (t, \y) \right],
\label{eqn:hB2}
\end{eqnarray}
where the multilinearity of determinants has been used.

Now we use the following determinant identity.
\begin{lem}
\label{thm:det_iden}
For $\x=(x_1, \dots, x_N) \in \W_N$,
$\z=(z_1, \dots, z_N) \in \C^{N}$,
\begin{equation}
\frac{h_N(\z)}{h_N(\x)}
=\det_{1 \leq i, j \leq N}
\Big[ \Phi_{\xi}^{x_i}(z_j) \Big],
\label{eqn:det_iden}
\end{equation}
where
\begin{equation}
\Phi_{\xi}^{u}(z)=\prod_{\substack{1 \leq k \leq N, \cr x_k \not= u}}
\frac{z-x_k}{u-x_k}, \quad
\mbox{for $\displaystyle{
\xi=\sum_{i=1}^N \delta_{x_i} \in \Conf_0(\R), z, u \in \C}$}. 
\label{eqn:Phi1}
\end{equation}
\end{lem}

Then (\ref{eqn:hB2}) is written as
\begin{eqnarray}
\frac{h_N(\y)}{h_N(\x)}
&=& \cI \left[
\det_{1 \leq i, j \leq N}
[\left. \Phi_{\xi}^{x_i}(W_j)] \right| (t, \y) \right]
\nonumber\\
&=& \det_{1 \leq i, j \leq N}
[ \cM_{\xi}^{x_i}(t, y_j) ],
\label{eqn:hB3}
\end{eqnarray}
where
\begin{equation}
\cM_{\xi}^{x}(t, y)
:=\cI[ \Phi_{\xi}^{x}(W)|(t, y)], \quad
x, y \in \R, \quad t \geq 0.
\label{eqn:M1}
\end{equation}

\begin{prop}
\label{thm:DMR}
Assume that $\xi=\sum_{i=1}^N \delta_{x_i} \in \Conf_0(\R)$.
The following are satisfied by (\ref{eqn:M1}). 
\begin{description}
\item{\rm (i)} \, 
$(\cM_{\xi}^{x_i}(t, B(t)))_{t \geq 0}, i=1, \dots N$
are continuous martingales.

\item{\rm (ii)} \, For any time $t \geq 0$, $\cM_{\xi}^{x_i}(t, y),
i=1, \dots, N$
are linearly independent functions of $y$.

\item{\rm (iii)} \, For $1 \leq i, j \leq N$,
$\lim_{t \downarrow 0} \rE^{x_i}[\cM_{\xi}^{x_j}(t, B(t))]=\delta_{ij}$.
\end{description}
Then for any $(\cF_{\Xi})_t$-measurable bounded function $F$, 
$0 \leq t \leq T < \infty$, the equality
\begin{equation}
\E^{\xi}[F(\Xi(\cdot))]
=\rE^{\x} \left[ F\left( \sum_{i=1}^N \delta_{B_i(\cdot)} \right)
\cD_{\xi}(T, \B(T)) \right]
\label{eqn:DMR1}
\end{equation}
holds, where
\begin{equation}
\cD_{\xi}(t, \y)=\det_{1 \leq i, j \leq N}
[ \cM_{\xi}^{y_j}(t, y_i)], \quad
\y=(y_1, \dots, y_N) \in \W_N, \quad t \geq 0.
\label{eqn:D1}
\end{equation}
\end{prop}

We remark that $\cD_{\xi}(t, \B(t)), \ t \geq 0$ is indeed a
continuous martingale by part (i) and
is not identically zero by part (ii) of Proposition \ref{thm:DMR}.
We call $\cD_{\xi}(t, \B(t)), \ t \geq 0$ a 
\textbf{determinantal martingale}
and the equality (\ref{eqn:DMR1}) 
the 
\textbf{determinantal martingale representation} (DMR) 
of DYS$_2$ \cite{Kat14}.

For $n \in \N$, an index set $\{1,2, \dots, n\}$
is denoted by $\I_{n}$.
Fixing $N \in \N$ with $N' \in \I_N$, 
we write $\J \subset \I_N, \sharp \J=N'$,
if $\J=\{j_1, \dots, j_{N'}\},
1 \leq j_1 < \dots < j_{N'} \leq N$.
For  $\x=(x_1, \dots, x_N) \in \R^N$, 
put $\x_{\J}=(x_{j_1}, \dots, x_{j_{N'}})$.
In particular, we write
$\x_{N'}=\x_{\I_{N'}}, 1 \leq N' \leq N$.
(By definition $\x_N=\x$.)
A collection of all permutations of 
elements in $\J$ is denoted by $\S(\J)$.
In particular, we write $\S_{N'}=\S(\I_{N'}), 1 \leq N' \leq N$.

The following shows the
\textbf{reducibility} of the determinantal martingale
in the sense that,
if we observe a symmetric function depending
on $N'$ variables, $N' \leq N$,
then the size of determinantal
martingale can be reduced from $N$ to $N'$.

\begin{lem}
\label{thm:reduce}
Assume that
$\xi=\sum_{i=1}^N \delta_{x_i}=\sum_{i \in \I_N} \delta_{x_i}$
with $\x \in \W_N$.
Let $1 \leq N' \leq N$.
For $0 < t \leq T < \infty$ and an 
$(\cF_{\Xi})_t$-measurable symmetric function
$F_{N'}$ on $\R^{N'}$,
\begin{eqnarray}
&& \sum_{\J \subset \I_N, \sharp \J=N'}
\rE^{\x} \left[
F_{N'}(\B_{\J}(t))
\cD_{\xi}(T, \B(T)) \right]
\nonumber\\
&& \quad
= \int_{\W_{N'}} \xi^{\otimes N'} (d\bv)
\rE^{\bv} \left[
F_{N'}(\B_{N'}(t))
\cD_{\xi}(T, \B_{N'}(T)) \right].
\label{eqn:reducibility}
\end{eqnarray}
\end{lem}

\subsubsection{Time-dependent density function $\rho_{\xi}(t, x)$}
\label{sec:density}

The 
\textbf{density function} at a single time for 
$(\Xi, \P^{\xi}), \xi \in \Conf_0(\R)$ is denoted by
$\rho_{\xi}(t, x)$.
It is defined as a continuous function of $x \in \R$ for
$0 \leq t \leq T < \infty$ such that for any 
\textbf{test function}, 
$\chi \in \rC_{\rm c}(\R)$,
\begin{equation}
\E^{\xi} \left[ \int_{\R} \chi(x) \Xi(t, dx) \right]
=\int_{\R} dx \, \chi(x) \rho_{\xi}(t, x).
\label{eqn:rho1}
\end{equation}
The test function $\chi$ is symmetrized as
$\displaystyle{g(\x)=\sum_{i=1}^N \chi(x_i)}$,
which is applied as $F$ to the DMR (\ref{eqn:DMR1}),
and we obtain the equality
\begin{equation}
\E^{\xi} \left[
\sum_{i=1}^N \chi(X_i(t)) \right]
=\rE^{\x} \left[ \sum_{i=1}^N \chi(B_i(t)) \cD_{\xi}(T, \B(T)) \right],
\quad 0 \leq t \leq T < \infty.
\label{eqn:DMR_ex1}
\end{equation}
The LHS of (\ref{eqn:DMR_ex1}) gives
$$
\E^{\xi} \left[ \sum_{i=1}^N \chi(X_i(t)) \right]
=\E^{\xi} \left[ \int_{\R} \chi(x) \Xi(t, dx) \right]
$$
by (\ref{eqn:Xit}).
On the other hand, RHS of (\ref{eqn:DMR_ex1}) 
is reduced by Lemma \ref{thm:reduce} as
\begin{eqnarray}
\sum_{i=1}^N \rE^{\x} [ \chi(B_i(t)) \cD_{\xi}(T, \B(T)) ]
&=& \int_{\R} \xi(dv) \rE^v [ \chi(B(t)) \cM_{\xi}^v(t, B(t))]
\nonumber\\
&=& \int_{\R} \xi(dv) \int_{\R} dx \, \chi(x) p(t, x|v) \cM_{\xi}^v(t, x).
\nonumber
\end{eqnarray}
By Fubini's theorem, we can rewrite it as
$
\int_{\R} dx \, \chi(x) \cG_{\xi}(t, x; t, x),
$
where
\begin{equation}
\cG_{\xi}(s,x;t,y)
=\int_{\R} \xi(dv) p(s, x|v) \cM_{\xi}^{v}(t, y).
\label{eqn:cG1}
\end{equation}
Then (\ref{eqn:rho1}) gives
\begin{equation}
\rho_{\xi}(t, x)=\cG_{\xi}(t,x;t,x), \quad
x \in \R, \quad t \geq 0.
\label{eqn:rho1b}
\end{equation}

\subsubsection{Two-time correlation function $\rho_{\xi}(s, x; t, y)$}
\label{sec:two-time}

For $0 \leq t_1 < t_2 \leq T < \infty$, set
$$
g_1(\x)=\sum_{i=1}^N \chi_1(x_i), \quad
g_2(\x)=\sum_{i=1}^N \chi_2(x_i),
$$
where $\chi_m \in \rC_{\rm c}(\R), m=1,2$, and put
$$
F(\Xi(\cdot))=\prod_{m=1}^2 g_m(\X(t_m)).
$$
If we apply this to DMR, (\ref{eqn:DMR1}), 
we obtain the equality 
\begin{eqnarray}
&& \E^{\xi} \left[ \sum_{i=1}^N \sum_{j=1}^N \chi_1(X_i(t_1)) \chi_2(X_j(t_2)) \right]
\nonumber\\
&& \quad
= \rE^{\x} \left[ \sum_{i=1}^N \sum_{j=1}^N
\chi_1(B_i(t_1)) \chi_2(B_j(t_2)) 
\cD_{\xi}(T, \B(T)) \right]\hspace{-0.1cm},
\, 0 \leq t \leq T < \infty.
\label{eqn:DMR_ex2}
\end{eqnarray}
The LHS of (\ref{eqn:DMR_ex2}) 
defines the two-time correlation function
$\rho_{\xi}(s,x;t,y)$ as
\begin{equation}
\E^{\xi} \left[
\sum_{i=1}^N \sum_{j=1}^N \chi_1(X_i(t_1)) \chi_2(X_j(t_2)) \right]
=\int_{\R^2} dx_1 dx_2 \,
\chi_1(x_1) \chi_2(x_2) \rho_{\xi}(t_1, x_1; t_2, x_2).
\label{eqn:two-time}
\end{equation}
On the other hand, RHS of (\ref{eqn:DMR_ex2}) gives
\begin{eqnarray}
&& \sum_{i=1}^N \sum_{j=1}^N \rE^{\x}
[ \chi_1(B_i(t_1)) \chi_2(B_j(t_2)) \cD_{\xi}(T, \B(T)) ]
\nonumber\\
&& \quad
= \sum_{\substack{1 \leq i, j \leq N, \cr i \not=j}}
\rE^{\x}[\chi_1(B_i(t_1)) \chi_2(B_j(t_2)) \cD_{\xi}(T, \B(T))]
\nonumber\\
&& \quad
+ \sum_{1 \leq i \leq N}
\rE^{\x}[\chi_1(B_i(t_1)) \chi_2(B_i(t_2)) \cD_{\xi}(T, \B(T))].
\nonumber
\end{eqnarray}
By the reducibility of DMR given by Lemma \ref{thm:reduce},
the last expression becomes
\begin{eqnarray}
&& \int_{\R^2} \xi^{\otimes 2}(d \bv) 
\nonumber\\
&& \quad \times \ 
\rE^{(v_1, v_2)}
\left[ \chi_1(B_1(t_1)) \chi_2(B_2(t_2))
\det \left( \begin{array}{ll}
\cM_{\xi}^{v_1}(T, B_1(T)) & \cM_{\xi}^{v_1}(T, B_2(T)) \cr
\cM_{\xi}^{v_2}(T, B_1(T)) & \cM_{\xi}^{v_2}(T, B_2(T))
\end{array} \right) \right]
\nonumber\\
&& + \int_{\R} \xi(dv) 
\rE^v [ \chi_1(B(t_1)) \chi_2(B(t_2)) \cM_{\xi}^v(T, B(T)) ]
\nonumber
\end{eqnarray}
If we use the martingale property (i) of $\cM_{\xi}^v$ 
in Proposition \ref{thm:DMR}, it is written as
\begin{eqnarray}
&& \int_{\R^2} \xi^{\otimes 2}(d \bv)
\nonumber\\
&& \quad \times \ 
\rE^{(v_1, v_2)} \left[ \chi_1(B_1(t_1)) \chi_2(B_2(t_2))
\det \left( \begin{array}{ll}
\cM_{\xi}^{v_1}(t_1, B_1(t_1)) & \cM_{\xi}^{v_1}(t_2, B_2(t_2)) \cr
\cM_{\xi}^{v_2}(t_1, B_1(t_1)) & \cM_{\xi}^{v_2}(t_2, B_2(t_2))
\end{array} \right) \right]
\nonumber\\
&& + \int_{\R} \xi(dv) 
\rE^v [ \chi_1(B(t_1)) \chi_2(B(t_2)) \cM_{\xi}^v(t_2, B(t_2)) ].
\nonumber
\end{eqnarray}
By Fubini's theorem, this is equal to
\begin{eqnarray}
&& \int_{\R^2} dx_1 dx_2 \,
\chi_1(x_1) \chi_2(x_2) 
\det \left( \begin{array}{ll}
\cG_{\xi}(t_1, x_1; t_1, x_1) & \cG_{\xi}(t_1, x_1; t_2, x_2) \cr
\cG_{\xi}(t_2, x_2; t_1, x_1) & \cG_{\xi}(t_2, x_2; t_2, x_2)
\end{array} \right)
\nonumber\\
&& + \int_{\R^2} dx_1 dx_2 \,
\chi_1(x_1) \chi_2(x_2) 
\cG_{\xi}(t_1, x_1; t_2, x_2) p(t_2-t_1, x_2|x_1)
\nonumber\\
&=&\int_{\R^2} dx_1 dx_2 \,
\chi_1(x_1) \chi_2(x_2) 
\nonumber\\
&& \qquad \times
\det \left( \begin{array}{ll}
\cG_{\xi}(t_1, x_1; t_1, x_1) & \cG_{\xi}(t_1, x_1; t_2, x_2) \cr
\cG_{\xi}(t_2, x_2; t_1, x_1)-p(t_2-t_1, x_2|x_1) & \cG_{\xi}(t_2, x_2; t_2, x_2)
\end{array} \right).
\nonumber
\end{eqnarray}
Since this is equal to (\ref{eqn:two-time}),
the two-time correlation function is determined as
\begin{equation}
\rho_{\xi}(s, x; t, y)
=\det \left( \begin{array}{ll}
\mbK_{\xi}(s, x; s, x) & \mbK_{\xi}(s, x; t, y) \cr
\mbK_{\xi}(t, y; s, x) & \mbK_{\xi}(t, y; t, y)
\end{array} \right)
\label{eqn:rho_2time}
\end{equation}
for $0 \leq s < t < \infty$, $x, y \in \R$,
where
\begin{equation}
\mbK_{\xi}(s, x; t, y)
= \cG_{\xi}(s, x; t, y)-\1_{(s>t)} p(s-t, x|y).
\label{eqn:K_a1}
\end{equation}

\subsection{Determinantal stochastic processes (DSPs)}
\label{sec:det_process}
\subsubsection{From DMR to DSP}
\label{sec:DMR_DP}

In the previous section, the density function
at a single time $\rho_{\xi}(t, x)$ and the two-time
(and two-point) correlation function
$\rho_{\xi}(s,x;t,y)$ were defined by
(\ref{eqn:rho1}) and (\ref{eqn:two-time}), respectively.
In order to give a general definition of 
\textbf{spatio-temporal correlation functions}
here we consider the 
\textbf{Laplace transformations of
the multitime joint distribution functions} of $(\Xi, \P^{\xi})$.
For any integer $M \in \N$,
a sequence of times
$\t=(t_1,\dots,t_M) \in [0, \infty)^M$ with 
$0 \leq t_1 < \cdots < t_M < \infty$, 
and a sequence of functions
$\f=(f_{t_1},\dots,f_{t_M}) \in \rC_{\rm c}(\R)^M$,
let
\begin{equation}
\Psi_{\xi}^{\t}[\f]
\equiv \E^{\xi} \left[ \exp \left\{ \sum_{m=1}^{M} 
\int_{\R} f_{t_m}(x) \Xi(t_m, dx) \right\} \right].
\label{eqn:GF1}
\end{equation}
By (\ref{eqn:Xit}), if we set 
test functions as
\begin{equation}
\chi_{t_m}(\cdot) = e^{f_{t_m}(\cdot)}-1, \quad
1 \leq m \leq M,
\label{eqn:test1}
\end{equation}
we can rewrite (\ref{eqn:GF1}) in the form
\begin{equation}
\Psi_{\xi}^{\t}[\f]
=\E^{\xi} \left[ \prod_{m=1}^M \prod_{i=1}^N
\{1+\chi_{t_m}(X_i(t_m)) \} \right].
\label{eqn:GF1b}
\end{equation}
We expand this with respect to test functions
and define the spatio-temporal correlation functions
$\{\rho_{\xi} \}$ as coefficients, 
\begin{equation}
\Psi_{\xi}^{\t}[\f]
=\sum_
{\substack
{0 \leq N_m \leq N, \\ 1 \leq m \leq M} }
\int_{\prod_{m=1}^{M} \W_{N_{m}^{\rm A}}}
\prod_{m=1}^{M} 
d \x_{N_m}^{(m)}
\prod_{i=1}^{N_{m}} 
\chi_{t_m} \Big(x_{i}^{(m)} \Big) 
\rho_{\xi} 
\Big( t_{1}, \x^{(1)}_{N_1}; \dots ; t_{M}, \x^{(M)}_{N_M} \Big),
\label{eqn:GF1c}
\end{equation}
where $\x^{(m)}_{N_m}$ denotes
$(x^{(m)}_1, \dots, x^{(m)}_{N_m})$
and
$d \x^{(m)}_{N_m}= \prod_{i=1}^{N_m} dx^{(m)}_i$,
$1 \leq m \leq M$.
Here the empty products equal 1 by convention
and the term with $N_m=0, 1 \leq \forall m \leq M$ 
is considered to be 1.
The previous two examples $\rho_{\xi}(t,x)$ and
$\rho_{\xi}(s,x; t, y)$ are special cases in which we set
$M=1, \ t_1=t, \ N_1=1, \ x_1^{(1)}=x$, and
$M=2, \ t_1=s, \ t_2=t, \ N_1=N_2=1, \ x_1^{(1)}=x, \ x_1^{(2)}=y$,
respectively.
The function $\Psi_{\xi}^{\t}[\f]$ is a 
\textbf{generating function of correlation functions}.

Given an integral kernel, 
$\bK(s,x;t,y)$, 
$(s,x), (t,y) \in [0, \infty) \times \R$,
and a sequence of functions
$(\chi_{t_1}, \dots, \chi_{t_M}) \in \rC_{\rm c}(\R)^M, \ M \in \N$,
the
\textbf{Fredholm determinant} associated with $\bK$ 
and $(\chi_{t_m})_{m=1}^M$ is defined as
\begin{eqnarray}
&& \mathop{{\rm Det}}_
{\substack{
(s,t)\in \{t_1, \dots, t_M\}^2, \\
(x,y)\in \R^2}
}
 \Big[\delta_{st} \delta_x(\{y\})
+ \bK(s,x;t,y) \chi_{t}(y) \Big]
\nonumber\\
&& 
=\sum_
{\substack
{0 \leq N_m \leq N, \\ 1 \leq m \leq M} }
\int_{\prod_{m=1}^{M} \W_{N_{m}}^{\rm A}}
\prod_{m=1}^{M} 
d \x_{N_m}^{(m)}
\prod_{k=1}^{N_{m}} 
\chi_{t_m} \Big(x_{k}^{(m)} \Big) 
\det_{\substack
{1 \leq i \leq N_{m}, 1 \leq j \leq N_{n}, \\
1 \leq m, n \leq M}
}
\Bigg[
\bK(t_m, x_{i}^{(m)}; t_n, x_{j}^{(n)} )
\Bigg].
\nonumber\\
\label{eqn:Fredholm1}
\end{eqnarray}

\begin{df}
\label{thm:determinantal}
If any moment generating function (\ref{eqn:GF1}) 
is given by a Fredholm determinant,
the process $(\Xi, \P^{\xi})$ is said to be a
\textbf{determinantal stochastic process} (DSP). 
In this case
all spatio-temporal correlation functions
are given by determinants as 
\begin{equation}
\rho_{\xi} \Big(t_1,\x^{(1)}_{N_1}; \dots;t_M,\x^{(M)}_{N_M} \Big) 
=\det_{\substack
{1 \leq i \leq N_{m}, 1 \leq j \leq N_{n}, \\
1 \leq m, n \leq M}
}
\Bigg[
\mbK_{\xi}(t_m, x_{i}^{(m)}; t_n, x_{j}^{(n)} )
\Bigg],
\label{eqn:rho1d}
\end{equation}
$0 \leq t_1 < \cdots < t_M < \infty$, 
$1 \leq N_m \leq N$,
$\x^{(m)}_{N_m} \in \R^{N_m}, 1 \leq m \leq M \in \N$.
Here the integral kernel, $\mbK_{\xi}:([0, \infty) \times \R)^2 \mapsto \R$,
is a function of the initial configuration
$\xi$ and is called
the 
\textbf{correlation kernel}.
\end{df}
\vskip 0.3cm
\begin{rem}
\label{eqn:remark1}
If the process $(\Xi, \P^{\xi})$ is DSP,
then, at each time $0 \leq t < \infty$,
all spatial correlation functions are
given by determinants as
\begin{equation}
\rho_{\xi}(\x_{N'})
=\det_{1 \leq i, j \leq N'} [ \rK(x_i, x_j)],
\quad 1 \leq N' \leq N,
\label{eqn:dpp1}
\end{equation}
with $\rK(x, y)=\mbK_{\xi}(t, x; t, y)$.
In general a random integer-valued
Radon measure in $\Conf(\R)$ (resp. $\Conf_0(\R)$) is called a
\textbf{point process}
(resp. \textbf{simple point process}).
A simple point process is said to be a
\textbf{determinantal point process} (DPP)
(or 
\textbf{Fermion point process})
with kernel $\rK$,
if its spatial correlation functions exist
and are given in the form (\ref{eqn:dpp1}).
When $\rK$ is symmetric,
{\it i.e.},
$\rK(x,y)=\rK(y,x), \ x, y \in \R$,
Soshnikov \cite{Sos00}
and Shirai and Takahashi \cite{ST00,ST03a,ST03b} 
gave sufficient
conditions for $\rK$ to be a correlation kernel
of a determinantal point process.
See also \cite{KS22}. 
The notion of DSP
given by Definition \ref{thm:determinantal}
is a dynamical extension of 
the determinantal point process \cite{BR05,KT07}.
\end{rem}
\vskip 0.3cm

Let $\Xi$ and $\widetilde{\Xi}$
be determinantal processes with
correlation kernels $\mbK$ and $\widetilde{\mbK}$,
respectively.
If there is a function $G(s,x)$, which is continuous 
with respect to $x \in \R$ for any fixed $s \in [0, \infty)$,
such that
\begin{equation}
\mbK(s,x;t,y)=\frac{G(s,x)}{G(t,y)}
\widetilde{\mbK}(s,x;t,y),
\quad (s, x), (t, y) \in [0, \infty) \times \R,
\label{eqn:gauge1}
\end{equation}
then we say that
$\widetilde{\Xi}$ is a 
\textbf{gauge transform} of $\Xi$.

\begin{lem}
\label{thm:factor_K}
Under (\ref{eqn:gauge1}), 
\[
(\Xi(t))_{t \geq 0} \law= (\widetilde{\Xi}(t))_{t \geq 0}.
\]
In other words, DSP is \textbf{gauge invariant}. 
\end{lem}

By Proposition \ref{thm:DMR},
(\ref{eqn:GF1b}) has the DMR 
\begin{equation}
\Psi_{\xi}^{\t}[\f]
=\rE^{\x} \left[
\prod_{m=1}^N \prod_{i=1}^N
\{1+\chi_{t_m}(B_i(t_m)) \} 
\cD_{\xi}(T, \B(T)) \right],
\label{eqn:GF1e}
\end{equation}
where $T \geq t_M$ and $\xi=\sum_{i=1}^N \delta_{x_i}$.
The following equality is established \cite{Kat14}.

\begin{lem}
\label{thm:equality}
Let $\x \in \W_N^{\rm A}$ and
$\xi=\sum_{i=1}^N \delta_{x_i}$.
Then for any $M \in \N$, $0 \leq t_1 < \cdots < t_M \leq T < \infty$,
$\chi_{t_m} \in \rC_{\rm c}(\R), 1 \leq m \leq M$,
the equality
\begin{eqnarray}
&& \rE^{\x} \left[
\prod_{m=1}^M \prod_{i=1}^N
\{1+\chi_{t_m}(B_i(t_m)) \}
\cD_{\xi}(T, \B(T)) \right]
\nonumber\\
&& \qquad =
\mathop{{\rm Det}}_
{\substack{
(s,t)\in \{t_1, \dots, t_M\}^2, \\
(x,y)\in \R^2}
}
 \Big[\delta_{st} \delta_x(\{y\})
+ \mbK_{\xi}(s,x;t,y) \chi_{t}(y) \Big]
\label{eqn:Fredholm2}
\end{eqnarray}
holds, 
where $\mbK_{\xi}$ is given by (\ref{eqn:K_a1})
with (\ref{eqn:cG1}). 
\end{lem}

Now we arrive at the following full characterization of
DYS$_2$, $(\Xi, \P^{\xi})$.

\begin{thm}
\label{thm:det_process_main}
For any finite and fixed initial configuration without multiple points,
that is, for $\xi \in \Conf_0(\R), \xi(\R) = N \in \N$,
DYS$_2$ is determinantal.
Its correlation kernel is given by
\begin{equation}
\mbK_{\xi}(s, x; t, y)
= \cG_{\xi}(s, x; t, y)-\1_{(s>t)} p(s-t, x|y),
\quad (s,x), (t,y) \in [0, \infty) \times \R 
\label{eqn:Kernel1}
\end{equation}
with
\begin{equation}
\cG_{\xi}(s,x;t,y)
=\int_{\R} \xi(dv) p(s, x|v) \cM_{\xi}^{v}(t, y).
\label{eqn:Kernel2}
\end{equation}
\end{thm}

\subsubsection{Generalization for initial configuration
with multiple points}
\label{sec:multi_point}

For general $\xi=\sum_{i=1}^N \delta_{x_i} \in \Conf(\R)$
with $\xi(\R)=N < \infty$, define
$\supp \xi=\{x \in \R : \xi(x) > 0\}$ and let
$
\xi_{*}(\cdot)=\sum_{v \in \supp \xi} \delta_{v}(\cdot).
$
For $s \in [0, \infty)$, $v, x \in \R$, $z, \zeta \in \C$, let
\begin{equation}
\phi_{\xi}^{v}((s,x);z,\zeta)
=
\frac{p(s,x|\zeta)}{p(s,x|v)} \frac{1}{z-\zeta}
\prod_{i=1}^N
\frac{z-x_{i}}{\zeta-x_{i}},
\label{eqn:phi}
\end{equation}
and
\begin{eqnarray}
\Phi_{\xi}^{v}((s,x); z) &=& \frac{1}{2 \pi \sqrt{-1}}
\oint_{C(\delta_{v})} d \zeta \, 
\phi^v_{\xi}((s,x); z, \zeta)
\nonumber\\
&=& {\rm Res} \,\Big[\phi^v_{\xi}((s,x); z, \zeta); \zeta=v \Big],
\label{eqn:PhiB}
\end{eqnarray}
where
$C(\delta_{v})$ is a closed contour on the complex plane $\C$
encircling the point $v$ 
once in the positive direction. 
This function (\ref{eqn:PhiB}) is entire with respect to $z \in \C$
parameterized by $(s, x) \in [0, \infty) \times \R$
in addition to $v \in \C, \xi \in \Conf(\R)$.
Remark that the polynomial function
$\Phi_{\xi}^{u}(z)$ defined by (\ref{eqn:Phi1}) is parameterized
only by $u \in \C$ and $\xi \in \Conf_0(\R)$.
Here we start from this entire function and consider 
its $\cI$-transformation, 
\begin{equation}
\cM_{\xi}^{v}((s, x)|(t,y))
=\cI \left[\Phi_{\xi}^{v}((s,x); W) \Big| (t,y) \right],
\quad (s,x), (t,y) \in [0, \infty) \times \R,
\label{eqn:cMB}
\end{equation}
which provides a martingale, if we put
$y=B(t), \ t \geq 0$.
Then it is easy to see that (\ref{eqn:Kernel1}) with (\ref{eqn:Kernel2}) 
is rewritten as
\begin{equation}
\mbK_{\xi}(s,x;t,y)
=\int_{S} \xi_{*}(dv) p(s, x|v) \cM_{\xi}^{v}((s,x)|(t,y))
- \1(s>t) p(s-t,x|y),
\label{eqn:KernelB1}
\end{equation}
$(s,x), (t,y) \in [0, \infty) \times \R$.

We note that
the kernel (\ref{eqn:KernelB1}) with (\ref{eqn:cMB})
is bounded and integrable also for $\xi \in \Conf(\R) \setminus \Conf_0(\R)$.
Therefore, the spatio-temporal correlations are given
by (\ref{eqn:rho1d}) for any $0 \leq t_1 < \dots < t_M < \infty, M \in \N$
and the finite-dimensional distributions 
are determined.
\begin{prop}
\label{thm:entrance_law}
Also for $\xi \in \Conf(\R) \setminus \Conf_0(\R)$,
the DSPs
with the correlation kernels (\ref{eqn:KernelB1})
are well-defined.
\end{prop}
The complete proof of this proposition was given
in Section 4.1 of \cite{KT10}.
See also \cite{Osa12,Osa13,Tsa16,OT20}.
The above extension will provide the 
\textbf{entrance laws} for
the processes $(\Xi(t),  t>0, \P^{\xi})$
in the sense of Section XII.4 in \cite{RY05}.

\subsubsection{DSP with extended Hermitian kernel}
\label{sec:extHermite}

In order to give an example
of Proposition \ref{thm:entrance_law}, 
here we study the extreme case where
all $N$ points are concentrated on an origin,
\begin{equation}
\xi=N \delta_0 \quad
\Longleftrightarrow \quad
\xi_{*}=\delta_0 \quad
\mbox{with} \quad
\xi(\{0\})=N.
\label{eqn:xi01}
\end{equation}
For (\ref{eqn:xi01}), (\ref{eqn:phi}) and (\ref{eqn:PhiB})
become
\begin{eqnarray}
\phi_{N \delta_0}^0((s,x); z, \zeta)
&=& \frac{p(s,x|\zeta)}{p(s,x|0)}
\frac{1}{z-\zeta} \left( \frac{z}{\zeta} \right)^N
\nonumber\\
&=& \frac{p(s,x|\zeta)}{p(s,x|0)}
\sum_{\ell=0}^{\infty} \frac{z^{N-\ell-1}}{\zeta^{N-\ell}},
\nonumber
\end{eqnarray}
and
\begin{eqnarray}
\Phi_{N \delta_0}^0((s,x);z) &=& 
\frac{1}{p(s,x|0)} \sum_{\ell=0}^{\infty} z^{N-\ell-1}
\frac{1}{2 \pi \sqrt{-1}} \oint_{C(\delta_0)} d \zeta \,
\frac{p(s,x|\zeta)}{\zeta^{N-\ell}}
\nonumber\\
&=& 
\frac{1}{p(s,x|0)} \sum_{\ell=0}^{N-1} z^{N-\ell-1}
\frac{1}{2 \pi \sqrt{-1}} \oint_{C(\delta_0)} d \zeta \,
\frac{p(s,x|\zeta)}{\zeta^{N-\ell}},
\label{eqn:xi03}
\end{eqnarray}
since the integrands are holomorphic when $\ell \geq N$.

For BM with the transition probability density (\ref{eqn:pt1}),
(\ref{eqn:xi03}) gives
\begin{eqnarray}
\Phi_{N \delta_0}^0((s,x);z)
&=& \sum_{\ell=0}^{N-1} z^{N-\ell-1}
\frac{1}{2 \pi \sqrt{-1}}
\oint_{C(\delta_0)} d \zeta \,
\frac{e^{x\zeta/s-\zeta^2/2s}}{\zeta^{N-\ell}}
\nonumber\\
&=& \sum_{\ell=0}^{N-1} 
\left( \frac{z}{\sqrt{2s}} \right)^{N-\ell-1}
\frac{1}{2 \pi \sqrt{-1}}
\oint_{C(\delta_0)} d \eta \,
\frac{e^{2(x/\sqrt{2s}) \eta-\eta^2}}{\eta^{N-\ell}}
\nonumber\\
&=& \sum_{\ell=0}^{N-1} 
\left( \frac{z}{\sqrt{2s}} \right)^{N-\ell-1}
\frac{1}{(N-\ell-1)!} H_{N-\ell-1} \left(\frac{x}{\sqrt{2s}} \right),
\nonumber
\end{eqnarray}
where we have used the contour integral representation
of the Hermite polynomials (\ref{eqn:Hncontour}).
Thus its integral transformation is calculated as
\begin{eqnarray}
&& \cI \left[ \left.
\Phi_{N \delta_0}^0((s,x);W) \right| (t, y) \right]
\nonumber\\
&& \quad 
= \sum_{\ell=0}^{N-1} \frac{1}{(N-\ell-1)!} H_{N-\ell-1}
\left( \frac{x}{\sqrt{2s}} \right)
\frac{1}{(2s)^{(N-\ell-1)/2}}
\cI[W^{N-\ell-1}|(t,y)]
\nonumber\\
&& \quad 
= \sum_{\ell=0}^{N-1} \frac{1}{(N-\ell-1)!} H_{N-\ell-1}
\left( \frac{x}{\sqrt{2s}} \right)
\frac{1}{(2s)^{(N-\ell-1)/2}} 
m_{N-\ell-1}(t,y)
\nonumber\\
&& \quad 
= \sum_{\ell=0}^{N-1} \frac{1}{(N-\ell-1)! 2^{N-\ell-1}} 
\left( \frac{t}{s} \right)^{(N-\ell-1)/2}
H_{N-\ell-1}
\left( \frac{x}{\sqrt{2s}} \right)
H_{N-\ell-1}
\left( \frac{y}{\sqrt{2t}} \right),
\nonumber
\end{eqnarray}
where we have used Lemma \ref{thm:mn} and (\ref{eqn:mn3})
in Section \ref{sec:poly_mar}.
Then we obtain the following,
\begin{eqnarray}
&& \cM_{N \delta_0}^0((s,x)|(t,B(t)))
=\sum_{n=0}^{N-1} \frac{1}{n! 2^n} m_n(s,x) m_n(t,B(t))
\nonumber\\ 
&& \qquad \quad =
\sqrt{\pi} e^{x^2/4s+B(t)^2/4t}
\sum_{n=0}^{N-1} 
 \left(\frac{t}{s}\right)^{n/2}
\varphi_n \left( \frac{x}{\sqrt{2s}} \right)
\varphi_n \left( \frac{B(t)}{\sqrt{2t}} \right),
\label{eqn:Hermite3B}
\end{eqnarray}
where
\begin{equation}
\varphi_n(x)=\frac{1}{\sqrt{ \sqrt{\pi} 2^n n!}}
H_n(x) e^{-x^2/2}, 
\quad x \in \R,  \quad n \in \N_0, 
\label{eqn:phi1}
\end{equation}
are the 
\textbf{Hermite orthonormal functions} on $\R$,
\begin{equation}
\int_{\R} dx \, \varphi_n(x) \varphi_m(x)
= \delta_{n m}, \quad n, m \in \N_0.
\label{eqn:orthonormal1}
\end{equation}
The following expression for the transition probability density (\ref{eqn:pt1})
of BM is known as 
\textbf{Mehler's formula},
\begin{equation}
p(s-t, x|y)
=\frac{e^{-x^2/4s}}{e^{-y^2/4t}}
\frac{1}{\sqrt{2s}} \sum_{n=0}^{\infty} \left( \frac{t}{s} \right)^{n/2}
\varphi_n \left( \frac{x}{\sqrt{2s}} \right)
\varphi_n \left( \frac{y}{\sqrt{2t}} \right).
\label{eqn:Mehler}
\end{equation}

Since $m_n$, $n \in \N_0$ are 
the fundamental martingale polynomials
associated with BM, 
the process (\ref{eqn:Hermite3B}) is 
a continuous martingale.
Then we see that 
$$
\rE \left[ \cM_{N \delta_0}^0((s,x)|(t, \B(t))) \right]
= \rE \left[ \cM_{N \delta_0}^0((s,x)|(0, \B(0))) \right]=1
$$
for $ (s,x) \in [0, \infty) \times \R$, $0 \leq t < \infty$.

By the formula (\ref{eqn:KernelB1}), we obtain the
correlation kernels as
\begin{eqnarray}
\mbK_{N \delta_0}(s,x;t,y) &=& p(s,x|0) \cM_{N \delta_0}^0((s,x)| (t,y))
-\1(s>t) p(s-t, x|y)
\nonumber\\
&=& \frac{e^{-x^2/4s}}{e^{-y^2/4t}}
\bK_{\rm Hermite}^{(N)}(s,x;t,y)
\label{eqn:KKH}
\end{eqnarray}
with
\begin{eqnarray}
&& \bK_{\rm Hermite}^{(N)}(s,x;t,y)
= \frac{1}{\sqrt{2s}}
\sum_{n=0}^{N-1} \left( \frac{t}{s} \right)^{n/2}
\varphi_n \left( \frac{x}{\sqrt{2s}} \right)
\varphi_n \left( \frac{y}{\sqrt{2t}} \right)
\nonumber\\
&& \qquad \qquad \qquad - \ \1(s>t) 
\frac{1}{\sqrt{2s}}
\sum_{n=0}^{\infty} \left( \frac{t}{s} \right)^{n/2}
\varphi_n \left( \frac{x}{\sqrt{2s}} \right)
\varphi_n \left( \frac{y}{\sqrt{2t}} \right).
\nonumber\\
&& \qquad = \left\{ \begin{array}{ll}
\displaystyle{
\frac{1}{\sqrt{2s}}
\sum_{n=0}^{N-1} \left( \frac{t}{s} \right)^{n/2}
\varphi_n \left( \frac{x}{\sqrt{2s}} \right)
\varphi_n \left( \frac{y}{\sqrt{2t}} \right)
}
& \quad \mbox{for $s \leq t$},
\cr
\displaystyle{
-\frac{1}{\sqrt{2s}}
\sum_{n=N}^{\infty} \left( \frac{t}{s} \right)^{n/2}
\varphi_n \left( \frac{x}{\sqrt{2s}} \right)
\varphi_n \left( \frac{y}{\sqrt{2t}} \right)
}
& \quad \mbox{for $s > t$},
\end{array}
\right.
\label{eqn:K_Hermite}
\end{eqnarray}
where Mehler's formula (\ref{eqn:Mehler}) was used.
By the \textbf{gauge invariance} (Lemma \ref{thm:factor_K}),
the factor $e^{-x^2/4s}/e^{-y^2/4t}$ in (\ref{eqn:KKH})
is irrelevant for DSPs.
The kernel $\bK_{\rm Hermite}^{(N)}$ 
is known as the 
\textbf{extended Hermite kernel}
(see, for instance, Exercise 11.6.3 in \cite{For10}).

The equal-time correlation kernel
\begin{eqnarray}
\rK_{\rm Hermite}^{(N, \, t)}(x,y)
&\equiv&
\bK_{\rm Hermite}^{(N)}(t, x; t, y)
\nonumber\\
&=& \frac{1}{\sqrt{2t}} \sum_{n=0}^{N-1}
\varphi_n \left( \frac{x}{\sqrt{2t}} \right)
\varphi_n \left( \frac{y}{\sqrt{2t}} \right)
\nonumber
\end{eqnarray}
has the following expression, 
\begin{eqnarray}
&& K_{\rm Hermite}^{(N, \, t)}(x,y) = \sqrt{\frac{N}{2}} \,
\frac{\varphi_{N}(x/\sqrt{2t}) \varphi_{N-1}(y/\sqrt{2t})
- \varphi_{N-1}(x/\sqrt{2t}) \varphi_{N}(y/\sqrt{2t})}
{x-y}, \nonumber\\
\label{eqn:KNt3a}
&& \hskip 7cm
\mbox{if $x \not= y$},  \\
&& K_{\rm Hermite}^{(N, \, t)}(x,x) 
\nonumber\\
&& \quad = \frac{1}{\sqrt{2t}} \left[
N \left\{\varphi_{N} \left( \frac{x}{\sqrt{2t}} \right) \right\}^2
-\sqrt{N(N+1)} \varphi_{N-1} \left( \frac{x}{\sqrt{2t}} \right)
\varphi_{N+1}\left( \frac{x}{\sqrt{2t}} \right) \right].
\nonumber\\
\label{eqn:KNt3b}
\end{eqnarray}
This spatial correlation kernel is a special case
of the
\textbf{Christoffel--Dorboux kernel} 
(see, for instance, Chapter 9 in \cite{For10} and
Chapter 3 in \cite{AGZ10}).
It is called
the \textbf{Hermite kernel} and defines
the determinantal point process 
\cite{Sos00,ST03a} on $\R$
such that a spatial correlation function is
given by
\begin{equation}
\rho_{\rm Hermite}^{(N, \, t)}(\x_{N'})
=\det_{1 \leq i, j \leq N'} \Big[\rK_{\rm Hermite}^{(N, \, t)}(x_i, x_j) \Big]
\label{eqn:rho_Hermite}
\end{equation}
for any $1 \leq N' \leq N$ and
$\x_{N'}=(x_1, \dots, x_{N'}) \in \R^{N'}$,
$ t >0$.
We write the probability measure of this determinantal point
process as $\rP_{\rm Hermite}^{(N, \, t)}$.

\subsection{Exercises 1}
\label{sec:exercises1}
\subsubsection{Exercise 1.1}
\label{sec:ex1_1}
For $t >0$, $x \in \R$, and $\alpha \in \C$, let
\begin{equation}
G_{\alpha}(t, x)= e^{\alpha x - \alpha^2 t/2}.
\label{eqn:GB1}
\end{equation}
Let $(B(t))_{t \geq 0}$ be the 1 dimensional standard 
Brownian motion (BM) in 
$(\Omega, \rP, \cF)$ with the filtration
$\cF_t=\sigma(B(s) : 0 \leq s \leq t), t \geq 0$.
\begin{description}
\item{(1)} \,
Prove the equality
\begin{equation}
\frac{e^{\alpha B(t)}}{\rE[e^{\alpha B(t)}]}
=G_{\alpha}(t, B(t)), 
\quad t \geq 0.
\label{eqn:ex_G1}
\end{equation}

\item{(2)} \,
Prove that $(G_{\alpha}(t, B(t)))_{t \geq 0}$ is $\cF_t$-martingale;
\begin{equation}
\rE[ G_{\alpha}(t, B(t)) | \cF_s]
=G_{\alpha}(s, B(s)), \quad 0 < s < t.
\label{eqn:ex_G2}
\end{equation}

\item{(3)} \,
Let $\widehat{G}_{w}(t, x)$ be the Fourier transformation of 
$G_{\alpha}(t, x)$;
\begin{equation}
\widehat{G}_{w}(t, x) :=
\int_{-\infty}^{\infty} \frac{e^{-\sqrt{-1} \alpha w}}{2 \pi}
G_{\alpha}(t, x) d \alpha.
\end{equation}
Show that
\begin{equation}
\widehat{G}_{w}(t, x)
= \frac{e^{-(\sqrt{-1} x+ w)^2/2t}}{\sqrt{2 \pi t}}.
\label{eqn:ex_G4}
\end{equation}

\item{(4)} \,
Define the integral transformation by
\begin{equation}
\cI[f(W)|(t, x)]
= \int_{-\infty}^{\infty} f(\sqrt{-1} w) \widehat{G}_w(t,x) dw.
\label{eqn:ex_G5}
\end{equation}
For $n \in \N_0 := \{0,1,\dots\}$, let
\begin{align}
m_n(t,x) &= \cI[W^n|(t,x)]
\nonumber\\
&=\int_{-\infty}^{\infty} (\sqrt{-1} w)^n 
\widehat{G}_w(t,x) dw, 
\quad t \geq 0.
\label{eqn:ex_G6}
\end{align}
Express $m_n(t,x), n \in \N_0$ using the
Hermite polynomials,
\begin{align}
H_n(x) 
&:=(-1)^n e^{x^2} \frac{d^n e^{-x^2}}{d x^n}
\nonumber\\
&=
\sum_{k=0}^{[n/2]} (-1)^k \frac{n!}{k! (n-2k)!} (2x)^{n-2k}.
\label{eqn:ex_G7}
\end{align}

\item{(5)} \,
Prove the following for 
$\{m_n(t, x) \}_{n \in \N_0}$.
\begin{description}
\item{(i)} \, 
They are monic polynomials of degrees $n \in \N_0$
with time-dependent coefficients: 
\[
m_n(t, x)=x^n+\sum_{k=0}^{n-1} c_n^{(k)}(t) x^k, \quad t \geq 0.
\]
\item{(ii)} \,
For $0 \leq k \leq n-1$, $c_n^{(k)}(0)=0$. That is,
\[
m_n(0, x)=x^n, \quad n \in \N_0.
\]
\item{(iii)} \, If we set $x=B(t)$, they provide martingales: 
\begin{equation}
\rE[m_n(t, B(t)) |\cF_s]=m_n(s, B(s)), \quad 0 \leq s \leq t,
\quad n \in \N_0.
\label{eqn:ex_G8}
\end{equation}
\end{description}
\end{description}

\subsubsection{Exercise 1.2}
\label{sec:ex1_2}
Let
\begin{equation}
h_N(\x) := 
\prod_{1 \leq i < j \leq N} (x_j-x_i)
=\det_{1 \leq i, j \leq N} [x_j^{i-1}].
\label{eqn:ex_h1}
\end{equation}
\begin{description}
\item{(1)} \,
Using the basic properties of determinant,
verify the equalities,
\begin{eqnarray}
\frac{h_N(\y)}{h_N(\x)}
&=& \frac{1}{h_N(\x)} \det_{1 \leq i, j \leq N} [y_j^{i-1}]
\nonumber\\
&=& \frac{1}{h_N(\x)} \det_{1 \leq i, j \leq N}[ m_{i-1}(t, y_j)],
\label{eqn:ex_h2}
\end{eqnarray}
where $\{m_n(t, x)\}_{n \in \N_0}$ are given by
(\ref{eqn:ex_G6}) in Exercise 1.1.

\item{(2)} \,
We extend the integral transformation $\cI$
defined by (\ref{eqn:ex_G5}) in Exercise 1.1
to a linear integral transformation of multivariate functions
as follows: 
When $F^{(i)}(\x)=\prod_{j=1}^N f_j^{(i)}(x_j)$, $i=1,2$
are given for $\x=(x_1, \dots, x_N) \in \R^N$, then we define
\begin{equation}
\cI \left[ F^{(i)}(\bW) \left| \{(t_{\ell}, x_{\ell})\}_{\ell=1}^N  \right. \right]
:= \prod_{j=1}^N \cI \left[ \left. f^{(i)}_j(W_j) \right| (t_{j}, x_{j}) \right],
\quad i=1,2,
\label{eqn:ex_h3}
\end{equation}
and
\begin{align}
& \cI \Big[c_1 F^{(1)}(\bW) +c_2 F^{(2)}(\bW) \left| \{(t_{\ell}, x_{\ell})\}_{\ell=1}^N 
\right. \Big]
\nonumber\\
& \quad
:= c_1 \cI \Big[F^{(1)}(\bW) \left| \{(t_{\ell}, 
x_{\ell})\}_{\ell=1}^N \right. \Big]
+ c_2 \cI \Big[F^{(2)}(\bW) \left| \{(t_{\ell}, 
x_{\ell})\}_{\ell=1}^N \right. \Big],
\label{eqn:ex_h4}
\end{align}
$c_1, c_2 \in \C$, 
for $0 < t_i < \infty, i=1, \dots, N$, 
where $\bW=(W_1, \dots, W_N) \in \R^N$.
In particular, if $t_{\ell}=t, 1 \leq {^{\forall}\ell} \leq N$, we write
$\cI[\, \cdot \,| \{(t_{\ell}, x_{\ell})\}_{\ell=1}^N]$ simply as
$\cI[\, \cdot \,|(t, \x)]$ with $\x=(x_1, \dots, x_N)$.

Using multilinearity of determinants, verify the following 
equalities,
\begin{align}
\frac{h_N(\y)}{h_N(\x)}
&= \frac{1}{h_N(\x)}
\det_{1 \leq i, j \leq N}
\Big[ \cI [(W_j)^{i-1} | (t, y_j) ] \Big]
\nonumber\\
&= \cI \left[ \left. \frac{1}{h_N(\x)}
\det_{1 \leq i, j \leq N} [(W_j)^{i-1}] \right| (t, \y) \right]
\nonumber\\
&= \cI \left[ \left.
\frac{h_N(\bW)}{h_N(\x)} \right| (t, \y) \right].
\label{eqn:ex_h5}
\end{align}

\item{(3)} \,
Prove the equality,
\begin{equation}
\frac{h_N(\z)}{h_N(\x)}
=\det_{1 \leq i, j \leq N}
\Big[ \Phi_{\xi}^{x_i}(z_j) \Big],
\label{eqn:ex_h6}
\end{equation}
where
\begin{equation}
\Phi_{\xi}^{u}(z)=\prod_{\substack{1 \leq k \leq N, \cr x_k \not= u}}
\frac{z-x_k}{u-x_k}, 
\label{eqn:ex_h7}
\end{equation}
$\xi=\sum_{i=1}^N \delta_{x_i} \in \Conf_0(\R), z, u \in \C$,
following the steps given below.
\begin{description}
\item{(i)} \,
Let 
\begin{equation}
H(\x, \z)= \det_{1\leq i,j \leq N} 
\left[ \prod_{1\leq k \leq N,  k \not= i}(z_j-x_k) \right].
\label{eqn:ex_h8}
\end{equation}
Show that (a) it is a polynomial function of $x_i, z_i, 1 \leq i \leq N$
with degree $N(N-1)$, 
(b) $H(\x, \x)=h(\x)^2$, and 
(c) $H(\x, \z)=0$, 
if $x_i=x_j$ or $z_i=z_j$ for some $i, j$ with $1\leq i<j\leq N$.
It is implied by (a)--(c) that
$H(\x, \z)=h_N(\x) h_N(\z)$.

\item{(ii)} \,
By the definition (\ref{eqn:ex_h7}), hence 
the RHS of (\ref{eqn:ex_h6})
is equal to 
$H(\x, \z)/h(\x)^2$
for $\xi \in \Conf_0(\R)$.
\end{description}
\end{description}

\subsubsection{Exercise 1.3}
\label{sec:ex1_3}
For a pair of space-time coordinates
$(s, x), (t, y) \in [0, \infty) \times \R$, consider
the function,
\begin{equation}
\cM_N((s,x)|(t, y))
=\sum_{n=0}^{N-1} \frac{1}{n! 2^n} m_n(s,x) m_n(t, y),
\label{eqn:ex_K1}
\end{equation}
where $\{m_n(t, x)\}_{n \in \N_0}$ are given by
(\ref{eqn:ex_G6}) in Exercise 1.1.
\begin{description}
\item{(1)} \,
The 
\textbf{Hermite orthonormal functions} are defined 
by
\begin{equation}
\varphi_n(x)=\frac{1}{\sqrt{ \sqrt{\pi} 2^n n!}}
H_n(x) e^{-x^2/2}, 
\quad x \in \R,  \quad n \in \N_0, 
\label{eqn:ex_K2}
\end{equation}
where $\{H_n(x)\}_{n \in \N_0}$ are 
the Hermite polynomials defined by
(\ref{eqn:ex_G7}). 
Prove the following expression of $\cM_N((s,x)|(t, y)$,
\begin{equation}
\cM_N((s,x)|(t, y))
=\sqrt{\pi} e^{x^2/4s+y^2/4t}
\sum_{n=0}^{N-1} 
 \left(\frac{t}{s}\right)^{n/2}
\varphi_n \left( \frac{x}{\sqrt{2s}} \right)
\varphi_n \left( \frac{y}{\sqrt{2t}} \right).
\label{eqn:ex_K3}
\end{equation}

\item{(2)} \,
Let $(B(t))_{t \geq 0}$ be BM
adapted to the filtration
$\cF_t, t \geq 0$.
Then prove that, for fixed $(s, x) \in [0, \infty) \times \R$, 
$(\cM_N((s,x)|(t, B(t))))_{t \geq 0}$ is 
an $\cF_t$-martingale, and 
\begin{equation}
\rE \left[ \cM_{N}((s,x)|(t, \B(t))) \right]
= \rE \left[ \cM_{N}((s,x)|(0, \B(0))) \right]=1.
\label{eqn:ex_K4}
\end{equation}

\item{(3)} \,
The following identity is known as
\textbf{Mehler's formula},
\begin{equation}
\sum_{n=0}^{\infty} 
\frac{H_n(x) H_n(y)}{2^n n!} a^n
=\frac{e^{[2xy a-(x^2+y^2)a^2]/(1-a^2)}}{\sqrt{1-a^2}},
\quad |a| < 1.
\label{eqn:ex_K5}
\end{equation}
Prove the following limit,
\begin{equation}
\lim_{N \to \infty}
\cM_N((s,x)|(t,y))
=\frac{p(s-t, x|y)}{p(s, x|0)},
\quad 0 < t < s, \quad x, y \in \R,
\label{eqn:ex_K6}
\end{equation}
where $p(t, y|x)$ denotes the
transition probability density of BM,
\begin{equation}
p(t, y|x)=\frac{e^{-(y-x)^2/2t}}{\sqrt{2 \pi t}},
\quad t >0, \quad x, y \in \R.
\label{eqn:ex_K7}
\end{equation}

\item{(4)} \,
Define the spatio-temporal kernel as follows,
\begin{equation}
K(s,x; t, y)
:= p(s, x|0) \cM_N((s,x)|(t,y))
-{\bf 1}(s > t) p(s-t, x|y),
\label{eqn:ex_K8}
\end{equation}
$(s,x), (t, y) \in [0, \infty) \times \R$,
where ${\bf 1}(\omega)$ is an indicator;
${\bf 1}(\omega)=1$, if $\omega$ is satisfied,
and ${\bf 1}(\omega)=0$, otherwise.
Then prove the equality
\begin{equation}
\mbK_{N \delta_0}(s,x;t,y)= \frac{e^{-x^2/4s}}{e^{-y^2/4t}}
\bK_{\rm Hermite}^{(N)}(s,x;t,y),
\label{eqn:e_K9}
\end{equation}
where
\begin{align}
& \bK_{\rm Hermite}^{(N)}(s,x;t,y)
\nonumber\\
& \qquad = \left\{ \begin{array}{ll}
\displaystyle{
\frac{1}{\sqrt{2s}}
\sum_{n=0}^{N-1} \left( \frac{t}{s} \right)^{n/2}
\varphi_n \left( \frac{x}{\sqrt{2s}} \right)
\varphi_n \left( \frac{y}{\sqrt{2t}} \right)
}
& \quad \mbox{for $s \leq t$},
\cr
\displaystyle{
-\frac{1}{\sqrt{2s}}
\sum_{n=N}^{\infty} \left( \frac{t}{s} \right)^{n/2}
\varphi_n \left( \frac{x}{\sqrt{2s}} \right)
\varphi_n \left( \frac{y}{\sqrt{2t}} \right)
}
& \quad \mbox{for $s > t$}.
\end{array}
\right.
\label{eqn:ex_K10}
\end{align}
The spatio-temporal kernel 
$\bK_{\rm Hermite}^{(N)}(s,x;t,y)$
is known as the \textbf{extended Hermite kernel}
giving the correlation kernel 
for DYS$_2$ starting from 
the configuration such that
all $N$ particles are at the origin,
$N \delta_0 \in \Conf(\R) \setminus \Conf_0(\R)$.
\end{description}

\SSC
{Gaussian Analytic Functions and Their Zero-Point Processes}
\label{sec:GAF}
\subsection{Reproducing kernel Hilbert spaces}
\label{sec:reproducing}
\subsubsection{Bergman spaces and Hardy spaces}
\label{sec:examples}
A \textbf{Hilbert function space of holomorphic functions}
is a Hilbert space $\cH$ of functions
on a domain $D$ in $\C$ 
equipped with the inner product $\bra \cdot, \cdot \ket_{\cH}$,
\[
\begin{array}{ll}
\bra af+bg, h \ket_{\cH}= a \bra f, h \ket_{\cH}+b \bra g, h \ket_{\cH},
& \cr
\bra h, a f + b g \ket_{\cH}
=\overline{a} \bra h, f \ket_{\cH} 
+\overline{b} \bra h, g \ket_{\cH}, 
&
\quad f, g, h \in \cH, \quad a, b \in \C
\cr
\bra g, f \ket_{\cH}=\overline{ \bra f, g \ket_{\cH}}
=\bra \overline{f}, \overline{g} \ket_{\cH},
&
\end{array} 
\]
and the norm $\|f\|_{\cH} :=\sqrt{\bra f, f \ket_{\cH}}$.
The evaluation at each point 
of $D$ is a continuous functional on $\cH$.
Therefore, for each point $w \in D$, 
there is an element of $\cH$,
which is called the 
\textbf{reproducing kernel at $w$}
and denote by $k_{w}$, with the property
\[
\bra f, k_w \ket_{\cH}=f(w), \quad \forall f \in \cH.
\]
Because $k_{w} \in \cH$, it is itself a function on $D$,
\[
k_w(z)=\bra k_w, k_z \ket_{\cH}.
\]
We write
\begin{equation}
k_{\cH}(z, w):=k_w(z) = \bra k_w, k_z \ket_{\cH}
\label{eqn:reproducing_kernel}
\end{equation}
and call it the \textbf{reproducing kernel} for $\cH$.
By definition, it is \textbf{Hermitian}; 
\[
\overline{k_{\cH}(z,w)}=k_{\cH}(w, z), \quad z, w \in D.
\]
If $\cH$ is a holomorphic Hilbert function space, then
$k_{\cH}$ is holomorphic in the first variable and
anti-holomorphic in the second.
We see that $k_{\cH}(z,w)$ is a 
\textbf{positive semi-definite kernel}: 
for any $n \in \N:=\{1,2,\dots\}$, 
for any points $z_i \in D$ and $\xi_i \in \C$, $i=1,2, \dots, n$,
\begin{align}
\sum_{i=1}^n \sum_{j=1}^n k_{\cH}(z_i, z_j) \xi_i \overline{\xi_j} 
&= \sum_{i=1}^n \sum_{j=1}^n 
\bra k_{z_j}, k_{z_i} \ket_{\cH} \xi_i \overline{\xi_j} 
=\sum_{i=1}^n \sum_{j=1}^n 
\bra k_{\cH}(\cdot, z_j), k_{\cH}(\cdot, z_i) 
\ket_{\cH} \xi_i \overline{\xi_j}
\nonumber\\
&=  \left\bra \sum_{i=1}^n k_{\cH}(z_i, \cdot) \xi_i,
\sum_{j=1}^n k_{\cH}(z_j, \cdot) \xi_j \right\ket_{\cH}
=\Big\|\sum_{i=1}^n \xi_i k_{\cH}(z_i, \cdot) 
\Big\|_{\cH}^2 \ge 0. 
\label{eqn:positive_kernel}
\end{align}
Let $\{e_n : n \in \cI \}$ be any 
\textbf{complete orthonormal systme} 
(CONS) for $\cH$,
where $\cI$ is an index set,
\begin{align*}
&\bra e_n, e_m \ket_{\cH}=\delta_{nm}, 
\, n, m \in \cI,
\nonumber\\
&f \in \cH \iff 
f=\sum_{n \in \cI} c_n e_n \, \, \mbox{with} \, \, \,
(c_n)_{n \in \cI} \in \ell^2(\cI).
\end{align*}
Then one can prove that
the reproducing kernel for $\cH$ is written in the form
\begin{equation}
k_{\cH}(z, w) = \sum_{n \in \cI} e_n(z) \overline{e_n(w)}.
\label{eqn:reproducing}
\end{equation}
Actually, this gives
\begin{align*}
\bra f(\cdot), k_{\cH}(\cdot, w) \ket_{\cH}
&=\left\bra \sum_{m \in \cI} c_m e_m(\cdot), 
\sum_{n \in \cI} e_n(\cdot) \overline{e_n(w)} \right\ket_{\cH}
=\sum_{m \in \cI} \sum_{n \in \cI} c_m \bra e_m, e_n \ket_{\cH} e_n(w)
\nonumber\\
&=\sum_{n \in \cI} c_n e_n(w) = f(w), \, \, \, \forall f \in \cH, \, w \in D. 
\end{align*}
We note that the positive definiteness of the kernel 
(\ref{eqn:positive_kernel})
is equivalent with the situation such that, 
for any points $z_i \in D$, $i \in \N$, the matrix 
$(k_{\cH}(z_i, z_k))_{1 \leq i,j \leq n}$ has a nonnegative determinant, 
\[
\det_{1\leq i, j \leq n}[k_{\cH}(z_i, z_j)] \geq 0
\quad \mbox{for any $n \in \N$}.
\]

First we show two examples
of holomorphic Hilbert function spaces,
the \textbf{Bergman space} and the \textbf{Hardy space}, 
for a \textbf{unit disk} $\D:=\{z \in \C: |z| < 1\}$ 
and the domains which are conformally 
transformed from $\D$.

\begin{example}
\label{thm:ex_Bargman}
The \textbf{Bergman space} on $\D$, denoted by $L^2_{\rB}(\D)$, is
the Hilbert space of holomorphic functions on $\D$
which are square-integrable with respect to
the Lebesgue measure on $\C$.
The inner product for $L^2_{\rB}(\D)$ is given by
\[
\bra f, g \ket_{L^2_{\rB}(\D)}
:= \frac{1}{\pi} \int_{\D} f(z) \overline{g(z)} m(dz)
= \sum_{n=0}^{\infty} \frac{\widehat{f}(n) 
\overline{\widehat{g}(n)}}{n+1},
\]
where the $n$th Taylor coefficient of $f$ at 0 is denoted by
$\widehat{f}(n)$;
$f(z) = \sum_{n=0}^{\infty} \widehat{f}(n) z^n$.
Let 
\[
\widetilde{e_n}(z)
:=\sqrt{n+1} z^{n}, n \in \N_0.
\]
Then $\{\widetilde{e_n}(z) \}_{n \in \N_0}$ form a
CONS for $L^2_{\rB}(\D)$
and the reproducing kernel (\ref{eqn:reproducing}) is given by
\begin{align}
K_{\D}(z, w) 
&:=k_{L^2_{\rB}(\D)}(z, w)
\nonumber\\
&= \sum_{n \in \N_0}
(n+1) (z \overline{w})^n
=\frac{1}{(1-z \overline{w})^2},
\quad z, w \in \D.
\label{eqn:K_D}
\end{align}
This kernel is called the \textbf{Bergman kernel} of $\D$.
\end{example}

\begin{example}
\label{thm:ex_Hardy}
The \textbf{Hardy space} on $\D$, $H^2(\D)$, consists of 
holomorphic functions on $\D$ such that the Taylor coefficients
form a square-summable series;
\[
\| f \|_{H^2(\D)}^2
:= \sum_{n \in \N_0} |\widehat{f}(n)|^2 < \infty,
\quad f \in H^2(\D).
\]
For every $f \in H^2(\D)$, 
the non-tangential limit
$\lim_{r \uparrow 1} f(r e^{\sqrt{-1} \phi})$ exists
a.e.~by Fatou's theorem and we write it as
$f(e^{\sqrt{-1} \phi})$.
It is known that $f(e^{\sqrt{-1} \phi}) \in L^2(\partial \D)$. 
Then one can prove that the inner product of $H^2(\D)$ 
is given by the following
three different ways,
\begin{equation}
\bra f, g \ket_{H^2(\D)}
= \begin{cases}
\displaystyle{
\sum_{n \in \N_0} \widehat{f}(n) \overline{\widehat{g}(n)}},
\cr
\displaystyle{
\lim_{r \uparrow 1} \frac{1}{2\pi}
\int_{0}^{2 \pi} f(r e^{\sqrt{-1} \phi})
\overline{g(r e^{\sqrt{-1} \phi})} d \phi},
\quad f, g \in H^2(\D), 
\cr
\displaystyle{
\frac{1}{2\pi}
\int_{0}^{2 \pi} f(e^{\sqrt{-1} \phi})
\overline{g(e^{\sqrt{-1} \phi})} d \phi},
\end{cases}
\label{eqn:IP_H2D}
\end{equation}
with $\|f\|_{H^2(\D)}^2=\bra f, f \ket_{H^2(\D)}$.
Let $\sigma$ be the measure on the boundary
of $\D$ which is the usual arc length measure.
Then the last expression of the inner product (\ref{eqn:IP_H2D})
is written as
$\bra f, g \ket_{H^2(\D)}
=(1/2 \pi) \int_{\gamma_1} f(z) \overline{g(z)} \sigma(dz)$,
where $\gamma_1$ is a unit circle 
$\{ e^{\sqrt{-1} \phi} : \phi \in [0, 2\pi) \}$ giving the boundary 
of $\D$.
If we set $e_n(z) :=e^{(0,0)}_n(z)=z^n, n \in \N_0$,
then $\{e_n(z)\}_{n \in \N_0}$ form
CONS for $H^2(\D)$. The reproducing kernel (\ref{eqn:reproducing})
is given by
\begin{align}
S_{\D}(z, w) 
&:= k_{H^2(\D)}(z, w)
\nonumber\\
&=
\sum_{n \in \N_0} (z \overline{w})^n
=\frac{1}{1-z \overline{w}},
\quad z, w \in \D, 
\label{eqn:S_D}
\end{align}
which is called the \textbf{Szeg\H{o} kernel} of $\D$.
\end{example}

Let $f: D \to \widetilde{D}$ be a 
\textbf{conformal transformation} 
(the angle-preserving one-to-one map) 
between two bounded domains $D, \widetilde{D} \subsetneq \C$
with $\cC^{\infty}$ smooth boundary (analytic boundary).
We find an argument in Chapter 12 of \cite{Bel16} 
concluding that 
the derivative of the transformation $f$ denoted by
$f'$ has a single valued square root on $D$.
We let $\sqrt{f'(z)}$ denote one of the 
square roots of $f'$.
The Szeg\H{o} kernel and the Bergman kernel 
are then transformed by $f$ as 
\cite[Chapters 12, 16]{Bel16}
\begin{align}
S_D(z, w) &= 
\sqrt{f'(z)} 
\overline{\sqrt{f'(w)}}
S_{\widetilde{D}}(f(z), f(w)),
\nonumber\\
K_D(z, w) &= |f'(z)|
|f'(w)| K_{\widetilde{D}}(f(z), f(w)), 
\quad z, w \in D.
\label{eqn:SD_KD_conformal}
\end{align}
Consider the special case in which
$D \subsetneq \C$ is a simply connected domain
with analytic boundary 
and $\widetilde{D}=\D$.
For each $\alpha \in D$, 
\textbf{Riemann's mapping theorem} 
gives a unique conformal transformation;
\begin{equation}
\mbox{
$h_{\alpha} : D \to \D$ \quad conformal 
such that \, 
$h_{\alpha}(\alpha)=0, \, \, h_{\alpha}'(\alpha) >0$.
}
\label{eqn:RiemannA}
\end{equation}
Such $h_{\alpha}$ is called the
\textbf{Riemann mapping function}. 
By (\ref{eqn:S_D}), 
the first equation in 
(\ref{eqn:SD_KD_conformal}) gives 
\begin{align*}
S_D(z, w) &= 
\sqrt{h_{\alpha}'(z)} 
\overline{\sqrt{h_{\alpha}'(w)}}
S_{\D}(h_{\alpha}(z), h_{\alpha}(w))
\nonumber\\
&=\sqrt{h_{\alpha}'(z)} 
\overline{\sqrt{h_{\alpha}'(w)}}
\frac{1}{1-h_{\alpha}(z) \overline{h_{\alpha}(w)}},
\quad z, w \in D.
\end{align*}
Since $h_{\alpha}(\alpha)=0$, if we put
$z=\alpha$, or $w=\alpha$, or
$z=w=\alpha$ in this formula, 
then we readily obtain
\begin{align*}
S_D(\alpha, w) &= \sqrt{h_{\alpha}'(\alpha)} 
\overline{\sqrt{h_{\alpha}'(w)}}
\frac{1}{1-0 \times \overline{h_{\alpha}(w)}}
= \sqrt{h_{\alpha}'(\alpha)} 
\overline{\sqrt{h_{\alpha}'(w)}},
\nonumber\\
S_D(z, \alpha) &=\sqrt{h_{\alpha}'(z)} 
\overline{\sqrt{h_{\alpha}'(\alpha)}},
\quad
S_D(\alpha, \alpha) =\sqrt{h_{\alpha}'(\alpha)} 
\overline{\sqrt{h_{\alpha}'(\alpha)}}
\nonumber\\
\Longrightarrow \quad
& \frac{S_D(z, \alpha) S_D(\alpha, w)}{S_D(\alpha, \alpha)}
=\frac{\sqrt{h_{\alpha}'(\alpha)} 
\overline{\sqrt{h_{\alpha}'(w)}}
\sqrt{h_{\alpha}'(z)} 
\overline{\sqrt{h_{\alpha}'(\alpha)}} }
{\sqrt{h_{\alpha}'(\alpha)} 
\overline{\sqrt{h_{\alpha}'(\alpha)}} }
=\sqrt{h_{\alpha}'(z)} 
\overline{\sqrt{h_{\alpha}'(w)}}
\end{align*}
Hence the above formula is written as
\cite{Bel95}, 
\begin{equation}
S_D(z, w)= \frac{S_D(z, \alpha) \overline{S_D(w, \alpha)}}
{S_D(\alpha, \alpha)} 
\frac{1}{1-h_{\alpha}(z) \overline{h_{\alpha}(w)}},
\quad z, w, \alpha \in D.
\label{eqn:SD_formula}
\end{equation}
Similarly, we have
\begin{equation}
K_D(z, w)= 
\frac{S_D(z, \alpha)^2 \overline{S_D(w, \alpha)^2}}
{S_D(\alpha, \alpha)^2}
\frac{1}{(1-h_{\alpha}(z) \overline{h_{\alpha}(w)})^2},
\quad z, w, \alpha \in D.
\label{eqn:KD_formula}
\end{equation}
Hence the following relationship is established, 
\begin{align}
S_{D}(z, w)^2 &= K_{D}(z, w),
\quad z, w \in D.
\label{eqn:S_K_D}
\end{align}

Here we notice that (\ref{eqn:SD_formula}) is rewritten as
\begin{equation}
S_D(z, w)- \frac{S_D(z, \alpha) S_D(\alpha, w)}{S_D(\alpha, \alpha)}
=S_D(z,w) h_{\alpha}(z) \overline{h_{\alpha}(w)},
\quad z, w, \alpha \in D.
\label{eqn:condition}
\end{equation}
LHS is regarded as the reproducing
kernel for the Hilbert subspace
\begin{equation}
H^2_a(D) := \{ f \in H^2(D) : f(\alpha)=0 \},
\label{eqn:H2a}
\end{equation}
and we will denote such a
\textbf{conditional kernel} with a deterministic zero
at $\alpha$ by $S_D^{\alpha}(z, w)$.
Hence we obtain the formula
\begin{equation}
S_D^{\alpha}(z, w)=S_D(z,w) h_{\alpha}(z) \overline{h_{\alpha}(w)},
\quad z, w, \alpha \in D.
\label{eqn:SaD}
\end{equation}

Let $q \in (0, 1)$ be a fixed number and 
we consider the \textbf{annulus}
\[
\A_q:=\{z \in \C: q < |z| < 1 \}. 
\]
As the third example of the 
Hilbert spaces of holomorphic functions,
we consider the Hardy space on $\A_q$.

\begin{example}
\label{thm:Hardy_Aq}
The \textbf{Hardy space for $\A_q$},
$H^2(\A_q)$, is the Hilbert space of holomorphic functions 
on $\A_q$ 
equipped with the inner product
\begin{equation}
\bra f, g \ket_{H^2(\A_q)} =
\frac{1}{2\pi}
\int_{0}^{2 \pi} f(e^{\sqrt{-1} \phi})
\overline{g(e^{\sqrt{-1} \phi})} d \phi
+ \frac{1}{2\pi} \int_{0}^{2 \pi} f(q e^{\sqrt{-1} \phi})
\overline{g(q e^{\sqrt{-1} \phi})} \, q d \phi,
\, \, f, g \in H^2(\A_q).
\label{eqn:e_qq}
\end{equation}
A CONS of $H^2(\A_q)$ 
is given by $\{e^{(q, q)}_n\}_{n \in \Z}$ with
\[
e^{(q, q)}_n(z) = \frac{z^n}{\sqrt{1+q^{2n+1}}},
\quad z \in \A_q, \quad n \in \Z,
\]
and the reproducing kernel is given by \cite{MS94}
\begin{equation}
S_{\A_q}(z, w)
= \sum_{n \in \Z} e^{(q, q)}_n(z) \overline{e^{(q, q)}_n(w)}
= \sum_{n=-\infty}^{\infty} 
\frac{(z \wbar)^n}{1+q^{2n+1}}. 
\label{eqn:SAq1}
\end{equation}
This infinite series converges absolutely for $z,w \in \A_q$
and 
is called the
\textbf{Szeg\H{o} kernel of $\A_q$}.
\end{example}

The nonexistence of zero in $\D$ 
of $S_{\D}(\cdot, \alpha), \alpha \in \D$ and the
uniqueness of zero in $\A_q$ of $S_{\A_q}(\cdot, \alpha)$, 
$\alpha \in \A_q$ are concluded from a general consideration
(see, for instance, \cite[Section 27]{Bel16}).
Define
\begin{equation}
\alphahat := - \frac{q}{\alphabar}, \quad \alpha \in \A_q.
\label{eqn:alphabar}
\end{equation}
The fact
\begin{equation}
S_{\A_q}(\alphahat, \alpha)=S_{\A_q}(\alpha, \alphahat)
=0, \quad
\alpha \in \A_q
\label{eqn:zeros_S}
\end{equation}
was proved 
as Theorem 1 in \cite{TT99} by direct calculation, 
for which a simpler proof can be given 
using theta functions (see Lemma \ref{thm:zero_S} below). 

Again we consider the conformal transformation 
$f: D \to \widetilde{D}$ between two bounded domains $D, 
\widetilde{D} \subsetneq \C$
with analytic boundaries.
Here we consider the case in which
$D$ is a {\bf 2-connected domain} and $\widetilde{D}=\D$.
For each $\alpha \in D$, 
there exists a unique function
$f_a$ giving 
a {\bf branched 2 to 1 covering map} of $D$ to $\D$, 
in which a unique point 
$\alphahat=\alphahat(\alpha) \not= \alpha$, 
$\alphahat \in D$ exists
such that 
\[
f_{\alpha}(\alpha)=f_{\alpha}(\alphahat)=0.
\]
The map $f_{\alpha}, \alpha \in D$
is called the \textbf{Ahlfors mapping function}
\cite[Chapter 13]{Bel16}.
By the fact (\ref{eqn:zeros_S}) with (\ref{eqn:alphabar}),
for an arbitrary, but fixed $\alpha \in D$, 
we can verify that the following set forms
a CONS for $H^2(\A_q)$ which is different
from $\{e^{(q,q)}_n(z)\}_{n \in \Z}$ given by (\ref{eqn:e_qq}),
\begin{equation}
\widehat{e}_{jn}(z)
:=\frac{S_{\A_q}(z, \alpha_j)}{\sqrt{S_{\A_q}(\alpha_j, \alpha_j)}}
f_{\alpha}(z)^n,
\quad j=0,1, \quad n \in \N_0,
\label{eqn:hat_e}
\end{equation}
where
\[
\alpha_0:=\alpha, \quad
\alpha_1:= \alphahat.
\]
Hence the Szeg\H{o} kernel of $\A_q$ can be expressed as
\begin{align}
S_{\A_q}(z, w)
&=\sum_{j=0}^{1} \sum_{n=0}^{\infty}
\widehat{e}_{j n}(z) \overline{\widehat{e}_{j n}(w)}
\nonumber\\
& 
=\left(
\frac{S_{\A_q}(z, \alpha) 
\overline{S_{\A_q}(w, \alpha)}}{S_{\A_q}(\alpha, \alpha)}
+ \frac{S_{\A_q}(z, \alphahat) 
\overline{S_{\A_q}(w, \alphahat)}}{S_{\A_q}(\alphahat, \alphahat)}
\right) \sum_{n=0}^{\infty} f_{\alpha}(z)^n \overline{f_{\alpha}(z)^n}
\nonumber\\
&
=\left(
\frac{S_{\A_q}(z, \alpha) 
\overline{S_{\A_q}(w, \alpha)}}{S_{\A_q}(\alpha, \alpha)}
+ \frac{S_{\A_q}(z, \alphahat) 
\overline{S_{\A_q}(w, \alphahat)}}{S_{\A_q}(\alphahat, \alphahat)}
\right) \frac{1}{1-f_{\alpha}(z) \overline{f_{\alpha}(w)}},
\, z, w \in \A_q.
\label{eqn:S_new}
\end{align}
Then, if we consider the Hilbert subspace of
$\A_q$ such that
\[
H^2_{\alpha, \alphahat}(\A_q)
:=\{f \in H^2(\A_q) : f(\alpha)=f(\alphahat)=0\},
\]
then its reproducing kernel should be given by
\begin{equation}
S_{\A_q}^{\alpha, \alphahat}(z, w)
:=S_{\A_q}(z, w)
-\left(
\frac{S_{\A_q}(z, \alpha) 
S_{\A_q}(\alpha, w)}{S_{\A_q}(\alpha, \alpha)}
+ \frac{S_{\A_q}(z, \alphahat) 
S_{\A_q}(\alphahat, w)}{S_{\A_q}(\alphahat, \alphahat)}
\right).
\label{eqn:SAq_a_ahat}
\end{equation}
It is easy to verify that (\ref{eqn:S_new}) gives 
\cite{Bel16}
\begin{equation}
S_{\A_q}^{\alpha, \alphahat}(z, w)
=S_{\A_q}(z, w) f_{\alpha}(z) \overline{f_{\alpha}(w)},
\quad z, w \in \A_q.
\label{eqn:SAq_a_ahat2}
\end{equation}

Now we ask the following question.
\begin{description}
\item{\bf (Q)} \,
For $\alpha \in \A_q$, consider the
conditional Hilbert space 
\[
H^2_{\alpha}(\A_q)
:= \{ f \in H^2(\A_q):  f(\alpha)=0\},
\]
whose reproducing kernel is given by
\begin{equation}
S_{\A_q}^{\alpha}(z, w)
:= S_{\A_q}(z, w) -
\frac{S_{\A_q}(z, \alpha) S_{\A_q}(\alpha, w)}
{S_{\A_q}(\alpha, \alpha)}, \quad
z, w \in \A_q.
\label{eqn:conditionSat}
\end{equation}
Is it possible to factorize this 
\textbf{conditional kernel} in the form
as (\ref{eqn:SaD}) and (\ref{eqn:SAq_a_ahat2})?
\end{description}

\subsubsection{Theta function}
\label{sec:theta}

Assume that $p \in \C$ is a fixed number such that 
$0 < |p| < 1$.
We use the following standard notation,
\begin{align}
&(a; p)_{n} := \prod_{i=0}^{n-1} (1-a p^i), \qquad
(a; p)_{\infty} := \prod_{i=0}^{\infty} (1-a p^i),
\nonumber\\
&(a_1, \dots, a_k; p)_{\infty}
: =(a_1; p)_{\infty} \cdots (a_k; p)_{\infty}.
\label{eqn:Pochhammer}
\end{align}
The \textbf{theta function} with argument $z$ and
\textbf{nome} $p$ is defined by
\begin{equation}
\theta(z; p)  :=(z, p/z; p)_{\infty}.
\label{eqn:theta}
\end{equation}
We often use the shorthand notation
$\theta(z_1, \dots, z_n; p)
: = \prod_{i=1}^n \theta(z_i; p)$.

As a function of $z$, the theta function $\theta(z; p)$ is
holomorphic in $\C^{\times}$ and has single zeros
precisely at $p^{i}$, $i \in \Z$, that is,
\begin{equation}
\{ z \in \C^{\times} : \theta(z; p)=0 \}
=\{p^i : i \in \Z \}.
\label{eqn:theta_zero}
\end{equation}
We will use the \textbf{inversion formula} 
\begin{equation}
\theta(1/z; p) = - \frac{1}{z} \theta(z; p)
\label{eqn:theta_inversion}
\end{equation}
and the \textbf{quasi-periodicity property}
\begin{equation}
\theta(pz; p) = -\frac{1}{z} \theta(z; p)
\label{eqn:theta_qp}
\end{equation}
of the theta function.
By comparing (\ref{eqn:theta_inversion}) and 
(\ref{eqn:theta_qp}) and performing 
the transformation $z \mapsto 1/z$, 
we immediately see the \textbf{periodicity property}, 
\begin{align}
\theta(p/z; p) 
= \theta(z; p).
\label{eqn:theta_period}
\end{align}
By \textbf{Jacobi's triple product identity}
(see, for instance, \cite[Section 1.6]{GR04}), 
we have 
the Laurent expansion
\[
\theta(z; p)=\frac{1}{(p; p)_{\infty}}
\sum_{n \in \Z} (-1)^n p^{\binom{n}{2}} z^n.
\]
One can show that \cite[Chapter 20]{NIST10}
\begin{align}
&\lim_{p \to 0} \theta(z; p)=1-z,
\label{eqn:theta_p0}
\\
&\theta'(1; p) 
:= \frac{\partial \theta(z; p)}{\partial z} \Big|_{z=1}
=- (p; p)_{\infty}^2.
\label{eqn:theta_prime}
\end{align}
The theta function satisfies the following 
\textbf{Weierstrass' addition formula} \cite{Koo14},
\begin{equation}
\theta(xy, x/y, uv, u/v; p)
- \theta(xv, x/v, uy, u/y; p)
=\frac{u}{y} \theta(yv, y/v, xu, x/u; p).
\label{eqn:Weierstrass_add1}
\end{equation}

When $p$ is real and $p \in (0, 1)$, we see that
\begin{equation}
\overline{\theta(z; p)}= \theta( \zbar; p).
\label{eqn:theta_real}
\end{equation}
In this case the definition (\ref{eqn:theta})
with (\ref{eqn:Pochhammer}) implies that
\begin{align}
&\left.
\begin{array}{ll}
\theta(x; p) > 0, \, & x \in (p^{2i+1}, p^{2i})
\cr
\theta(x; p) =0, \, & x = p^{i}
\cr
\theta(x; p) < 0, \, & x \in (p^{2i}, p^{2i-1})
\end{array}
\right\}
\quad i \in \Z,
\nonumber\\
&\quad \theta(x; p) >0, \quad x \in (-\infty, 0).
\label{eqn:theta_signs}
\end{align}
Moreover, we can prove the following: 
In the interval $x \in (-\infty, 0)$,
$\theta(x) :=\theta(x; p)$ is strictly convex with
\begin{equation}
\min_{x \in (-\infty, 0)} \theta(x)
=\theta(-\sqrt{p}) =
\prod_{n=1}^{\infty}(1+p^{n-1/2})^2 > 0,
\label{eqn:min_theta}
\end{equation}
and 
$\lim_{x \downarrow -\infty} \theta(x)
=\lim_{x \uparrow 0} \theta(x)=+\infty$,
and in 
the interval $x \in (p, 1)$,
$\theta(x)$ is strictly concave with
\begin{equation}
\max_{x \in (p, 1)} \theta(x)
=\theta(\sqrt{p}) =
\prod_{n=1}^{\infty}(1-p^{n-1/2})^2,
\label{eqn:max_theta}
\end{equation}
$\theta(x) \sim (p; p)_{\infty}^2 (x-p)/p$
as $x \downarrow p$, and
$\theta(x) \sim (p; p)_{\infty}^2 (1-x)$
as $x \uparrow 1$,
where (\ref{eqn:theta_qp}) and 
(\ref{eqn:theta_prime}) were used.

\subsubsection{Weighted Szeg\H{o} kernels of annulus 
and conformal transformations}
\label{sec:Ahlfors}
We introduce a positive parameter $r>0$.
Consider the Hilbert space of holomorphic functions on $\A_q$ 
equipped with the inner product
\[
\bra f, g \ket_{H^2_r(\A_q)} = 
\frac{1}{2\pi} \int_{\gamma_1 \cup \gamma_q}
f(z) \overline{g(z)} \sigma_r(dz),
\quad f, g \in H^2_r(\A_q)
\]
with
\[
\sigma_r(d z)
=\begin{cases}
d \phi, & \mbox{if $z \in \gamma_1
:=\{e^{\sqrt{-1} \phi}: \phi \in [0, 2 \pi) \}$},
\cr
r d \phi, & \mbox{if $z \in \gamma_q 
:=\{q e^{\sqrt{-1} \phi}: \phi \in [0, 2 \pi) \}$},
\end{cases}
\]
which we write as $H_r^2(\A_q)$.
A CONS of $H_r^2(\A_q)$ 
is given by $\{e^{(q, r)}_n\}_{n \in \Z}$ with
\[
e^{(q, r)}_n(z) = \frac{z^n}{\sqrt{1+r q^{2n}}},
\quad z \in \A_q, \quad n \in \Z,
\]
and the reproducing kernel is given by \cite{MS94}
\begin{equation}
S_{\A_q}(z, w; r)
= \sum_{n \in \Z} e^{(q,r)}_n(z) \overline{e^{(q,r)}_n(w)}
= \sum_{n=-\infty}^{\infty} 
\frac{(z \wbar)^n}{1+r q^{2n}}. 
\label{eqn:SAqr1}
\end{equation}
This infinite series converges absolutely for $z,w \in \A_q$.

When $r=q$, this Hilbert function space isreduced to
the Hardy space for $\A_q$, $H^2(\A_q)$,
introduced in Example \ref{thm:Hardy_Aq} 
in Section \ref{sec:examples}, associated with the
reproducing kernel 
\[
S_{\A_q}(\cdot, \cdot) :=S_{\A_q}(\cdot, \cdot; q).
\]
The kernel (\ref{eqn:SAqr1}) 
with a parameter $r >0$ is called the 
\textbf{weighted Szeg\H{o} kernel} of $\A_q$
and $H_r^2(\A_q)$ is the
\textbf{reproducing kernel Hilbert space}
(RKHS) with respect to 
$S_{\A_q}(\cdot, \cdot; r)$ \cite{MS94}.
We call $r$ the \textbf{weight parameter} \cite{KS22b}.
Notice that (\ref{eqn:SAqr1}) implies that
$S_{\A_q}(z, z; r)$ is a monotonically decreasing
function of the weight parameter
$r \in (0, \infty)$ for each fixed
$z \in \A_q$.

From now on, we put
\begin{equation}
p=q^2
\label{eqn:pq}
\end{equation}
for the nome $p$ of the theta function (\ref{eqn:theta}),
and assume that
$\theta(z) = \theta(z; q^)$,
where $q \in (0, 1)$ specifies the
radius of the inner circle of the annulus $\A_q$.
\begin{prop}[McCullough and Shen \cite{MS94}]
\label{thm:S_Aq_theta}
For $r > 0$
\begin{equation}
S_{\A_q}(z, w; r)
= \frac{q_0^2 \theta(-r z \wbar)}
{\theta(-r, z \wbar)},
\quad z, w \in \A_q, 
\label{eqn:S_qt_theta}
\end{equation}
where 
\begin{equation}
q_0:=\prod_{n \in \N}(1-q^{2n})=(q^2; q^2)_{\infty}.
\label{eqn:q0}
\end{equation}
In particular, 
\begin{equation}
S_{\A_q}(z, w) =S_{\A_q}(z, w; q)
= \frac{q_0^2 \theta(-q z \wbar)}
{\theta(-q, z \wbar)},
\quad z, w \in \A_q.
\label{eqn:S_q_theta}
\end{equation}
\end{prop}
\noindent{\it Proof.} 
The following function has been studied 
in \cite{MS94,Coo00,Ven12}, 
\begin{equation}
f^{\rm JK}(z, a; q) 
:= \sum_{n \in \Z} \frac{z^n}{1-aq^{2n}}, 
\label{eqn:JK1}
\end{equation}
with $q^2 < |z| <1$, 
$a \notin \{ q^{2 i} : i \in \Z\}$,
which is called
the \textbf{Jordan--Kronecker function} \cite[p.59]{Ven12}.
We readily find that
\begin{equation}
S_{\A_q}(z, w; r)=f^{\rm JK}(z \wbar, -r; q),
\quad z, w \in \A_q.
\label{eqn:S_qt_JK}
\end{equation}
The \textbf{bilateral basic hypergeometric series}
in base $p$ with one numerator parameter $a$
and one denominator parameter $b$ is defined by 
\cite{GR04}
\[
{_1}\psi_1(a; b; p, z) 
= {_1}\psi_1 \Big[ \begin{array}{l} a \cr b \end{array}; p, z \Big]
:= \sum_{n \in \Z} \frac{(a;p)_{n}}{(b;p)_{n}} z^n,
\quad |b/a| < |z| < 1.
\]
The Jordan--Kronecker function (\ref{eqn:JK1}) 
is a special case
of the ${_1}\psi_1$ function \cite{Coo00,Ven12};
\[
f^{\rm JK}(z, a; q) = \frac{1}{1-a} {_1}\psi_1 (a; aq^2; q^2, z).
\]
The following equality is known as 
\textbf{Ramanujan's ${_1}\psi_1$ summation formula} 
\cite{Coo00,GR04,Ven12}, 
\[
\sum_{n \in \Z} \frac{(a;p)_{n}}{(b;p)_{n}} z^n
=\frac{(az, p/(az), p, b/a;p)_{\infty}}
{(z, b/(az), b, p/a; p)_{\infty}},
\quad |b/a| < |z| < 1.
\]
Combining the above two equalities
with an appropriate change of variables, 
we obtain \cite{Coo00,Ven12}
\begin{equation}
f^{\rm JK}(z, a; q) = 
\frac{(az, q^2/(az), q^2, q^2; q^2)_{\infty}}
{(z, q^2/z, a, q^2/a; q^2)_{\infty}}
= \frac{q_0^2 \theta(za; q^2)}
{\theta(z, a; q^2)}.
\label{eqn:f2}
\end{equation}
By (\ref{eqn:S_qt_JK}), (\ref{eqn:f2}) proves the
equality (\ref{eqn:S_qt_theta}). The proof is complete.
\qed

Here we notice that the following symmetries 
of $f^{\rm JK}$ are readily verified by (\ref{eqn:f2})
using (\ref{eqn:theta_inversion}) and (\ref{eqn:theta_qp})
\cite{Coo00,Ven12}.
\begin{align}
f^{\rm JK}(z, a) &= f^{\rm JK}(a, z),
\label{eqn:JK2a}
\\
f^{\rm JK}(z, a) &=-f^{\rm JK}(z^{-1}, a^{-1}),
\label{eqn:JK2b}
\\
f^{\rm JK}(z, a) &= z f^{\rm JK}(z, aq^2) = a f^{\rm JK}(zq^2, a).
\label{eqn:JK2c}
\end{align}

Since $\theta$ is holomorphic in $\C^{\times}$, 
the expression (\ref{eqn:S_qt_theta}) implies 
that $S_{\A_q}(z, \alpha; r)$ is 
meromorphic in $\C^{\times}$. 
By the basic properties of $\theta$ 
given in Section \ref{sec:theta}, the
following is easily verified.
\begin{lem}
\label{thm:zero_S}
Assume that $\alpha \in \A_q$.
Then $S_{\A_q}(z, \alpha; r)$ has zeros at
$z=-q^{2i}/(\alphabar r)$, $i \in \Z$ in $\C^{\times}$. 
In particular, 
$S_{\A_q}(z, \alpha)$ has a unique zero in $\A_q$ at
$z=\alphahat$ given by (\ref{eqn:alphabar}).
\end{lem}
We have found the answer to the 
question {\bf (Q)} addressed in Section \ref{sec:examples}
in the paper by McCullough and Shen \cite{MS94}.

\begin{prop}[Mccullough and Shen \cite{MS94}]
\label{thm:MS}
The equality 
\begin{equation}
S_{\A_q}^{\alpha}(z, w; r)
=S_{\A_q}(z, w; r |\alpha|^2)
h_{\alpha}^q(z) \overline{h_{\alpha}^q(w)},
\quad z, w, \alpha \in \A_q, 
\label{eqn:SaAq}
\end{equation}
holds with 
\begin{equation}
h_{\alpha}^q(z):=
z \frac{\theta(\alpha/z)}
{\theta(\alphabar z)}
=- \alpha \frac{\theta(z/\alpha)}{\theta(z \alphabar)}, 
\quad z, \alpha \in \A_q.
\label{eqn:haq}
\end{equation}
The special case with $r=q$ answers 
the question {\bf (Q)} as
\begin{equation}
S_{\A_q}^{\alpha}(z, w)
=S_{\A_q}(z, w; q |\alpha|^2)
h_{\alpha}^q(z) \overline{h_{\alpha}^q(w)},
\quad z, w, \alpha \in \A_q.
\label{eqn:SaAq_ans}
\end{equation}
\end{prop}
\noindent{\it Proof.} 
We put (\ref{eqn:conditionSat}) with (\ref{eqn:S_qt_theta})
and (\ref{eqn:haq}) to (\ref{eqn:SaAq}), then
the equality is expressed by
theta functions. After multiplying both sides by
the common denominator, we see that the
equality (\ref{eqn:SaAq}) is equivalent to the following,
\begin{align}
&
\theta(-r z \wbar, -r |\alpha|^2, \alphabar z, \alpha \wbar)
- \theta(-r \alphabar z, -r \alpha \wbar, z \wbar, |\alpha|^2)
\nonumber\\
& \qquad 
=z \wbar \theta(-r z \wbar |\alpha|^2, \alpha z^{-1}, 
\alphabar \, \wbar^{-1}, -r).
\label{eqn:Weierstrass_add2}
\end{align}
Now we change the variables from
$\{z, \wbar, \alpha, r\}$ to
$\{x, y, u, v \}$ as
$\alphabar z =x/y$, 
$\alpha \wbar=u/v$, 
$z \wbar = x/v$, 
$|\alpha|^2= u/y$, and
$r=-yv$.
Then LHS 
of (\ref{eqn:Weierstrass_add2}) becomes 
$\theta(xy, x/y, uv, u/v)
-\theta(xv, x/v, uy, u/y)$, 
and RHS side becomes 
$(x/v) \theta(yv, (y/v)^{-1}, xu, (x/u)^{-1})$
which is equal to 
$(u/y) \theta(yv, y/v, xu, x/u)$
by (\ref{eqn:theta_inversion}).
Hence Weierstrass' addition formula
(\ref{eqn:Weierstrass_add1}) proves the equality
(\ref{eqn:Weierstrass_add2}). 
The proof is complete.
\vskip 0.3cm

Following the formula (\ref{eqn:condition})
for $S_{\A_q}^{\alpha}$, 
conditional kernels
$S_{\A_q}^{\alpha_1, \dots, \alpha_n}$ are inductively defined as
\begin{equation}
S_{\A_q}^{\alpha_1, \dots, \alpha_n}(z, w)
=(S_{\A_q}^{\alpha_1, \dots, \alpha_{n-1}})^{\alpha_n}(z,w),
\quad
z, w, 
\alpha_1, \dots, \alpha_n \in D,
\quad n=2,3,\dots.
\label{eqn:conditionS2}
\end{equation}
The kernels $S_{\A_q}^{\alpha_1, \dots, \alpha_n}$,
$n=2, 3, \dots$, 
will construct Hilbert subspaces of the Hardy space 
for $H^2(\A_q)$, 
\[
H^2_{\alpha_1, \dots, \alpha_n} (\A_q)
:= \{f \in \cH_{k} :f(\alpha_1)=\cdots=f(\alpha_n)=0\}.
\]
For $n \in \N$, $\alpha_1, \dots, \alpha_n \in \A_q$, define
\begin{equation}
\gamma^q_{\{\alpha_{\ell}\}_{\ell=1}^n}(z) 
:= \prod_{\ell=1}^n h_{\alpha_{\ell}}^q(z),
\quad z \in \A_q.
\label{eqn:gamma}
\end{equation}
Then the Mccullough and Shen formula (\ref{eqn:SaAq}) 
\cite{MS94} is generalized as 
\begin{equation}
S_{\A_q}^{\alpha_1, \dots, \alpha_n}(z, w; r)
=S_{\A_q} \Big(z, w; r \prod_{\ell=1}^n |\alpha_{\ell}|^2 \Big)
\gamma^q_{\{\alpha_{\ell}\}_{\ell=1}^n}(z) 
\overline{
\gamma^q_{\{\alpha_{\ell}\}_{\ell=1}^n}(w)
}, \quad z, w \in \A_q,
\label{eqn:MS_general}
\end{equation}
for $n \in \N$, $\alpha_1, \dots, \alpha_n \in \A_q$. 
As a special case of (\ref{eqn:MS_general}), we have
\begin{equation}
S_{\A_q}^{\alpha, \alphahat}(z, w)
=S_{\A_q}(z, w) f^q_{\alpha}(z) \overline{f^q_{\alpha}(w)},
\quad z, w \in \A_q
\label{eqn:A1}
\end{equation}
with
\begin{align}
f^q_{\alpha}(z) &:= \frac{1}{z} h^q_{\alpha}(z) h^q_{\alphahat}(z)
\nonumber\\
&=z 
\frac{\theta(-q z \alphabar, \alpha/z)}
{\theta(-q z/\alpha, \alphabar z)}
=-\alpha \frac{\theta(-q z \alphabar, z/\alpha)}
{\theta(-qz/\alpha, z \alphabar)}.
\label{eqn:Ahlfors1}
\end{align}
The equation (\ref{eqn:A1}) is nothing but 
(\ref{eqn:SAq_a_ahat2}), and 
$f^q_{\alpha}$ given by (\ref{eqn:Ahlfors1})
is identified with the
\textbf{Ahlfors map} from $\A_q$ to $\D$.

We can prove the following.
\begin{lem}
\label{thm:Blaschke}
For $\alpha \in \A_q$,
\begin{align}
\mbox{(i)} \quad &
h_{\alpha}^q(\alpha) =0,
\label{eqn:h_aq_1}
\\
\mbox{(ii)} \quad &
0< |h_{\alpha}^q(z)| < 1
\quad \forall z \in \A_q \setminus \{\alpha\},
\label{eqn:h_aq_2}
\\
\mbox{(iii)} \quad &
|h_{\alpha}^q(z)| =
\begin{cases}
1, & \quad \mbox{if $z \in \gamma_1 :=\{z \in \C: |z|=1\}$},
\cr
|\alpha|, & \quad \mbox{if $z \in \gamma_q := \{z \in \C: |z|=q \}$},
\end{cases}
\label{eqn:h_aq_3}
\\
\mbox{(iv)} \quad &
{h_{\alpha}^q}'(\alpha)
= - \frac{\theta'(1)}{\theta(|\alpha|^2)}
=\frac{q_0^2}{\theta(|\alpha|^2)}>0, 
\label{eqn:h_aq_4}
\\
\mbox{(v)} \quad &
\lim_{q \to 0} h_{\alpha}^q(z)
=\frac{z-\alpha}{1-\overline{\alpha} z}
=z \frac{1-\alpha/z}{1-\alphabar z}
=: h_{\alpha}(z).
\label{eqn:Mob1}
\end{align}
\end{lem}
Hence, 
$h_{\alpha}^q$ is identified with a \textbf{conformal map}
from $\A_{q}$ to the unit disk with a \textbf{circular slit} in it,
in which $\alpha \in \A_q$ is sent to the origin \cite{MS94}.
Notice that $h_{\alpha}$ is 
the \textbf{M\"{o}bius transformation}
$\D \to \D$ sending $\alpha$ to the origin.

\begin{figure}[ht]
\begin{center}
\includegraphics[scale=0.4]{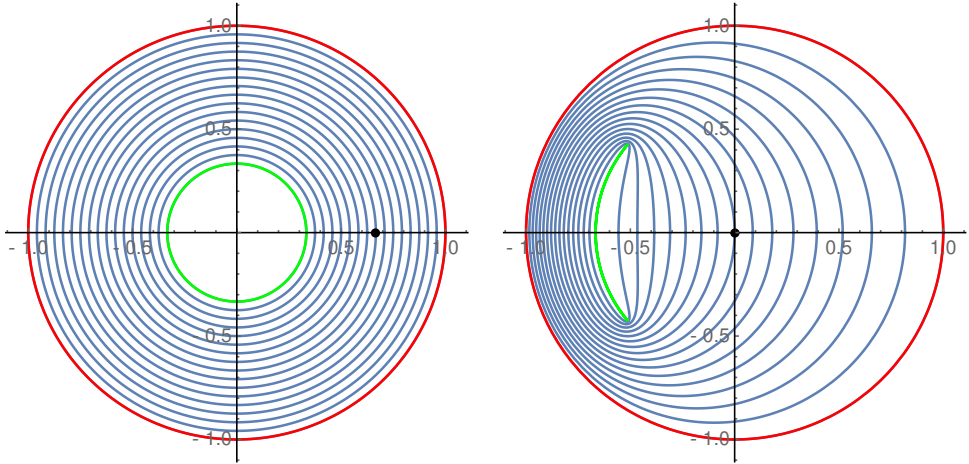}
\end{center}
\caption{Conformal map $h^{q}_{\alpha} : 
\A_{q} \to \D \setminus \{\mbox{a circular slit}\}$
is illustrated for $q=1/3$ and $\alpha=2/3$. 
The point $\alpha=2/3$ in $\A_{1/3}$ 
is mapped to the origin.
The outer boundary $\gamma_1$ of $\A_{1/3}$ 
(denoted by a red circle) is mapped to a unit circle
(a red circle) making the boundary of $\D$. 
The inner boundary $\gamma_{1/3}$ of $\A_{1/3}$
(a green circle) is mapped to a circular slit (denoted by a green arc)
which is a part of the circle with radius $\alpha=2/3$,
where the map is two-to-one
except the two points on $\gamma_{1/3}$
mapped to the two edges of the circular slit.}
\label{fig:h_alpha}
\end{figure}

\subsection{I.i.d.Gaussian power series, Laurent series, 
and their zero point processes}
\label{sec:Laurent}
For a domain $D \subset \C$, let $X$ be a random variable
on a probability space which takes values 
in the space of analytic functions on $D$. 
If $(X(z_1), \dots, X(z_n))$ follows a centered (mean zero) 
complex Gaussian distribution for every $n \in \N$
and every $z_1, \dots, z_n \in D$,
$X$ is said to be a \textbf{Gaussian analytic function} (GAF)
\cite{HKPV09}.
The zero set of $X$ is regarded as 
a point process on $D$ 
denoted by a nonnegative-integer-valued 
Radon measure
\[
\cZ_{X}=\sum_{z \in D: X(z)=0} \delta_z, 
\]
and it is simply called the
\textbf{zero point process} of the GAF.
Zero-point processes of GAFs have been extensively studied 
in quantum and statistical physics as solvable models of 
quantum chaotic systems and 
interacting particle systems. 
Many important characterizations of their probability laws 
have been reported in probability theory.

A typical example of GAF 
is provided by the 
\textbf{i.i.d.~Gaussian power series}
defined on the unit disk $\D$: 
Let $\N_0 :=\{0,1,2,\dots\}$ and $\{\zeta_n\}_{n \in \N_0}$ 
be i.i.d.~standard complex Gaussian random variables 
with density 
\[
{\rm p}(z)=\frac{1}{\pi} e^{-|z|^2}
=\frac{1}{\sqrt{\pi}} e^{-x^2}
\cdot
\frac{1}{\sqrt{\pi}} e^{-y^2},
\quad z= x+ \sqrt{-1} y, \quad x, y \in \R,
\]
and consider a random power series, 
\begin{equation}
X_{\D}(z)=\sum_{n=0}^{\infty} \zeta_n z^n,
\label{eqn:GAF_D}
\end{equation}
which defines an analytic function on $\D$ a.s.
This gives a GAF on $\D$ with a \textbf{covariance kernel}
\begin{equation}
\bE[X_{\D}(z) \overline{X_{\D}(w)}]
= \frac{1}{1-z \wbar}=:S_{\D}(z, w), \quad z, w \in \D.
\label{eqn:SD}
\end{equation}
This kernel is identified with 
the \textbf{Szeg\H{o} kernel} of $\D$ (\ref{eqn:S_D})
of the Hardy space $H^2(\D)$ introduced
in Section \ref{sec:examples}). 
Peres and Vir\'ag \cite{PV05} proved that 
$\cZ_{X_{\D}}$ is a \textbf{determinantal point process} (DPP) 
such that the \textbf{correlation kernel}
is given by 
$S_{\D}(z,w)^2=(1-z \wbar)^{-2}$, 
$z,w \in \D$ 
with respect to the reference measure 
$\lambda=m/\pi$.
Here $m$ represents the Lebesgue measure on $\C$; 
$m(dz) :=dx dy$, $z=x+\sqrt{-1} y \in \C$.
This correlation kernel is identified with
the reproducing kernel of the Bergman space on $\D$
(\ref{eqn:K_D}), which is 
called \textbf{Bergman kernel} of $\D$ in Section \ref{sec:examples}.
Thus the study of Peres and Vir\'ag on
$X_{\D}$ and $\cZ_{X_{\D}}$
is associated with the following relationship 
between kernels on $\D$ \cite{PV05},
\begin{equation}
\bE[X_{\D}(z) \overline{X_{\D}(w)}]^2
=S_{\D}(z, w)^2
=K_{\D}(z, w), 
\quad z, w \in \D.
\label{eqn:SKC_D}
\end{equation}

Associated with the RKHS $H_r^2(\A_q)$
studied in Section \ref{sec:Ahlfors}
we consider the 
\textbf{Gaussian Laurent series}
\begin{equation}
X_{\A_q}^r(z):= \sum_{n \in \Z} \zeta_n e^{(q, r)}_n(z)
=\sum_{n=-\infty}^{\infty}
\zeta_n \frac{z^n}{\sqrt{1+r q^{2n}}},
\label{eqn:GAF_Aqr}
\end{equation}
where $\{\zeta_n\}_{n \in \Z}$ are i.i.d.~standard complex
Gaussian random variables with density $e^{-|z|^2}/\pi$. 
Since $\lim_{n \to \infty} |\zeta_n|^{1/n} = 1$ a.s., 
we apply the Cauchy--Hadamard
criterion to the positive 
and negative powers of $X_{\A_q}^r(z)$ separately
to conclude that this random Laurent series converges 
a.s. whenever $z \in \A_q$. 
Moreover, since the distribution $\zeta_n$ is symmetric, 
both of $\gamma_1$ and $\gamma_q$ are natural 
boundaries.
Hence $X_{\A_q}^r$ provides a GAF on $\A_q$ whose 
covariance kernel is given by 
the \textbf{weighted Szeg\H{o} kernel of $\A_q$}, 
\[
\bE[X_{\A_q}^r(z) \overline{X_{\A_q}^r(w)} ]
=S_{\A_q}(z, w; r), \quad z, w \in \A_q,
\]
and the zero point process is denoted by
$\cZ_{X_{\A_q}^r} :=\sum_{z \in \A_q: X^r_{\A_q}(z)=0} \delta_{z}$.
In particular, we write
$X_{\A_q}(z) :=X_{\A_q}^q(z), z \in \A_q$ and
$\cZ_{X_{\A_q}} := \cZ_{X_{\A_q}^q}$
as mentioned above. 

We recall \textbf{Schottky's theorem}: 
The group of conformal transformations
from $\A_q$ to itself 
is generated by the rotations 
and the $q$-inversions 
\[
T_q(z) := \frac{q}{z}. 
\]
The invariance of the present 
GAF and its zero point process under rotation is obvious. 
Using the properties of $S_{\A_q}$, 
we can prove the following.
\begin{prop}
\label{thm:conformal_Aq}
\begin{description}
\item{\rm (i)} 
The GAF $X_{\A_q}^r$ given by (\ref{eqn:GAF_Aqr}) 
has the 
\textbf{$(q, r)$-inversion symmetry} in the sense that
\[
\Big\{ (T_q'(z))^{1/2} X_{\A_q}^r (T_q(z)) \Big\}
\dis=
\Big\{ \sqrt{\frac{q}{r}} X_{\A_q}^{q^2/r}(z) \Big\},
\quad z \in \A_q,
\]
where $T_q'(z) := \frac{d T_q}{dz}(z)=-q/z^2$.
\item{\rm (ii)} 
For $\cZ_{X_{\A_q}^r} = \sum_i \delta_{Z_i}$,
let 
$T_q^{*}\cZ_{X_{\A_q}^r} := \sum_i \delta_{T_q^{-1}(Z_i)}$.
Then
$T_q^{*}\cZ_{X_{\A_q}^r} \dis= \cZ_{X_{\A_q}^{q^2/r}}$.
\item{\rm (iii)} 
In particular, when $r=q$, 
the GAF $X_{\A_q}$ is invariant under
conformal transformations which preserve $\A_q$,
and so is its zero point process $\cZ_{X_{\A_q}}$.
\end{description}
\end{prop}
\noindent{\it Proof.} 
The equivalence in distribution asserted by (i) 
is implied by proving the following equality
of the covariance kernel given by the 
weighted  Szeg\H{o} kernel of $\A_q$,
\begin{equation}
\sqrt{T_q'(z)} \overline{\sqrt{T_q'(w)}}
S_{\A_q}(T_q(z), T_q(w) ; r) 
= \frac{q}{r} S_{\A_q} (z, w; q^2/r ).
\label{eqn:equal_S}
\end{equation}
If we use the expression (\ref{eqn:S_qt_JK}) 
of $S_{\A_q}(z, w; r)$,
the symmetries 
of the Jordan--Kronecker function 
(\ref{eqn:JK2b}) and (\ref{eqn:JK2c})
prove (\ref{eqn:equal_S}). 
By the definition of zero point processes, 
(ii) is concluded from (i). 
When $r=q$, (\ref{eqn:equal_S}) becomes
\[
\sqrt{T_q'(z)} \overline{\sqrt{T_q'(w)}} S_{\A_q}(T_q(z), T_q(w); q) 
=\sqrt{T_q'(z)} \overline{\sqrt{T_q'(w)}} S_{\A_q}(T_q(z), T_q(w)) 
=  S_{\A_q}(z, w), 
\]
which implies the invariance of the GAF $X_{\A_q}$ under conformal 
transformations preserving $\A_q$
by Schottkey's theorem and (iii) is proved. \qed

We can give 
\textbf{probabilistic interpretations}
of the above facts as follows.
\begin{prop}
\label{thm:equivalence2}
For any $\alpha_1, \dots, \alpha_n \in \A_q$, $n \in \N$, 
the following hold. 
\begin{description}
\item{\rm (i)} \, 
The following equality is established,
\begin{align*}
&\mbox{$\{X_{\A_q}^r(z) : z \in \A_q \}$ 
given $\{X_{\A_q}^r(\alpha_1)= \cdots =X_{\A_q}^r(\alpha_n)=0\}$}
\nonumber\\
& \hskip 4cm
\dis= 
\left\{ \gamma^q_{\{ \alpha_{\ell} \}_{\ell=1}^n}(z) 
X_{\A_q}^{r \prod_{\ell=1}^n |\alpha_{\ell}|^2}(z)
 : z \in \A_q \right\}.
\end{align*}
\item{\rm (ii)} \, 
Let $\cZ_{X_{\A_q}^r}^{\alpha_1, \dots, \alpha_n}$ 
denote the zero point process
of the GAF $X_{\A_q}^r(z)$
given $\{X_{\A_q}^r(\alpha_1)=\cdots = X_{\A_q}^r(\alpha_n)=0\}$. 
Then, 
\[
\cZ_{X_{\A_q}^r}^{\alpha_1, \dots, \alpha_n} 
\dis= \cZ_{X_{\A_q}^{r \prod_{\ell=1}^n |\alpha_{\ell}|^2}}
+ \sum_{i=1}^n \delta_{\alpha_i}.
\]
\end{description}
\end{prop}
\begin{rem}
\label{rem:Remark_condition1}
For the GAF on $\D$
studied by Peres and Vir\'ag \cite{PV05},
$\{X_{\D}(z) : z \in \D \}$ 
given $\{X_{\D}(\alpha)=0\}$
is equal in law to $\{h_{\alpha}(z) X_{\D}(z) : z \in \D \}$, 
$\forall \alpha \in \D$, 
where $h_{\alpha}$ is given by (\ref{eqn:Mob1}),
and then, in the notation 
used in Proposition \ref{thm:equivalence2}, 
\[
\cZ_{X_{\D}}^{\alpha} \dis= \cZ_{X_{\D}}+\delta_{\alpha}, 
\quad 
\forall \alpha \in \D.
\]
Hence, no new GAF nor new zero point process 
appear by conditioning of zeros.
For the present GAF on $\A_q$, however, 
conditioning of zeros induces new GAFs and 
new zero point processes
as shown by Proposition \ref{thm:equivalence2}.
Actually, by (\ref{eqn:SAqr1}) 
the covariance of the induced GAF
$X_{\A_q}^{r \prod_{\ell=1}^n |\alpha_{\ell}|^2}$ 
is expressed by
\[
S_{\A_q} \Big(z, w; r \prod_{\ell=1}^n |\alpha_{\ell}|^2 \Big)
=\sum_{n=-\infty}^{\infty} 
\frac{(z \wbar)^n}{
1+r \prod_{\ell=1}^n |\alpha_{\ell}|^2 q^{2n}}.
\]
Since $q < |\alpha_{\ell}| <1$,
as increasing the number of conditioning zeros, 
the variance of induced GAF monotonically
increases, in which the increment is a decreasing function of
$|\alpha_{\ell}| \in (q, 1)$. 
\end{rem}

\subsection{Permanental-determinantal point processes}
\label{sec:perdet}
\subsubsection{Correlation functions of point processes}
\label{sec:pp}

A \textbf{point process} is formulated as follows \cite{KS22}.
Let $S$ be a \textbf{base space}, 
which is locally compact Hausdorff space
with a countable base, and 
a \textbf{reference measure} $\lambda$ 
is given on $S$. 
The \textbf{configuration space} of a point process on $S$ is given by
\[
\Conf(S)
=\Big\{ \xi = \sum_i \delta_{x_i} : \mbox{$x_i \in S$,
$\xi(\Lambda) < \infty$ for all bounded set $\Lambda \subset S$} 
\Big\}.
\]
We say $\xi_n, n \in \N :=\{1, 2, \dots\}$ converges
to $\xi$ in the \textbf{vague topology}, if
\[
\int_{S} f(x) \xi_n(dx) \to
\int_{S} f(x) \xi(dx), \quad
\forall f \in \cC_{\rm c}(S), 
\]
where $\cC_{\rm c}(S)$ is the set of
all continuous real-valued functions 
with compact support.
A point process on $S$ is
a $\Conf(S)$-valued random variable $\Xi=\Xi(\cdot)$.
If $\Xi(\{ x \}) \in \{0, 1\}$ for any point $z \in S$,
then the point process is said to be \textbf{simple}.
Assume that $\Lambda_i, i=1, \dots, m$, 
$m \in \N$ are
disjoint bounded sets in $S$ and
$k_i \in \N_0, i=1, \dots, m$ satisfy
$\sum_{i=1}^m k_i = n \in \N_0$.
A symmetric measure $\lambda^n$ on $S^n$
is called the $n$-th \textbf{correlation measure},
if it satisfies
\[
\bE \Bigg[
\prod_{i=1}^m \frac{\Xi(\Lambda_i)!}
{(\Xi(\Lambda_i)-k_i)!} \Bigg]
=\lambda^n(\Lambda_1^{k_1} \times \cdots
\times \Lambda_m^{k_m}),
\]
where when $\Xi(\Lambda_i)-k_i < 0$,
we interpret $\Xi(\Lambda_i)!/(\Xi(\Lambda_i)-k_i)!=0$.
If $\lambda^n$ is absolutely continuous 
with respect to the $n$-product measure $\lambda^{\otimes n}$,
the Radon--Nikodym derivative
$\rho^n(x_1, \dots, x_n)$ is called
the \textbf{$n$-point correlation function}
with respect to the reference measure $\lambda$;
that is, 
\[
\lambda^n(dx_1 \cdots dx_n)
=\rho^n(x_1, \dots, x_n) 
\lambda^{\otimes n}(dx_1 \cdots dx_n).
\]

\subsubsection{Permanental-determinantal point processes
(PDPPs) with hierarchical structures}
\label{sec:PDPP}
\begin{figure}[ht]
\begin{center}
\includegraphics[scale=0.7]{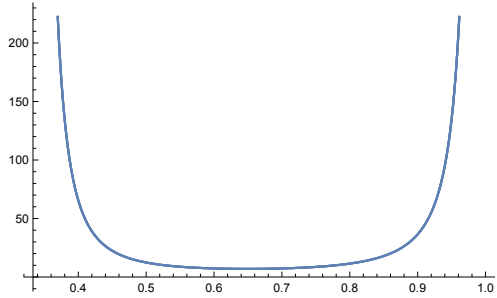}
\end{center}
\caption{Dependence on the radial coordinate $r=|z| \in (q, 1)$
of the density function $\rho^1_{\A_q}(z; r)$ 
in the case $q=r=1/3$. }
\label{fig:density_Aq}
\end{figure}

We introduce the following 
notation. For an $n \times n$ matrix
$M=(m_{ij})_{1 \leq i, j \leq n}$,
\begin{equation}
\perdet M = \perdet_{1 \leq i, j \leq n} [m_{ij}]
:=\per M \det M,
\label{eqn:perdet}
\end{equation}
that is, $\perdet M$ denotes
$\per M$ multiplied by $\det M$.
Note that
$\perdet$ is a special case
of \textbf{hyperdeterminants} introduced
by Gegenbauer following Cayley \cite{KS22b}. 
If $M$ is a positive semidefinite Hermitian matrix, then
$\per M \geq \det M \geq 0$,
and hence $\perdet M \geq 0$ by the definition (\ref{eqn:perdet}).
The following is the main theorem in Section \ref{sec:GAF}. 
\begin{thm}
\label{thm:mainA1}
The zero point process $\cZ_{X_{\A_q}^r}$ 
on $\A_q$ is a 
\textbf{permanental-determinantal point process}
(PDPP) in the sense that 
it has correlation functions
$\{\rho^{n}_{\A_q}\}_{n \in \N}$ 
given by
\begin{equation}
\rho^{n}_{\A_q}(z_1,\dots,z_n; r) 
=
\frac{\theta(-r)}{\theta( -r \prod_{k=1}^n |z_k|^4)}
\perdet_{1 \leq i, j \leq n}
\Big[
S_{\A_q} \Big(z_{i}, z_{j}; r \prod_{\ell=1}^n |z_{\ell}|^2 \Big) 
\Big]
\label{eqn:rho_mainA1}
\end{equation}
for every $n \in \N$ and
$z_1, \dots, z_n \in \A_q$
with respect to $m/\pi$.
\end{thm}
The \textbf{density of zeros on $\A_q$}
with respect to $m/\pi$ is given by 
\begin{equation}
\rho^1_{\A_q}(z; r) = 
\frac{\theta(-r)}{\theta(-r |z|^4)}
S_{\A_q}(z,z; r|z|^2)^2
=
\frac{q_0^4 \theta(-r, -r |z|^4)}
{\theta(-r |z|^2, |z|^2)^2},
\quad z \in \A_q,
\label{eqn:density_Aqt}
\end{equation}
which is always positive. 
Since $\rho^1_{\A_q}(z; r)$ depends only on 
the modulus of the coordinate $|z| \in (q, 1)$, 
the PDPP is rotationally invariant. 
As shown by 
(\ref{eqn:theta_signs})--(\ref{eqn:max_theta}) 
in Section \ref{sec:theta}, 
in the interval $x \in (-\infty, 0)$, $\theta(x)$ is positive and
strictly convex with $\lim_{x \downarrow -\infty} \theta(x)
=\lim_{x \uparrow 0} \theta(x)=+\infty$, while
in the interval $x \in (q^2,1)$, $\theta(x)$ is positive
and strictly concave 
with
$\theta(x) \sim
q_0^2(x-q^2)/q^2$ as $x \downarrow q^2$
and $\theta(x) \sim q_0^2(1-x)$ as $x \uparrow 1$.
Therefore, the density shows divergence both at the inner
and outer boundaries as
\begin{equation}
\rho^1_{\A_q}(z; r) \sim
\begin{cases}
\displaystyle{
\frac{q^2}{ (|z|^2-q^2)^2}
},
& |z| \downarrow q,
\cr
\displaystyle{
\frac{1}{(1-|z|^2)^2}
},
& |z| \uparrow 1,
\end{cases}
\label{eqn:density_edges}
\end{equation}
which is independent of $r$ and implies $\bE[\cZ_{X^r_{\A_q}}(\A_q)]=\infty$. See Fig.\ref{fig:density_Aq}.

If $M$ is a $2 \times 2$ matrix, we see that
$\perdet M=\det (M \circ M)$,
where $M \circ M$ denotes the Hadamard product of $M$, 
i.e., entrywise multiplication, $(M \circ M)_{ij}=M_{ij} M_{ij}$. 
Then the two-point correlation is expressed by a single 
determinant as
\begin{equation}
\rho^{2}_{\A_q}(z_1, z_2; r)
= \frac{\theta(-r)}{\theta( -r |z_1|^4|z_2|^4)}
\det_{1 \leq i, j \leq 2} 
\left[S_{\A_q}(z_i, z_j ; r |z_1|^2 |z_2|^2)^2
\right], 
\quad z_1, z_2 \in \A_q.
\label{eqn:rho2_Aqt}
\end{equation}

\begin{rem}\label{rem:PDPP_property}
(i) The PDPP with correlation functions \eqref{eqn:rho_mainA1} 
turns out to be a \textbf{simple point process}, 
i.e., there is no multiple point a.s., 
due to the existence of two-point correlation 
function with respect to
the Lebesgue measure $m/\pi$ \cite[Lemma 2.7]{Kal17}. 
(ii) Using the explicit expression (\ref{eqn:rho_mainA1}) together 
with the \textbf{Frobenius determinantal formula}
\eqref{eqn:Frobenius_formula} given below, 
we can verify that for every $n \in \N$, 
the $n$-point correlation $\rho^n_{\A_q}(z_1, \dots, z_n) >0$
if all coordinates $z_1, \dots, z_n \in \A_q$ are 
different from each other, 
and that $\rho_{\A_q}^n(z_1, \dots, z_n) =0$ if
some of  $z_1, \dots, z_n$ coincide; 
e.g., $z_i = z_j, i \not=j$,  
by the determinantal factor in $\perdet$ (\ref{eqn:perdet}).
\end{rem}

\begin{rem}
\label{rem:pdet_correlation}
The \textbf{determinantal point processes} (DPPs) 
and the \textbf{permanental point processes} (PPPs)
have the $n$-correlation functions of the forms  
\begin{align*}
\rho_{\mathrm{DPP}}^n(z_1,\dots,z_n) 
= \det_{1 \le i,j \le n}[K(z_i,z_j)], \quad 
\rho_{\mathrm{PPP}}^n(z_1,\dots,z_n) 
= \per_{1 \le i,j \le n}[K(z_i,z_j)], 
\end{align*}
respectively (cf. \cite{ST03a,HKPV09}). 
Due to \textbf{Hadamard's inequality} for the determinant \cite[Section II.4]{MM92}
and \textbf{Lieb's inequality for the permanent} \cite{Lie66}, 
we have 
\begin{align*}
 \rho_{\mathrm{DPP}}^2(z_1,z_2) \le \rho_{\mathrm{DPP}}^1(z_1) 
 \rho_{\mathrm{DPP}}^1(z_2), \quad 
 \rho_{\mathrm{PPP}}^2(z_1,z_2) \ge \rho_{\mathrm{PPP}}^1(z_1) 
 \rho_{\mathrm{PPP}}^1(z_2).
\end{align*}
These correlation inequalities suggest a repulsive nature
(negative correlation) for DPPs and 
an attractive nature (positive correlation) for PPPs. 
Since $\perdet$ is considered to have intermediate nature 
between determinant and permanent, 
PDPPs are expected to exhibit 
\textbf{both repulsive and attractive characters}, 
depending on the position of two points $z_1$ and $z_2$. 
For example, Remark~\ref{rem:PDPP_property} (ii) shows 
the repulsive nature inherited from the DPP side. 
The two-sidedness of the present PDPP will be summarized in Theorem~\ref{thm:unfolded} given below. 
\end{rem}
\begin{rem}
\label{rem:Remark_condition2}
The \textbf{$(q,r)$-inversionb symmetry} given by 
Proposition \ref{thm:conformal_Aq} (ii) is 
rephrased using correlation functions as
\begin{equation}
\rho^{n}_{\A_q}(T_q(z_1),\dots, T_q(z_n); r) 
\prod_{\ell=1}^n |T_q'(z_{\ell})|^2
= \rho^{n}_{\A_q}(z_1, \dots, z_n; q^2/r)
\label{eqn:rho_inversion0}
\end{equation}
for any $n \in \N$ and $z_1, \dots, z_n \in \A_q$,
where $T_q(z)=q/z$ and $|T_q'(z)|^2=q^2/|z|^4$.
In the correlation functions 
$\{\rho_{\A_q}^n\}_{n \in \N}$ given by
Theorem \ref{thm:mainA1}, we see 
an \textbf{inductive structure} such that
the functional form of the 
permanental-determinantal correlation kernel
$S_{\A_q}(\cdot, \cdot; r \prod_{\ell=1}^n |z_{\ell}|^2)$
is depending on the
points $\{z_1, \dots, z_n\}$, which we intend to measure
by $\rho_{\A_q}^n$, via the weight parameter
$r \prod_{\ell=1}^n |z_{\ell}|^2$.
This is due to the inductive structure 
of the induced GAFs generated in conditioning of zeros as
mentioned in Remark \ref{rem:Remark_condition1}.
In addition, the reference measure $m/\pi$ is
also weighted by 
$\theta(-r)/\theta(-r \prod_{k=1}^n |z_{k}|^4)$.
Such a \textbf{hierarchical structure}
of correlation functions and reference measures 
is necessary to realize
the $(q,r)$-inversion symmetry (\ref{eqn:rho_inversion0})
and the invariance under \textbf{conformal transformations}
preserving $\A_q$ when $r=q$. 
\end{rem}

\subsubsection{Simpler but still non-trivial PDPP in 
$q \to 0$ limit}
\label{sec:simplerPDPP}

The above GAF and the PDPP induce the following limiting cases.
With fixed $r >0$ we take the limit $q \to 0$.
By the reason explained in Remark \ref{rem:q0limit} below, 
in this limiting procedure, we should consider the 
point processes $\{\cZ_{X^r_{\A_q}} : q>0\}$ to be defined on 
the \textbf{punctured unit disk} 
\[
\D^{\times} :=\{z \in \C : 0 < |z| < 1\}
\]
instead of $\D$.
Although the limit point process is given on $\D^{\times}$ 
by definition, 
it can be naturally viewed as a point process
defined on $\D$, which we will introduce below.  
Let $H_r^2(\D)$ be the 
\textbf{Hardy space on $\D$ with the
weight parameter $r>0$},
whose inner product is given by
\[
\bra f, g \ket_{H^2_r(\D)} 
=\frac{1}{2\pi} \int_0^{2 \pi} f(e^{\sqrt{-1} \phi})
\overline{g(e^{\sqrt{-1} \phi})} d \phi
+ r f(0) g(0), 
\quad f, g \in H^2_r(\D). 
\]
The reproducing kernel of $H^2_r(\D)$ is given by
\begin{align}
S_{\D}(z, w; r)
&= \sum_{n=0}^{\infty} e^{(0,r)}_n(z) \overline{e^{(0,r)}_n(w)}
= \frac{1}{1+r} + \sum_{n=1}^{\infty} (z \wbar)^n 
\nonumber\\
&=
\frac{1+r z \wbar}{(1+r)(1-z \wbar)}, 
\quad z, w \in \D.
\label{eqn:SA0t}
\end{align}
The GAF associated with $H_r^2(\D)$ is then
defined by
\begin{equation}
X_{\D}^r(z)
= \frac{\zeta_0}{\sqrt{1+r}} 
+  \sum_{n=1}^{\infty} \zeta_n z^n,
\quad z \in \D
\label{eqn:GAF_A0t}
\end{equation}
so that the covariance kernel is given by 
$\bE[X_{\D}^r(z) 
\overline{ X_{\D}^r(w)} ]
=S_{\D}(z,w;r)$, 
$z, w \in \D$.
For the conditional GAF given a zero at $\alpha \in \D$, 
the covariance kernel is given by
\[
S_{\D}^{\alpha}(z, w; r)
=S_{\D}(z, w; r |\alpha|^2)
h_{\alpha}(z) \overline{h_{\alpha}(w)},
\quad z, w, \alpha \in \D,
\]
where the replacement of the weight parameter
$r$ by $r |\alpha|^2$ should be done, even though
the factor $h_{\alpha}(z)$ is simply given by
the M\"{o}bius transformation (\ref{eqn:Mob1}).

For the zero point process
Theorem \ref{thm:mainA1} is reduced to
the following by the formula 
\[
\lim_{q \to 0} \theta(z; q^2)=1-z.
\]
\begin{cor}
\label{thm:mainA2}
Assume that $r >0$.
Then $\cZ_{X_{\D}^r}$ is 
a PDPP on $\D$ with the correlation functions
\begin{equation}
\rho^{n}_{\D}(z_1,\dots,z_n; r) 
=
\frac{1+r}{1+r \prod_{k=1}^n |z_k|^4}
\perdet_{1 \leq i, j \leq n}
\Big[S_{\D}\Big(z_{i}, z_{j}; r \prod_{\ell=1}^n |z_{\ell}|^2 \Big) \Big]
\label{eqn:rho_mainA2}
\end{equation} 
for every $n \in \N$ and
$z_1, \dots, z_n \in \D$
with respect to $m/\pi$.
In particular, the density of zeros on $\D$
is given by
\begin{equation}
\rho_{\D}^1(z;r)
= \frac{(1+r)(1+r|z|^4)}{(1+r|z|^2)^2 (1-|z|^2)^2},
\quad z \in \D.
\label{eqn:density_A0t}
\end{equation}
\end{cor}
\noindent
As $r$ increases the first term in (\ref{eqn:GAF_A0t}), 
which gives the value of the GAF at the origin, 
decreases and hence the variance at the origin,
$\bE[|X_{\D}^r(0)|^2]=S_{\D}(0,0;r)
=(1+r)^{-1}$ decreases monotonically.
As a result the density of zeros 
in the vicinity of the origin 
increases as $r$ increases.
Actually we see that $\rho_{\D}^1(0; r)=1+r$.
\begin{rem}
\label{rem:q0limit}
The asymptotics \eqref{eqn:density_edges} show that 
the density of zeros of 
$\cZ_{X^r_{\A_q}}$ diverges at the inner boundary
$\gamma_q =\{z : |z|=q\}$ for each $q>0$ 
while the density of $\cZ_{X^r_{\D}}$ is finite at the origin 
as in \eqref{eqn:density_A0t}. Therefore infinitely many zeros near 
the inner boundary $\gamma_q$
seem to vanish in the limit as $q \to 0$. 
This is the reason why we regard the base space of
$\{\cZ_{X^r_{\A_q}} : q>0\}$ and the limit point process
$\cZ_{X^r_{\D}}$ as $\D^{\times}$ instead of $\D$ 
as mentioned before.
Indeed, in the vague topology, 
with which we equip a configuration space, 
we cannot see configurations outside each compact set, hence 
infinitely many zeros are not observed 
on each compact set in $\D^{\times}$ (not $\D$)
for any sufficiently small $q>0$ 
depending on the compact set that we take. 
\end{rem}

\subsubsection{Reduction to GAF and DPP of 
Peres and Vir\'ag}
\label{sec:PV}
If we take the further limit $r \to 0$ in 
(\ref{eqn:SA0t}), we obtain the Szeg\H{o} kernel of $\D$
given by (\ref{eqn:SD}).
Since the matrix
$
( S_{\D}(z_i, z_j)^{-1} )_{1 \leq i, j \leq n}
=(1-z_i \overline{z_j})_{1 \leq i, j \leq n}
$
has rank 2, the following equality called 
\textbf{Borchardt's identity} holds
(see Theorem 3.2 in \cite{Min78}, Theorem 5.1.5 in \cite{HKPV09}),
\begin{equation}
\perdet_{1 \leq i, j \leq n} 
\Big[ (1-z_i \overline{z_j})^{-1} \Big]
= \det_{1 \leq i, j \leq n}
\Big[ (1-z_i \overline{z_j})^{-2} \Big].
\label{eqn:Borchardt}
\end{equation}
By the relation (\ref{eqn:SKC_D}), the $r \to 0$ limit
of $\cZ_{X_{\D}^r}$ is identified with
the DPP on $\D$, $\cZ_{X_{\D}}$,
studied by Peres and Vir\'ag \cite{PV05},
whose correlation functions are given by
\[
\rho_{\D, \mathrm{PV}}^n(z_1, \dots, z_n)
=\det_{1 \leq i, j \leq n} [K_{\D}(z_i, z_j)],
\quad n \in \N, \quad z_1, \dots, z_n \in \D,
\]
with respect to $m/\pi$ 

Here we give the result by Peres and Vir\'ag as 
a theorem.
\begin{thm}[Peres and Vir\'ag \cite{PV05}]
\label{thm:Peres_Virag}
$\cZ_{X_{\D}}$ is a DPP on $\D$
such that the correlation kernel 
with respect to $m/\pi$ 
is given by
the Bergman kernel $K_{\D}$ of $\D$ given by
(\ref{eqn:K_D}).
\end{thm}

\begin{rem}
\label{rem:Remark_r_infinity}
We see from (\ref{eqn:SA0t}) that
\[
\lim_{r \to \infty} S_{\D}(z,w; r)
= (1-z \wbar)^{-1} -1, \quad z, w \in \D,
\]
which can be identified with 
the conditional kernel 
given a zero at the origin; 
$S_{\D}^0(z, w)
= S_{\D}(z, w)- S_{\D}(z,0) S_{\D}(0, w)/S_{\D}(0,0)$
for $S_{\D}(z,0) \equiv 1$. 
In this limit we can use Borchardt's identity again, since
the rank of the matrix
$(S_{\D}(z_i, z_j; \infty)^{-1})_{1 \leq i, j \leq n}
=(z_i^{-1} z_j^{-1}-1)_{1 \leq i, j \leq n}$ is two. 
Then, thanks to the proper limit of the prefactor 
of $\perdet$ in (\ref{eqn:rho_mainA2})
when $z_k \in \D^{\times}$ for all $k=1,2,\dots, n$; 
\[
\lim_{r \to \infty}
\frac{1+r}{1+r \prod_{k=1}^n |z_k|^4}
=\prod_{k=1}^n |z_k|^{-4}, 
\]
we can verify that 
\[
\lim_{r \to \infty} \rho_{\D}^n(z_1, \dots, z_n; r)
= \rho_{\D, \mathrm{PV}}^n(z_1, \dots, z_n)
\]
for every $n \in \N$, and every 
$z_1, \dots, z_n \in \D^{\times}$. 
On the other hand, taking \eqref{eqn:GAF_A0t} into account, we have 
\[
X^{\infty}_{\D}(z)
= z \sum_{n=1}^{\infty} \zeta_n z^{n-1}
\dis= z X_{\D}(z), 
\] 
from which,
we can see that 
\[
\lim_{r \to \infty} \cZ_{X_{\D}^r}
=: \cZ_{X_{\D}^{\infty}} \dis= \cZ_{X_{\D}}+\delta_0, 
\]
that is,
the DPP of Peres and Vir\'ag with
a deterministic zero added at the origin.
This is consistent with the fact that $\rho_{\D}^1(0;r) =
1+r$ diverges as $r \to \infty$.
Since 
\[
\lim_{r \to 0} \cZ_{X_{\D}^r} 
=: \cZ_{X_{\D}^0} 
\dis= \cZ_{X_{\D}}
\]
as mentioned above, 
the one-parameter family of 
PDPPs $\{\cZ_{X_{\D}^r} : r \in (0, \infty)\}$ 
can be regarded as an 
\textbf{interpolation} between 
the DPP of Peres and Vir\'ag and 
that DPP with a deterministic zero added at the origin. 
\end{rem}

The distribution of $\cZ_{X_{\D}}$ is invariant under
M\"{o}bius transformations that preserve $\D$ \cite{PV05}.
This invariance is a special case of the following,
which can be proved using
the conformal transformations 
of the Szeg\H{o} kernel and the Bergman kernel 
given by (\ref{eqn:SD_KD_conformal}) 
\cite{PV05,HKPV09}. 

\begin{prop}[Peres and Vir\'ag \cite{PV05}]
\label{thm:Peres_Virag2}
Let $\widetilde{D} \subsetneq \C$ be a 
simply connected domain wit
analytic boundary. 
Then there is a GAF $X_{\widetilde{D}}$ with 
covariance kernel
$\bE[X_{\widetilde{D}}(z) \overline{X_{\widetilde{D}}(w)}] 
= S_{\widetilde{D}}(z,w)$, $z, w \in \widetilde{D}$,
where $S_{\widetilde{D}}$ denotes 
the Szeg\H{o} kernel of $\widetilde{D}$.
The zero point process $\cZ_{X_{\widetilde{D}}}$ 
is the DPP such 
that the correlation kernel is given by
the Bergman kernel $K_{\widetilde{D}}$ of $\widetilde{D}$. 
This DPP is conformally invariant in the following sense.  
If $D \subsetneq \C$ is another 
simply connected domain with analytic
boundary, and 
\[
f: D \to \widetilde{D}
\] 
is a \textbf{conformal transformation}, then 
for $\cZ_{X_{D}}=\sum_i \delta_{Z_i}$, 
\[
\sum_i \delta_{f^{-1}(Z_i)}=:
f^{*}\cZ_{X_{\widetilde{D}}} \dis= \cZ_{X_{D}}.
\] 
That is, $f^{*}\cZ_{X_{\widetilde{D}}}$
is a DPP such that the correlation kernel
is equal to the Bergman kernel $K_{D}$ of $D$.
\end{prop}
We can say that the zero point process
of Peres and Vir\'ag \cite{PV05}
is a \textbf{conformal invariant DPP}.

\subsection{Proof of Theorem \ref{thm:mainA1}}
\label{sec:proof_main}
\subsubsection{Edelman-Kostlan formulas and 
second log-derivatives of Szeg\H{o} kernels}
\label{sec:EK}

For the DPP of
Peres and Vir\'ag the \textbf{Edelman--Kostlan formula}
\cite{EK95} gives the density of $\cZ_{X_{\D}}$ 
with respect to $m/\pi$ as
\[
\rho^1_{\D, \mathrm{PV}}(z)
=\frac{1}{4} \Delta \log S_{\D}(z,z),
\quad z \in \D,
\]
where $\Delta := 4 \partial_z \partial_{\zbar}$.
Moreover, we have the equality
\begin{equation}
K_{\D}(z, w)= 
\partial_z \partial_{\wbar} 
\log S_{\D}(z, w)=S_{\D}(z, w)^2,
\quad z, w \in \D.
\label{eqn:EK_like}
\end{equation}
As explained above (\ref{eqn:S_K_D}), 
this gives an example of the general formula
\begin{equation}
S_{D}(z, w)^2= K_{D}(z, w), \quad
z, w \in D,
\label{eqn:S_K_D_special}
\end{equation}
which is establised 
for the kernels on any simply connected domain
$D \subsetneq \C$.

Here we have reported our work to
generalize the above to 
a family of GAFs and their zero point processes
on the annulus $\A_q$.
We can prove the following \cite{KS22b}.
\begin{prop}
\label{thm:log_derivative}
For $r > 0$, 
the following equality holds,
\begin{equation}
\partial_z \partial_{\wbar} \log S_{\A_q}(z, w; r)
=\frac{\theta(-r)}{\theta(-r (z \wbar)^2)}
S_{\A_q}(z, w; r z \wbar)^2,
\quad z, w \in \A_q.
\label{eqn:log_deriv_Aq1}
\end{equation}
In particular, 
\begin{equation}
\Delta \log S_{\A_q}(z, z; r)
=4 \frac{\theta(-r)}{\theta(-r|z|^4)}
S_{\A_q}(z, z; r|z|^2)^2,
\quad z \in \A_q.
\label{eqn:log_deriv_Aq2}
\end{equation}
\end{prop}
By comparing the expression (\ref{eqn:density_Aqt}) 
for the density obtained from Theorem \ref{thm:mainA1} 
with (\ref{eqn:log_deriv_Aq2})
in Proposition \ref{thm:log_derivative}, 
we can recover the Edelman--Kostlan formula as follows,
\[
\rho^1_{\A_q}(z; r)
=\frac{\theta(-r)}{\theta(-r|z|^4)}
S_{\A_q}(z, z; r|z|^2)^2
=\frac{1}{4} \Delta \log S_{\A_q}(z, z; r),
\quad z \in \A_q.
\]
However, (\ref{eqn:EK_like}) does not
hold for the weighted Szeg\H{o} kernel for $H^2_r(\A_q)$.
As shown by (\ref{eqn:log_deriv_Aq1}), 
the second log-derivative of $S_{\A_q}(z, w; r)$
cannot be expressed by $S_{\A_q}(z, w; r)$ itself
but a new function $S_{\A_q}(z, w; r z \wbar)$
should be introduced. 

In addition the proportionality between the square of 
the Szeg\H{o} kernel
and the Bergman kernel (\ref{eqn:S_K_D_special}),
which holds as (\ref{eqn:S_K_D_special}) on $\D$, 
is no longer valid for the point processes on $\A_q$
The following are verified \cite{KS22b}.
A CONS for the \textbf{Bergman space on $\A_q$}
is given by 
$\{ \widetilde{e}^{(q)}_n(z)\}_{n \in \Z}$ where we set
\[
\widetilde{e}^{(q)}_n(z)
=\begin{cases}
\displaystyle{
\sqrt{\frac{n+1}{1-q^{2(n+1)}} } z^n,
}
& n \in \Z \setminus \{ -1 \},
\cr
\displaystyle{
\sqrt{\frac{1}{-2 \log q}} z^{-1}},
& n=-1.
\end{cases}
\]
The Bergman kernel of $\A_q$ is then given by
\begin{align}
K_{\A_q}(z, w)
&:= k_{L^2_{\rB}(\A_q)}(z, w)
= \sum_{n \in \Z} \widetilde{e}^{(q)}_n(z)
\overline{ \widetilde{e}^{(q)}_n(w) }
\nonumber\\
&= - \frac{1}{2 \log q} \frac{1}{z \overline{w}}
+\frac{1}{z \overline{w}} 
\sum_{n \in \Z \setminus \{0\}} 
\frac{n}{1-q^{2n}}(z \overline{w})^n,
\quad z, w \in \A_q.
\label{eqn:K_Aq}
\end{align}
The \textbf{Weierstrass $\wp$-function}
is defined by
\begin{align}
\wp(\phi)
&= - \frac{1}{12} + 2 \sum_{n=1}^{\infty}
\frac{q^{2n}}{(1-q^{2n})^2} 
- \sum_{n=-\infty}^{\infty} \frac{e^{\sqrt{-1} \phi} q^{2n}}
{(1-e^{\sqrt{-1} \phi} q^{2n})^2}.
\label{eqn:wp_expansion2}
\end{align}
Then (\ref{eqn:K_Aq}) is expressed as as \cite{Ber70}
\begin{equation}
K_{\A_q}(z, w)
= - \frac{1}{2 \log q} \frac{1}{z \wbar}
- \frac{1}{z \wbar}
\Big(
\wp(\phi_{z \wbar})
+ \frac{P(q)}{12}
\Big), \quad z, w \in \A_q,
\label{eqn:BergmanK3}
\end{equation}
where we have used the notation,
\begin{equation}
z=e^{\sqrt{-1} \phi_z} \iff
\phi_z= - \sqrt{-1} \log z,
\label{eqn:x_phi}
\end{equation}
and 
\[
P(q) = 1-24 \sum_{n=1}^{\infty} \frac{q^{2n}}{(1-q^{2n})^2}.
\]
We can prove the equality \cite{KS22b}
\begin{equation}
S_{\A_q}(z, w)^2
= K_{\A_q}(z, w) + \frac{a(q)}{z \wbar}, 
\quad z, w \in \A_q, 
\label{eqn:SK1}
\end{equation}
where
\begin{equation}
a(q) =
-2 \sum_{n \in \N}
\frac{(-1)^n n q^n}{1-q^{2n}}
+\frac{1}{2 \log q}.
\label{eqn:aq}
\end{equation}

\subsubsection{Hammersley formula and Shirai's Proposition}
\label{sec:Hammersley and Frobenius determinantal formula}

We recall a general formula for correlation functions
of zero point process of a GAF, which
is called the
\textbf{Hammersley formula} in \cite{PV05}, but here we use 
a slightly different expression given 
by Shirai (Proposition 6.1 of \cite{Shi12}). 
Let $\partial_z \partial_{\wbar} 
:= \frac{\partial^2}{\partial z \partial \wbar}$. 

\begin{prop}[Shirai \cite{Shi12}]
\label{thm:correlation_function}
The correlation functions of $\cZ_{X_D}$ of the GAF $X_D$ on 
$D \subsetneq \C$ with
covariance kernel $S_D(z,w)$ are given by 
\[
 \rho_D^n(z_1,\dots, z_n) 
= \frac{ \per_{1 \leq i, j \leq n} \big[(\partial_z \partial_{\wbar}
S_D^{z_1,\dots,z_n})(z_i, z_j) \big]}
{\det_{1 \leq i, j \leq n} \Big[
 S_D(z_i, z_j) \Big]}, 
\quad n \in \N, \quad
z_1, \dots, z_n \in D, 
\]
with respect to a reference measure $\lambda$, 
whenever $\det_{1 \leq i, j \leq n}[S_D(z_i, z_j)]> 0$. 
Here $\{S_D^{z_1,\dots,z_n}(z_i, z_j)\}$ denote 
the conditional kernels. 
\end{prop}
Here we abbreviate 
$\gamma^q_{\{z_{\ell}\}_{\ell=1}^n}$ given by (\ref{eqn:gamma})
to $\gamma^q_n$.
Then (\ref{eqn:MS_general}) gives
\[
S_{\A_q}^{z_1, \dots, z_n}(z, w; r)
=S_{\A_q}(
z, w; r \prod_{\ell=1}^n |z_{\ell}|^2)
\gamma^q_{n}(z)
\overline{ \gamma^q_{n}(w) }
\]
for $z, w, z_1, \dots, z_n \in \A_q$. 
By (\ref{eqn:h_aq_1}) of Lemma \ref{thm:Blaschke}, 
this formula gives
\[
(\partial_z \partial_{\wbar} 
S_{\A_q}^{z_1,\dots,z_n})(z_i, z_j; r) 
= S_{\A_q} \Big(z_i, z_j;  r \prod_{\ell=1}^n |z_{\ell}|^2 \Big) 
{\gamma^q_n}'(z_i)
\overline{{\gamma^q_n}'(z_j)}. 
\]
Therefore, 
Proposition \ref{thm:correlation_function} gives now
\begin{equation}
 \rho^n_{\A_q}(z_1, \dots, z_n; r) 
= \frac{\per_{1 \leq i, j \leq n}
 \left[S_{\A_q} \left(z_i, z_j;  r \prod_{\ell=1}^n |z_\ell|^2 \right) \right]
\prod_{k=1}^n |{\gamma^q_n}'(z_k)|^2 }
{\det_{1 \leq i, j \leq n} 
\Big[
S_{\A_q}(z_i, z_j; r) \Big]}. 
\label{eqn:rho_Z1}
\end{equation}
By (\ref{eqn:haq}) and 
(\ref{eqn:h_aq_1}) and (\ref{eqn:h_aq_4}) 
of Lemma \ref{thm:Blaschke}, 
we see that 
\begin{align*}
\prod_{i=1}^n |{\gamma_n^q}'(z_i)|^2 
&= \prod_{i=1}^n \Big| \Big(
\prod_{1 \leq j \leq n, j \not=i}
h^q_{z_j}(z_i) \Big)
{h^q_{z_i}}'(z_i) \Big|^2
= \prod_{i=1}^n 
\Big| \Big( \prod_{1 \leq j \leq n,  j \not=i}
z_i \frac{\theta(z_j/z_i)}{\theta(\overline{z_j} z_i)} \Big)
\frac{ q_0^2 }{ \theta(|z_i|^2)} \Big|^2 
\nonumber\\
&= \Bigg|
\frac{ q_0^{2n} \prod_{1 \le i < j \le n} z_i \theta(z_j/z_i)
\cdot \prod_{1 \le i' < j' \le n} z_{j'} \theta(z_{i'}/z_{j'}) }
{\prod_{i=1}^n \prod_{j=1}^n \theta(z_i \overline{z_j}) } \Bigg|^2. 
\end{align*}
By (\ref{eqn:theta_inversion}),
$z_i \theta(z_j/z_i)
=z_i ( - z_j/z_i ) \theta(z_i/z_j)
=-z_j \theta(z_i/z_j)$.
Hence this is written as 
\begin{align}
\prod_{i=1}^n |{\gamma_n^q}'(z_{i})|^2 
&= q_0^{4n} \left|
\frac{(-1)^{n(n-1)/2} \big(\prod_{1 \leq i  < j \leq n} z_j \theta(z_i/z_j) \big)^2}
{\prod_{i=1}^n \prod_{j=1}^n \theta(z_i \overline{z_j})} \right|^2
\nonumber\\
&= q_0^{4n}
\Bigg(
\frac{\prod_{1 \leq i < j \leq n} |z_j|^2 
\theta(z_i/z_j, \overline{z_i}/\overline{z_j})}
{\prod_{i=1}^n \prod_{j=1}^n \theta(z_i \overline{z_j})}
\Bigg)^2.
\label{eqn:gamma2}
\end{align}

The following identity is known as 
the \textbf{Frobenius determinantal formula},
which is 
an \textbf{elliptic extension} of the 
\textbf{Cauchy determinantal formula}
due Frobenius (see 
Theorem 1.1 in \cite{KN03}, 
Theorem 66 in \cite{Kra05}, 
Corollary 4.7 in \cite{RS06}),
\begin{align*}
\det_{1 \leq i, j \leq n}
\left[ \frac{\theta(t x_i a_j)}{\theta(t, x_i a_j)} \right]
&=\frac{\theta(t \prod_{k=1}^n x_k a_k)}{\theta(t)}
\frac{\prod_{1 \leq i < j \leq n}
x_j a_j
\theta(x_i/x_j, a_i/a_j)}
{\prod_{i=1}^n \prod_{j=1}^n \theta(x_i a_j)}.
\end{align*}
By (\ref{eqn:S_qt_theta}) in Proposition \ref{thm:S_Aq_theta}, 
we have
\begin{equation}
q_0^{2n}
\frac{\prod_{1 \leq i < j \leq n} |z_j|^2
\theta(z_i/z_j, \overline{z_i}/\overline{z_j})}
{\prod_{i=1}^n \prod_{j=1}^n \theta(z_i \overline{z_j})}
=\frac{\theta(-s)}{\theta( -s \prod_{\ell=1}^n |z_{\ell}|^2)}
\det_{1 \leq i, j \leq n}
\left[S_{\A_q}(z_i, z_j; s) \right], 
\quad \forall s >0.
\label{eqn:Frobenius_formula}
\end{equation}
Then (\ref{eqn:gamma2}) is written as
\begin{align*}
\prod_{i=1}^n |{\gamma_n^q}'(z_i)|^2 
&= \frac{\theta(-r)}{\theta(-r \prod_{\ell=1}^n |z_{\ell}|^2)}
\det_{1 \leq i, j \leq n} [S_{\A_q}(z_i, z_j; r) ]
\nonumber\\
& \quad \times
\frac{\theta(-r \prod_{\ell=1}^n |z_{\ell}|^2)}
{\theta(-r \prod_{\ell=1}^n |z_{\ell}|^4)}
\det_{1 \leq i, j \leq n} 
\Big[S_{\A_q} \Big(z_{i}, z_{j}; 
r \prod_{\ell=1}^n |z_{\ell}|^2 \Big) \Big]
\nonumber\\
&= \frac{\theta(-r)}{\theta(-r\prod_{\ell=1}^n |z_{\ell}|^4)}
\det_{1 \leq i, j \leq n} [S_{\A_q}(z_i, z_j; r) ]
\det_{1 \leq i, j \leq n} 
\Big[S_{\A_q} \Big(z_{i}, z_{j}; r \prod_{\ell=1}^n |z_{\ell}|^2 
\Big) \Big]. 
\end{align*}
Applying the above to (\ref{eqn:rho_Z1}),
the correlation functions
in Theorem \ref{thm:mainA1} are obtained.

\subsection{Unfolded 2-correlation function of
PDPP}
\label{sec:2_corr}

By the determinantal factor in 
$\perdet$ (\ref{eqn:perdet})
the PDPP shall be negatively correlated 
when distances of points
are short in the domain $\A_q$.
The effect of the permanental part
in $\perdet$ will appear in long distances.
In order to clarify this fact, 
we study the two-point correlation
function normalized by the product of
one-point functions,
\begin{align}
g_{\A_q}(z, w; r)
:= \frac{\rho_{\A_q}^2(z,w; r)}
{\rho_{\A_q}^1(z;r) \rho_{\A_q}^1(w;r)},
\quad 
(z, w) \in \A_q^2, 
\label{eqn:unfolded}
\end{align}
where $\rho^1_{\A_q}$ and $\rho^2_{\A_q}$
are explicitly given by (\ref{eqn:density_Aqt})
and (\ref{eqn:rho2_Aqt}), respectively.
This function is simply called an intensity ratio in \cite{PV05},
but here we call it an \textbf{unfolded 2-correlation function}
following a terminology used in 
random matrix theory \cite{For10}.

\begin{figure}[htbp]
\begin{center}
\includegraphics[width=0.4\hsize]{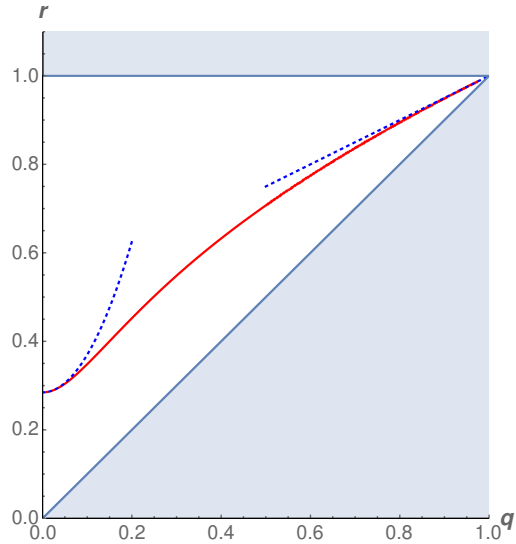}
\end{center}
\caption{
The critical curve $r=r_0(q)$ mentioned in 
Theorem \ref{thm:unfolded} (ii)
is numerically plotted 
(in red) in the fundamental cell
$\Omega:=\{(q, r): q \in (0, 1), q \leq r \leq 1\}$
in the parameter space,
which is located between the diagonal line
$r=q$ 
and the horizontal line $r=1$.
The parabolic curve $r_{\rm c}+ c q^2$
given in \cite{KS22b} 
and the line $1-(1-q)/2$ 
are also dotted,
which approximate $r=r_0(q)$ well for $q \gtrsim 0$ and
$q \lesssim1$, respectively.}
\label{fig:phase_diagram}
\end{figure}

We have proved the following 
(see Fig.\ref{fig:phase_diagram}) \cite{KS22b}.
\begin{thm}
\label{thm:unfolded}
\begin{description}
\item{(i)} \, When $0 < q < 1$, in the short distance,
the correlation is generally repulsive in common with DPPs.

\item{(ii)} \, There exists a critical value 
\[
r_0 =r_0(q) \in (q, 1)
\quad 
\mbox{for each $q \in (0, 1)$}
\]
such that if $r \in (r_0, 1)$ 
\textbf{positive correlation} emerges
between zeros when the distance between them 
is large enough within $\A_q$. 
\item{(iii)} \, The limits 
$g_{\D}(z, w; r) :=\lim_{q \to 0} g_{\A_q}(z, w; r)$,
$z, w \in \D^{\times}$ and
$r_{\rm c}:=\lim_{q \to 0} r_0(q)$
are well-defined, and $r_{\rm c}$ is positive.
When $r \in [0, r_{\rm c})$
all positive correlations vanish in $g_{\D}(z, w; r)$,
while when $r \in (r_{\rm c}, \infty)$
positive correlations can survive. 
\end{description}
\end{thm}

\subsection{Exercises 2}
\label{sec:exercises2}
\subsubsection{Exercise 2.1}
\label{sec:ex2_1}
\begin{description}
\item{(1)} \,
By the definition of the theta function
\begin{align}
\theta(z; p)  &:=(z, p/z; p)_{\infty}
=(z; p)_{\infty} (p/z; p)_{\infty}
\nonumber\\
&=
\prod_{i=1}^{\infty} \Big\{ (1-zp^i) (1-p^{i+1}/z \Big\},
\label{eqn:ex2_theta1}
\end{align}
prove the inversion formula
\begin{equation}
\theta(1/z; p) = - \frac{1}{z} \theta(z; p)
\label{eqn:ex2_theta2}
\end{equation}
and the quasi-periodicity property
\begin{equation}
\theta(pz; p) = -\frac{1}{z} \theta(z; p).
\label{eqn:ex2_theta3}
\end{equation}

\item{(2)} \,
By direct calculation following the
definition (\ref{eqn:ex2_theta1}), verify the equalities,
\begin{align}
&\theta(\zeta; p^3) \theta(\zeta p; p^3) \theta(\zeta p^2; p^3)
= \theta(\zeta; p),
\label{eqn:ex2_theta4}
\\
&
\theta(\zeta; p)
\theta(\zeta \omega_3; p) \theta(\zeta \omega_3^2; p)
= \theta(\zeta^3; p^3), 
\quad \zeta \in \C^{\times}, 
\label{eqn:ex2_theta5}
\end{align}
where $\omega_3$ is a primitive 3rd root of unity.

\item{(3)} \,
Let $x=e^{-2 i a}$, $y=e^{-2 i b}$, $u=e^{-2 i c}$,
and $v=e^{-2 i d}$ in 
the Weierstrass addition formula 
\begin{equation}
\theta(xy, x/y, uv, u/v; p)
- \theta(xv, x/v, uy, u/y; p)
=\frac{u}{y} \theta(yv, y/v, xu, x/u; p).
\label{eqn:ex2_theta6}
\end{equation}
Consider the limit $p \to 0$ and
derive the following trigonometric equality,
\begin{align}
&\sin(a+b) \sin(a-b) \sin(c+d) \sin(c-d)
\nonumber\\
& \quad 
-\sin(a+d) \sin(a-d) \sin(c+b) \sin(c-b)
\nonumber\\
& \qquad
=\sin(b+d) \sin(b-d) \sin(a+c) \sin (a-c).
\label{eqn:ex2_theta7}
\end{align}
Show that (\ref{eqn:ex2_theta7}) is readily verified if
we use the following addition formulas of 
trigonometric functions,
\begin{equation}
\sin(a \pm b)=\sin a \cos b \pm \cos a \sin b.
\label{eqn:ex2_theta8}
\end{equation}
\end{description}

\subsubsection{Exercise 2.2}
\label{sec:ex2_2}
Consider the weighted Szeg\H{o} kernel of
$\A_q:=\{z \in \C: q < |z| < 1\}$, 
$q \in (0,1)$ expressed using the theta function
\begin{equation}
S_{\A_q}(z, w; r)
= \frac{q_0^2 \theta(-r z \wbar)}
{\theta(-r, z \wbar)},
\quad z, w \in \A_q, \quad r >0,
\label{eqn:ex2_S1}
\end{equation}
where 
\begin{equation}
q_0:=\prod_{n \in \N}(1-q^{2n})=(q^2; q^2)_{\infty}.
\label{eqn:ex2_S2}
\end{equation}
Schottky's theorem asserts that 
the group of conformal transformations
from $\A_q$ to itself 
is generated by the rotations 
\begin{equation}
R_{\theta}(z) := e^{\sqrt{-1} \theta} z, \quad
\theta \in [0, 2 \pi),
\label{eqn:ex2_S3}
\end{equation}
and the $q$-inversions 
\begin{equation}
T_q(z) := \frac{q}{z}. 
\label{eqn:ex2_S4}
\end{equation}
\begin{description}
\item{(1)} \,
Show that the rotational invariance is obvious. 

\item{(2)} \,
Prove the following equality,
\begin{equation}
\sqrt{T_q'(z)} \overline{\sqrt{T_q'(w)}}
S_{\A_q}(T_q(z), T_q(w) ; r) 
= \frac{q}{r} S_{\A_q} (z, w; q^2/r ).
\label{eqn:ex2_S5}
\end{equation}
Note that if $r=q$ the weighted Szeg\H{o} kernel
is reduced to the original Szeg\H{o} kernel,
$S_{\A_q}(z, w)=S_{\A_q}(z, w; q)$.
In this case the above becomes
\begin{equation}
\sqrt{T_q'(z)} \overline{\sqrt{T_q'(w)}}
S_{\A_q}(T_q(z), T_q(w)) 
= S_{\A_q} (z, w).
\label{eqn:ex2_S6}
\end{equation}
Hence we can conclude that the
original Szeg\H{o} kernel,
$S_{\A_q}(z, w)=S_{\A_q}(z, w; q)$,
is invariant under conformal transformations
from $\A_q$ to itself.
\end{description}

\subsubsection{Exercise 2.3}
\label{sec:ex2_3}
Consider the DPP of Peres and Vir\'ag, $\cZ_{X_{\D}}$,
whose correlation kernel with respect to $m/\pi$ is given by the 
Bergman kernel of $\D$, $K_{\D}$;
\begin{align}
&\rho_{\D, \mathrm{PV}}^n(z_1, \dots, z_n)
=\det_{1 \leq i, j \leq n} [K_{\D}(z_i, z_j)],
\quad n \in \N, \quad z_1, \dots, z_n \in \D,
\nonumber\\
& \mbox{with} \quad
K_{\D}(z, w)=\frac{1}{(1-z \wbar)^2},
\quad z, w \in \D.
\label{eqn:ex2_PV1}
\end{align}
\begin{description}
\item{(1)} \,
Prove the following equality, 
\begin{equation}
\rho^1_{\D, \mathrm{PV}}(z)
=\frac{1}{4} \Delta \log S_{\D}(z,z),
\quad z \in \D,
\label{eqn:ex2_PV2}
\end{equation}
where $\Delta := 4 \partial_z \partial_{\zbar}$ and 
$S_{\D}$ denotes the Szeg\H{o} kernel of $D$,
\begin{equation}
S_{\D}(z, w) = \frac{1}{1-z \wbar}, \quad z, w \in \D.
\label{eqn:ex2_PV3}
\end{equation}
This equality is called the Edelman--Kostlan formula.
\item{(2)} \,
Show that the following general equalities are established,
\begin{equation}
K_{\D}(z, w)= 
\partial_z \partial_{\wbar} 
\log S_{\D}(z, w)=S_{\D}(z, w)^2,
\quad z, w \in \D.
\label{eqn:ex2_PV4}
\end{equation}
\end{description}

\subsubsection{Exercise 2.4}
\label{sec:ex2_4}
Consider the GAF and the permanenta-determinantal point process
on $\A_q$ discussed there.
\begin{description}
\item{(1)} \,
Prove the equalities,
\begin{equation}
\partial_z \partial_{\wbar} \log S_{\A_q}(z, w; r)
=\frac{\theta(-r)}{\theta(-r (z \wbar)^2)}
S_{\A_q}(z, w; r z \wbar)^2,
\quad z, w \in \A_q.
\label{eqn:ex2_Aq1}
\end{equation}
In particular, 
\begin{equation}
\Delta \log S_{\A_q}(z, z; r)
=4 \frac{\theta(-r)}{\theta(-r|z|^4)}
S_{\A_q}(z, z; r|z|^2)^2,
\quad z \in \A_q.
\label{eqn:ex_Aq2}
\end{equation}
\item{(2)} \,
The correlation function of the permanental-detaeminantal
point process on $\A_q$ is generally given by
\begin{equation}
\rho^{n}_{\A_q}(z_1,\dots,z_n; r) 
=
\frac{\theta(-r)}{\theta( -r \prod_{k=1}^n |z_k|^4)}
\perdet_{1 \leq i, j \leq n}
\Big[
S_{\A_q} \Big(z_{i}, z_{j}; r \prod_{\ell=1}^n |z_{\ell}|^2 \Big) 
\Big]
\label{eqn:ex_Aq3}
\end{equation}
for every $n \in \N$ and
$z_1, \dots, z_n \in \A_q$
with respect to $m/\pi$.
Show that the Edelman--Kostran formula is also 
satisfied.
\end{description}

\SSC
{Multiple Schramm--Loewner Evolution/Gaussain 
Free Field Coupling}
\label{sec:SLE_GFF}
\subsection{Multiple Schramm--Loewner 
Evolution} \label{sec:multiple_SLE}
\subsubsection{Loewner equations for single-slit and multi-slit}
\label{sec:chordal_LE}

Let $D$ be a simply connected domain in $\C$ which does not
complete the plane; $D \subsetneq \C$.
Its boundary is denoted by $\partial D$.
We consider a slit in $D$, which is defined as a trace
$\eta = \{\eta(t) : t \in (0, \infty) \}$
of a \textbf{simple curve} 
$\eta(t) \in D, 0 < t < \infty$; 
\[
\eta(s) \not= \eta(t) \quad \mbox{for $s \not=t$}.
\]
We assume that the initial point of the slit is located in $\partial D$,
$^{\exists}\eta(0) := \lim_{t \to 0} \eta(t) \in \partial D$.
Let 
\[
\eta(0,t]:=\{\eta(s) : s \in (0, t]\}
\quad \mbox{and} \quad
D^{\eta}_t := D \setminus \eta(0, t], t \in (0, \infty).
\]
The Loewner theory \cite{Low23} describes a slit $\eta$
by encoding the curve into a time-dependent
analytic function $g_{D^{\eta}_t} : t \in (0, \infty)$
such that 
\[
g_{D^{\eta}_t} : \mbox{conformal map} \, \, 
D^{\eta}_t \to D,
\quad t \in (0, \infty).
\]

By the \textbf{Riemann mapping theorem},
for $D \subsetneq \C$ and a point $z_0 \in D$,
there exists a unique analytic function
$\varphi(z)$ in $D$, specified by
$\varphi(z_0)=0, \varphi'(z_0) > 0$, such that
\[
\varphi : \mbox{conformal map} \, \, 
D \to \D,
\]
where $\D$ denotes a unit disk; $\D:=\{z \in \C: |z| < 1\}$
(see (\ref{eqn:RiemannA}) in the previous section).
Loewner gave a differential equation for $g_{D^{\eta}_t} $ 
in the case $D=\D$, which is called the Loewner equation \cite{Low23}.

Since a special case of the 
\textbf{M\"obius transformation}
\[
\sm(z) := \sqrt{-1} \frac{\alpha-z}{\alpha+z}, \quad |\alpha|=1,
\]
maps $\D$ to the \textbf{upper half complex plane}
$\HH:=\{z \in \C: \Im z > 0\}$
with $\sm(0)=\sqrt{-1}, \sm(\infty) = -\sqrt{-1}$,
we can apply the Loewner theory to the case
with $D=\HH$, in which  \cite{KSS68}
\begin{align*}
&\mbox{$\eta:=\{\eta(t) : t \in (0, \infty)\}$ is a \textbf{simple curve}},
\nonumber\\
&\eta(0) := \lim_{t \to 0} \eta(t) = 0 \in \R,
\nonumber\\
&\eta(0, t] \subset \HH, \quad \forall t \in (0, \infty),
\nonumber\\
& \lim_{t \to \infty} \eta(t)=\infty.
\end{align*}
For each time  $t \in (0, \infty)$,
\begin{equation}
\HH^{\eta}_t := \HH \setminus \eta(0, t]
\label{eqn:Heta1}
\end{equation}
is a simply connected
domain in $\C$ and there exists a unique analytic function
$g_{\HH^{\eta}_t}$ such that
\[
g_{\HH^{\eta}_t} : \mbox{conformal map} \, \, 
\HH^{\eta}_t \to \HH,
\]
which satisfies the condition
\[
g_{\HH^{\eta}_t}(z)=z+ \frac{c_t}{z} + \rO(|z|^{-2})
\quad \mbox{as $z \to \infty$}
\]
for some $c_t > 0$, 
in which the coefficient of $z$ is unity and no constant term appears.
This is called the
\textbf{hydrodynamic normalization}.
The coefficient $c_t$ gives the 
\textbf{half-plane capacity} of 
$\eta(0, t]$ and denoted by 
$\hcap(\eta(0, t])$.
The following has been shown (see \cite{KSS68,Law05}).

\begin{thm}
\label{thm:LE}
Let $\eta$ be a slit in $\HH$ for which 
the parameterization by $t$ is arranged so that
\[
c_t=\hcap(\eta(0, t])=2t, \quad t \in (0, \infty).
\]
Then there exists a unique continuous 
\textbf{driving function} $V(t) \in \R, t \in (0, \infty)$ such that the 
solution $g_t$ of the differential equation
\begin{equation}
\frac{d g_t(z)}{dt} = \frac{2}{g_t(z)-V(t)}, \quad
t \geq 0, 
\label{eqn:LE1}
\end{equation}
under the initial condition
\[
g_0(z)=z \in \HH
\]
gives $g_t=g_{\HH^{\eta}_t}, t \in (0, \infty)$.
\end{thm}

The equation (\ref{eqn:LE1}) is called
the \textbf{chordal Loewner equation}.
Note that at each time $t \in (0, \infty)$,
the \textbf{tip of slit} $\eta(t)$ and the value of $V(t)$
satisfy the following relations,
\begin{equation}
V(t)=\lim_{\substack{z \to 0, \cr \eta(t)+z \in \HH^{\eta}_t}}
g_{\HH^{\eta}_t}(\eta(t)+z)
\quad \Longleftrightarrow \quad
\eta(t)=\lim_{\substack{z \to 0, \cr z \in \HH}} 
g_{\HH^{\eta}_t}^{-1}(V(t)+z),
\quad t \geq 0.
\label{eqn:Ut1}
\end{equation}
Moreover, 
$\displaystyle{V(t)=\lim_{s < t, s \to t} g_{\HH^{\eta}_s}(\eta(t))}$
and $t \mapsto V(t)$ is continuous
(see, for instance, Lemma 4.2 in \cite{Law05}).
We write 
\begin{equation}
g_{\HH^{\eta}_t}(\eta(t))=V(t) \in \R, 
\quad \Longleftrightarrow \quad
\eta(t)=
g_{\HH^{\eta}_t}^{-1}(V(t)) \in \partial \HH^{\eta}_t,
\quad t \geq 0, 
\label{eqn:tip}
\end{equation}
in the sense of (\ref{eqn:Ut1}). 

\vskip 0.3cm
\begin{example}
\label{thm:SLE0}
When the driving function is identically zero;
$V(t) \equiv 0, t \in (0, \infty)$,
the chordal Loewner equation 
$dg_{\HH^{\eta}_t}(z)/dt=2/g_{\HH^{\eta}_t}(z), t \geq 0$
is solved under the initial condition
$g_{\HH^{\eta}_0}(z)=z \in \HH$ as
$g_{\HH^{\eta}_t}(z)^2=4t + z^2, t \geq 0$.
In this simple case, (\ref{eqn:Ut1}) gives
$\eta(t)=2 \sqrt{-1} t^{1/2}, t \geq 0$.
That is, the slit $\eta(0, t], t >0$
is a straight line along the imaginary axis 
starting from the origin, 
$\eta(0)=\lim_{t \to 0} \eta(t)=0$,
and growing upward as time $t$ is passing.
\end{example}
\vskip 0.3cm
\begin{example}
\label{thm:SLE0b}
The above example can be extended by introducing
one parameter $\alpha \in (0, 1)$ as follows.
Let
$\kappa=\kappa(\alpha)=4(1-2 \alpha)^2/\{\alpha(1-\alpha)\}$,
and consider the case such that
\[
V(t)=\begin{cases}
\sqrt{\kappa t}, & \mbox{if $\alpha \leq 1/2$}, \cr
-\sqrt{\kappa t}, & \mbox{if $\alpha > 1/2$}.
\end{cases}
\]
In this case, the inverse of $g_t$ is solved as
$g_{\HH^{\eta}_t}^{-1}(z) = \left( z + 2 \sqrt{\frac{\alpha}{1-\alpha}} \sqrt{t} \right)^{1-\alpha}
\left(z - 2 \sqrt{\frac{1-\alpha}{\alpha}} \sqrt{t} \right)^{\alpha}$,
and the slit is obtained as
\[
\eta(t)=g_{\HH^{\eta}_t}^{-1}(V(t))
=2 \left( \frac{1-\alpha}{\alpha} \right)^{1/2-\alpha}
e^{\sqrt{-1} \alpha \pi} t^{1/2},
\quad t \geq 0.
\]
The slit grows from the origin along a straight line in $\HH$
which makes an angle $\alpha \pi$ with respect to the positive
direction of the real axis.
When $\alpha=1/2$, this is reduced to the result mentioned
in Example \ref{thm:SLE0}.
More detail for this example, see Example 4.12 in \cite{Law05}
and Section 2.2 in \cite{Kat16_Springer}.
\end{example}

Theorem \ref{thm:LE} can be extended 
to the situation such that $\eta$ in $\HH$ is given by
a multi-slit \cite{RS17,KK21b}.
Let $N \in \N :=\{1,2, \dots\}$ and assume that we have $N$ slits 
$\eta_i =\{\eta_i(t): t \in (0, \infty)\} \subset \HH$, 
$i=1, \dots, N$, 
which are simple curves, disjoint with each other, 
$\eta_i \cap \eta_j = \emptyset, i \not= j$,
starting from $N$ distinct points 
$\lim_{t \to 0} \eta_i(t) =: \eta_i(0)$ on $\R$;
$\eta_1(0) < \cdots < \eta_N(0)$, 
and all going to infinity; $\lim_{t \to \infty} \eta_i(t)=\infty$, 
$i=1, \dots, N$. 
A \textbf{multi-slit} is defined as a union of them, 
$\bigcup_{i=1}^N \eta_i$, and 
\begin{equation}
\HH^{\eta}_t := \HH \setminus \bigcup_{i=1}^N \eta_i(0, t]
\label{eqn:Heta2}
\end{equation}
for each $t > 0$ with $\HH^{\eta}_0 := \HH$.
For each time $t \in (0, \infty)$, 
$\HH^{\eta}_t$
is a simply connected domain in $\C$ and then
there exists a unique analytic function $g_{\HH^{\eta}_t}$ such that
\[
g_{\HH^{\eta}_t} : \mbox{conformal map} \, \, 
\HH^{\eta}_t \to \HH,
\]
satisfying the hydrodynamic normalization condition
\[
g_{\HH^{\eta}_t}(z)=z+ 
\frac{\hcap(\bigcup_{i=1}^N \eta_i(0, t])}{z} + \rO(|z|^{-2})
\quad \mbox{as $z \to \infty$}.
\]
\begin{thm}
\label{thm:mLE}
For $N \in \N$, 
let $\bigcup_{i=1}^N \eta_i$ be a multi-slit in $\HH$ such that
\[
\hcap \Big(\bigcup_{i=1}^N \eta(0, t] \Big)=2 N t, 
\quad t \in (0, \infty).
\]
Then there exists a set of weight functions 
$w_i(t) \geq 0, t \geq 0, i=1, \dots, N$
satisfying $\sum_{i=1}^N w_i(t)=1, t \geq 0$ and an
$N$-variate continuous driving function 
$\V(t)=(V_1(t), \dots, V_N(t)) \in \R^N, t \in (0, \infty)$ such that the 
solution $g_t$ of the differential equation
\begin{equation}
\frac{d g_t(z)}{dt} = \sum_{i=1}^N \frac{2 N w_i(t)}{g_t(z)-V_i(t)}, \quad
t \geq 0, \quad 
g_0(z)=z,
\label{eqn:mLE1}
\end{equation}
gives $g_t=g_{\HH^{\eta}_t}, t \in (0, \infty)$.
\end{thm}

Roth and Schleissinger \cite{RS17} called (\ref{eqn:mLE1}) 
the \textbf{Loewner equation for the multi-slit} 
$\bigcup_{i=1}^N \eta_i$.
Similar to (\ref{eqn:Ut1}), the following relations hold,
\begin{align}
V_i(t) =
\lim_{\substack{z \to 0, \cr \eta_i(t)+z \in \HH^{\eta}_t}}
g_{\HH^{\eta}_t}(\eta_i(t)+z)
\quad &\Longleftrightarrow \quad
\eta_i(t)=\lim_{\substack{z \to 0, \cr z \in \HH}} g_{\HH^{\eta}_t}^{-1}(V_i(t)+z),
\quad i=1, \dots, N, \quad t \geq 0,
\label{eqn:Ut2}
\end{align}
and we write for the \textbf{multiple tips},
$\eta_i(t), i=1, \dots, N$, $t \geq 0$,
\[
g_{\HH^{\eta}_t}(\eta_i(t))=V_i(t) \in \R,
\quad i=1, \dots, N, \quad t \geq 0
\]
in the sense of (\ref{eqn:Ut2}).

\subsubsection{Schramm--Loewner evolution 
with parameter $\kappa$ (SLE$_{\kappa}$)}
\label{sec:SLE}

So far we have considered the problem such that, 
given time-evolution of a single slit $\eta(0, t], t \geq 0$
or a multi-slit $\bigcup_{i=1}^N \eta(0, t], t \geq 0$ in $\HH$,
time-evolution of the conformal map from $\HH^{\eta}_t$ to $\HH$,
$t \geq 0$ is asked.
The answers are given by the solution of the Loewner equation
(\ref{eqn:LE1}) in Theorem \ref{thm:LE} for a single slit
and by the solution of the multiple Loewner equation
(\ref{eqn:mLE1}) in Theorem \ref{thm:mLE} for a multi-slit,
which are driven by a function
$(V(t))_{t \geq 0}$ and by a multi-variate function
$\V(t)=(V_1(t), \dots, V_N(t)) \in \R^N, t \geq 0$, 
respectively. The both processes are defined in $\R$
and deterministic:
\begin{align*}
\mbox{single slit $\eta(0, t] \in \HH, t \geq 0$}
\quad &\Longrightarrow \quad
\mbox{driving function $(V(t))_{t \geq 0}$ on $\R$}
\nonumber\\
\mbox{multi-slit $\bigcup_{i=1}^N \eta(0, t] \in \HH, t \geq 0$}
\quad &\Longrightarrow \quad
\mbox{multi-variate driving function $(\V(t))_{t \geq 0}$ on $\R^N$}
\end{align*}

For $\HH$ with a single slit, Schramm considered
an \textbf{inverse problem in a probabilistic setting} \cite{Sch00}.
He first asked a suitable family of driving stochastic processes
$(Y(t))_{t \geq 0}$ on $\R$.
Then he asked the probability law of a random slit in $\HH$,
which will be determined by the relations (\ref{eqn:Ut1})
from $(Y(t))_{t \geq 0}$ and the solution 
$g_t=g_{\HH^{\eta}_t}, t \geq 0$
of the Loewner equation (\ref{eqn:LE1}):
\begin{align*}
\mbox{\textbf{random} curve $\eta(0, t] \in \HH, t \geq 0$}
\quad &\Longleftarrow \quad
\mbox{driving \textbf{stochastic} process 
$(Y(t))_{t \geq 0}$ on $\R$}
\end{align*}
Schramm argued that \textbf{conformal invariance} implies that the
driving process $(Y(t))_{t \geq 0}$ should be 
a \textbf{continuous Markov process} which has in a particular 
parameterization
\textbf{independent increments}.
Hence $Y(t)$ can be a constant time change of
a one-dimensional standard Brownian motion
$(B(t))_{t \geq 0}$, and it is expressed as
\begin{equation}
(\sqrt{\kappa} B(t))_{t \geq 0} \law= (B(\kappa t))_{t \geq 0}
\quad \mbox{with a parameter $\kappa >0$}.
\label{eqn:kappa1}
\end{equation}
The solution of the Loewner equation 
driven by $Y(t)=\sqrt{\kappa} B(t), t \geq 0$, 
\begin{equation}
\frac{d g_{\HH^{\eta}_t}(z)}{dt} 
= \frac{2}{g_{\HH^{\eta}_t}(z)-\sqrt{\kappa} B(t)},\quad t \geq 0, 
\quad g_{\HH^{\eta}_0}(z)=z \in \HH,
\label{eqn:SLE1}
\end{equation}
is called the
\textbf{chordal Schramm--Loewner evolution} (\text{chordal SLE})
with parameter $\kappa >0$ and is written as
\textbf{SLE$_{\kappa}$} for short. 

The following was proved by Lawler, Schramm, and Werner \cite{LSW04}
for $\kappa=8$ and by Rohde and Schramm \cite{RS05} for $\kappa \not=8$.
\begin{prop} 
\label{thm:SLE_curve}
By (\ref{eqn:Ut1}), a chordal SLE$_{\kappa}$
$g_{\HH^{\eta}_t}, t \in (0, \infty)$ determines 
a continuous curve $\eta=\{\eta(t): t \in (0, \infty)\} \subset \HH$
such that $\eta(0) :=\lim_{t \downarrow 0} \eta(t)=0$ and 
$\lim_{t \to \infty} |\eta(t)|=\infty$
with probability one.
\end{prop}

The continuous curve $\eta$ determined by an SLE$_{\kappa}$ 
is called an \textbf{SLE$_{\kappa}$ curve}
(or \textbf{SLE$_{\kappa}$ trace})
The probability law of an SLE$_{\kappa}$ curve depends on 
$\kappa$. As a matter of fact, SLE$_{\kappa}$ curve
becomes self-intersecting and can touch the real axis $\R$
when $\kappa > 4$,
so it is no more a slit, since a slit has been defined as
a trace of a continuous simple curve.
When $\kappa > 4$, the domain
$\HH \setminus \eta(0, t]$ is divided into 
many components, only one of which 
is unbounded. 
So here we change the definition (\ref{eqn:Heta1}) of
$\HH^{\eta}_t, t \geq 0$ as follows,
\begin{equation}
\mbox{
$\HH^{\eta}_t :=$ the unbounded component
of $\HH \setminus \eta(0, t], \quad t \geq 0$}. 
\label{eqn:HetaB1}
\end{equation}
Then 
\[
g_{\HH^{\eta}_t}(z)
 : \mbox{conformal map} \, \, 
\HH^{\eta}_t \to \HH, \quad t \geq 0.
\]
We also define
\begin{equation}
K^{\eta}_t := \overline{\HH \setminus \HH^{\eta}_t},
\quad t \geq 0,
\label{eqn:SLE_hull1}
\end{equation}
and call it the \textbf{SLE hull}.

\begin{figure}[t]
\includegraphics[scale=.80]{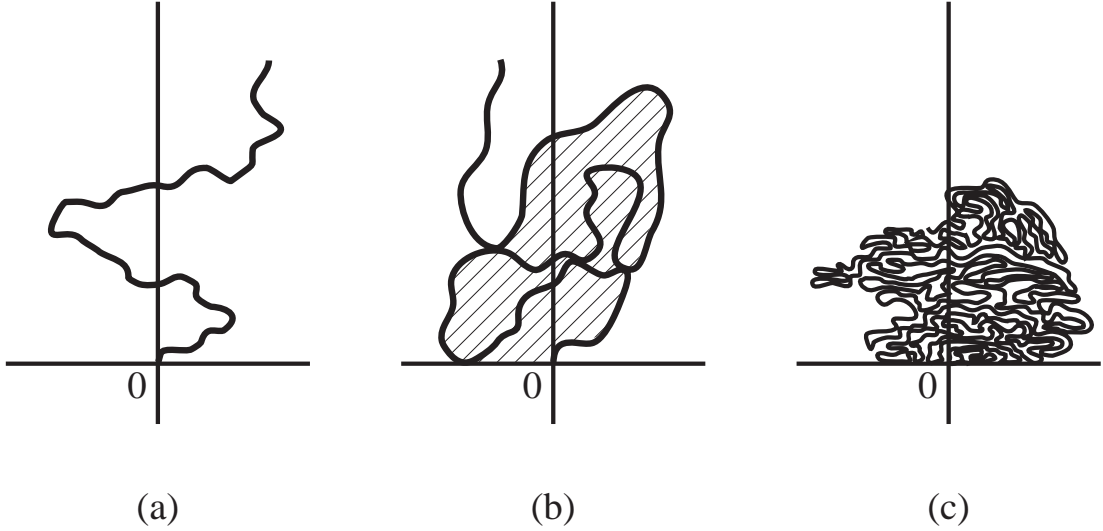}
\caption{Schematic pictures of SLE$_{\kappa}$ curves in
  (a) Phase 1 ($0 < \kappa \leq \kappa_{\rm c}=4$),
  (b) Phase 2 ($\kappa_{\rm c}=4 < \kappa < \overline{\kappa}_{\rm c}=8$), and
  (c) Phase 3 ($\kappa \geq \overline{\kappa}_{\rm c}=8$).
}
\label{fig:SLEphases}       
\end{figure}

There are \textbf{three phases of an SLE$_{\kappa}$ curve}
as shown by the follows (see Fig.\ref{fig:SLEphases}). 

\begin{prop}
\label{thm:SLE_3_phases}
There are two critical values of $\kappa$;
\[
\kappa_{\rm c}=4 \quad \mbox{and} \quad
\overline{\kappa}_{\rm c}=8. 
\]
\begin{description}
\item{\rm (a)} \quad
If $0 < \kappa \leq \kappa_{\rm c}=4$, then the SLE$_{\kappa}$ curve
is \textbf{simple} and $\eta=\eta(0, \infty) \subset \HH$
with probability one.

\item{\rm (b)} \quad
If $\kappa_{\rm c}=4 < \kappa < \overline{\kappa}_{\rm c}=8$,
the SLE$_{\kappa}$ curve is 
\textbf{self-intersecting}, $\eta \cap \R \not= \emptyset$,
and touch the real $\R$ with positive probability. 
Then
\[
\bigcup_{t \in [0, \infty)} K^{\eta}_t = \overline{\HH},
\]
but $\eta[0, \infty) \cap \HH \not= \HH$.

\item{\rm (c)} \quad
If $\kappa \geq \overline{\kappa}_{\rm c}=8$, 
then $\eta$ is a \textbf{space-filling} curve; 
\[
\eta[0, \infty) = \overline{\HH}.
\]
\end{description}
\end{prop}

\begin{figure}[t]
\begin{center}
\includegraphics[scale=.50]{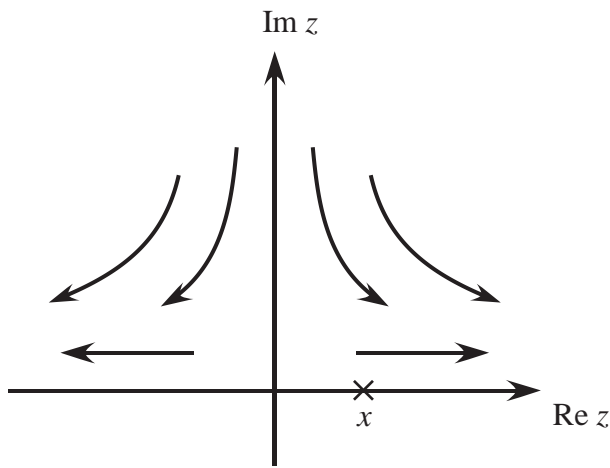}
\caption{A schematic picture of 
\textbf{complexificated Bessel flow}
on $\overline{\HH} \setminus \{0\}$ for $D > 2$.
}
\end{center}
\label{fig:SLEflow}       
\end{figure}
\begin{rem}
\label{thm:remark2_1}
Consider BES$_{D}$ introduced in Section \ref{sec:BES}.
We set 
$Z^{z}(t)=X^{z}(t)+\sqrt{-1} Y^{z}(t) \in 
\overline{\HH} \setminus \{0\}, \ t \geq 0$
and complexificate (\ref{eqn:BESeq1}) as
\begin{equation}
dZ^z(t)=dB(t)+ \frac{D-1}{2} \frac{dt}{Z^{z}(t)}
\label{eqn:SLE_BES1}
\end{equation}
with the initial condition
$$
Z^{z}(0)=z=x+\sqrt{-1} y \in \overline{\HH} \setminus \{0\}.
$$
The crucial point of this 
\textbf{complexification of
Bessel flow} is that the BM remains real,
$B(t) \in \R, \ t \geq 0$.
Then, there is an asymmetry between the real part
and the imaginary part of the flow in $\HH$,
\begin{eqnarray}
\label{eqn:SLE3a}
&& dX^{z}(t)=dB(t)+\frac{D-1}{2}
\frac{X^{z}(t)}{(X^z(t))^2+(Y^z(t))^2} dt, \\
\label{eqn:SLE3b}
&& dY^{z}(t)=-\frac{D-1}{2}
\frac{Y^{z}(t)}{(X^z(t))^2+(Y^z(t))^2} dt.
\end{eqnarray}
Assume $D >1$.
Then as indicated by the minus sign in RHS of (\ref{eqn:SLE3b}),
the flow is downward in $\overline{\HH}$.
If the flow goes down and arrives at the real axis,
the imaginary part vanishes, $Y^{z}(t)=0$,
then equation (\ref{eqn:SLE3a}) is reduced to be
the same equation as equation (\ref{eqn:BESeq1})
for the BES$_{D}$, which is now
considered for $\R \setminus \{0\} = \R_+ \cup \R_-$.
If $D > D_{\rm c}=2$, by Theorem \ref{thm:BESD1} (ii), the flow
on $\R \setminus \{0\}$ is asymptotically outward,
$X^{z}(t) \to \pm \infty$ as $t \to \infty$.
Therefore, the flow on $\overline{\HH}$ will be described
as shown by Fig. \ref{fig:SLEflow}.
The behavior of flow should be, however, more complicated
when $\overline{D}_{\rm c}=3/2 < D < D_{\rm c}$
and $1 < D < \overline{D}_{\rm c}=3/2$.
For $z \in \overline{\HH} \setminus \{0\}, \ t \geq 0$, put
\begin{equation}
\widehat{g}_t(z) :=Z^{z}(t)+B(t).
\label{eqn:SLE4}
\end{equation}
Then, Eq. (\ref{eqn:SLE_BES1}) is rewritten for 
$\widehat{g}_t(z)$ as
\begin{equation}
\frac{\partial \widehat{g}_t(z)}{\partial t}=\frac{D-1}{2}
\frac{1}{\widehat{g}_t(z)-B(t)}, 
\quad t \geq 0.
\label{eqn:SLE5}
\end{equation}
Put
\begin{equation}
\kappa = \frac{4}{D-1}
\quad \Longleftrightarrow \quad
D=1 + \frac{4}{\kappa},
\label{eqn:D_kappa}
\end{equation}
and set
\[
g_t(z)=\widehat{g}_{\kappa t}(z)
\]
in (\ref{eqn:SLE5}).
Then 
we have the equation in the form 
\begin{equation}
\frac{\partial g_t(z)}{\partial t}=
\frac{2}{g_t(z)-\sqrt{\kappa}B(t)},
\quad t \geq 0.
\label{eqn:Lc1}
\end{equation}
This is equal to the Schramm--Loewner evolution 
(\ref{eqn:SLE1}). 
As a matter of fact, identification of the SLE$_{\kappa}$
as a complexification of BES$_{1+4/\kappa}$ gives
the proof of Proposition \ref{thm:SLE_3_phases},
in which simple correspondence between two sets of
critical values;
\begin{equation}
D_{\rm c}=2 \, \iff \, \kappa_{\rm c}=4,
\qquad
\overline{D}_{\rm c}=\frac{3}{2} \,
\iff \, \overline{\kappa}_{\rm c}=8.
\label{eqn:critical_values}
\end{equation}
More detailed description of
the probability laws of an SLE$_{\kappa}$ curves
at special values of $\kappa$, see, for instance
\cite{Law05,Kat16_Springer,Kem17}.
\end{rem}

The highlight of the theory of SLE would be
that, if the value of $D$ is properly chosen,
the probability law of $\gamma$ realizes
that of the scaling limit of important lattice paths
studied in a 
\textbf{statistical mechanics model} exhibiting
\textbf{critical phenomena} 
or describing interesting 
\textbf{fractal geometry}
defined on
an infinite discrete lattice.
The following is a list of the correspondence (up to a conjecture) 
between the SLE$_{\kappa}$ paths with
specified values of $\kappa$, and
the names of lattice paths
(with the names of models studied
in statistical mechanics and fractal physics),
whose scaling limits are described by the SLE$_{\kappa}$ paths.
(See also \cite{DS12,Dum13}.)
\begin{eqnarray*}
\mbox{SLE$_2$}
&\Longleftrightarrow& 
\mbox{loop-erased random walk} \, \cite{LSW04}
\nonumber\\
\mbox{SLE$_{8/3}$}
&\Longleftrightarrow& 
\mbox{self-avoiding walk [conjecture]}
\nonumber\\
\mbox{SLE$_3$}
&\Longleftrightarrow& 
\mbox{Ising interface (critical Ising model)} \, \cite{CS12,CDHKS13}
\nonumber\\
\mbox{SLE$_4$}
&\Longleftrightarrow& 
\mbox{random contour curve (Gaussian free surface model)} \, \cite{SS05}
\nonumber\\
\mbox{SLE$_{16/3}$}
&\Longleftrightarrow& 
\mbox{FK--Ising interface (critical Ising model)} \, \cite{Smi10,CDHKS13}
\nonumber\\
\mbox{SLE$_6$}
&\Longleftrightarrow& 
\mbox{percolation exploration process 
(critical percolation model)} \,
\cite{Smi01}
\nonumber\\
\mbox{SLE$_8$}
&\Longleftrightarrow& 
\mbox{random Peano curve (uniform spanning tree)} \, \cite{LSW04}
\end{eqnarray*}

Moreover, Beffara \cite{Bef08}
determined the Hausdorff dimensions $d^{\rm H}_{\kappa}$
of the SLE$_{\kappa}$ curves as
\begin{equation}
d^{\rm H}_{\kappa}=
\begin{cases}
\displaystyle{1 + \frac{\kappa}{8}},
& \kappa \in (0, 8)
\cr
2,
& \kappa \geq 8.
\end{cases}
\label{eqn:Beffara}
\end{equation}

The relationship between the SLE$_{\kappa}$
and the \textbf{conformal field theory} 
has been clarified \cite{BB02}.
The \textbf{central charge} $c$ and
the \textbf{scaling dimension} 
(the highest weight) of representation
of the Virasoro algebra are related with the
parameter $\kappa$ of SLE$_{\kappa}$ as
\begin{equation}
c=c_{\kappa} :=\frac{(6-\kappa)(3\kappa-8)}{2 \kappa},
\quad
h=h_{\kappa} :=\frac{6-\kappa}{2 \kappa}.
\label{eqn:CFT}
\end{equation}

\subsubsection{Multiple SLE}
\label{sec:mSLE}

For simplicity, we assume that 
$w_i(t) \equiv 1/N, t \geq 0, i=1, \dots, N$
in (\ref{eqn:mLE1}) in Theorem \ref{thm:mLE}.
Then 
the Loewner equation for the multi-slit in $\HH$ is written as
\begin{equation}
\frac{d g_{\HH^{\eta}_t}(z)}{dt}
= \sum_{i=1}^N \frac{2}{g_{\HH^{\eta}_t}(z)-Y_i(t)}, \quad t \geq 0, \quad 
g_{\HH^{\eta}_0}(z)=z \in \HH.
\label{eqn:mSLE1}
\end{equation}
Then we ask what is the suitable family of driving
stochastic processes of $N$ particles on $\R$,
$\Y(t)=(Y_1(t), \dots, Y_N(t)), t \geq 0$
\cite{Car03,BBK05,KL07,Gra07};
\begin{align*}
&\mbox{\textbf{multiple} SLE$_{\kappa}$ curves
$\displaystyle{\bigcup_{i=1}^N \eta_i(0, t] \in \HH, t \geq 0}$}
\nonumber\\
& \qquad 
\qquad \Longleftarrow \quad
\mbox{driving \textbf{many-particle} stochastic process 
$(\Y(t))_{t \geq 0}$ on $\R^N$ ?}
\end{align*}

The same argument with Schramm \cite{Sch00} will give that
$\Y(t)$ should be a continuous Markov process.
Moreover, Bauer, Bernard, and Kyt\"ol\"a \cite{BBK05}, 
Graham \cite{Gra07}, and Dub\'edat \cite{Dub07} argued that 
$(Y_i(t))_{t \geq 0}, i=1, \dots, N$ are
semi-martingales and the quadratic variations 
should be given by 
$d \langle Y_i, Y_j \rangle_t= \kappa \delta_{ij} dt, t \geq 0$,
$1 \leq i, j \leq N$ with $\kappa > 0$.
Then we will be able to assume that the system of SDEs for
$(\Y(t))_{t \geq 0}$ is in the form,
\begin{equation}
d Y_i(t)=\sqrt{\kappa} dB_i(t)+
F_i(\Y(t)) dt, \quad t \geq 0,
\quad i=1, \dots, N, 
\label{eqn:SDE_A1}
\end{equation}
where $(B_i(t))_{t \geq 0}, i=1, \dots, N$ are
independent one-dimensional standard Brownian motions,
$\kappa > 0$, and
$\{F_i(\x)\}_{i=1}^N$ are suitable functions
of $\x=(x_1, \dots, x_N)$ which do not explicitly
depend on $t$. 

In the following, We will give a \textit{theory} so that
the driving process 
$(\Y(t))_{t \geq 0}$ should be a time change 
of DYS$_{\beta}$ with $\beta=8/\kappa$
to construct a \textit{proper} multiple SLE.
Here we define the Gaussian free field (GFF)
and its generalization called the
imaginary surface with parameter $\chi$,
which are considered as the
distribution-valued random fields on $\HH$. 
Under the relation 
$\chi=2/\sqrt{\kappa}-\kappa/\sqrt{2}$, 
we regard the SLE/GFF coupling
studied by Dub\'edat, Sheffield, and Miller 
as a temporally stationary field, and extend 
it to multiple cases. We prove that
the multiple SLE/GFF coupling is established,  
if and only if the driving $N$-particle process
on $\R$ is identified with DYS$_{8/\kappa}$.

\begin{figure}[t]
\begin{center}
\includegraphics[scale=.50]{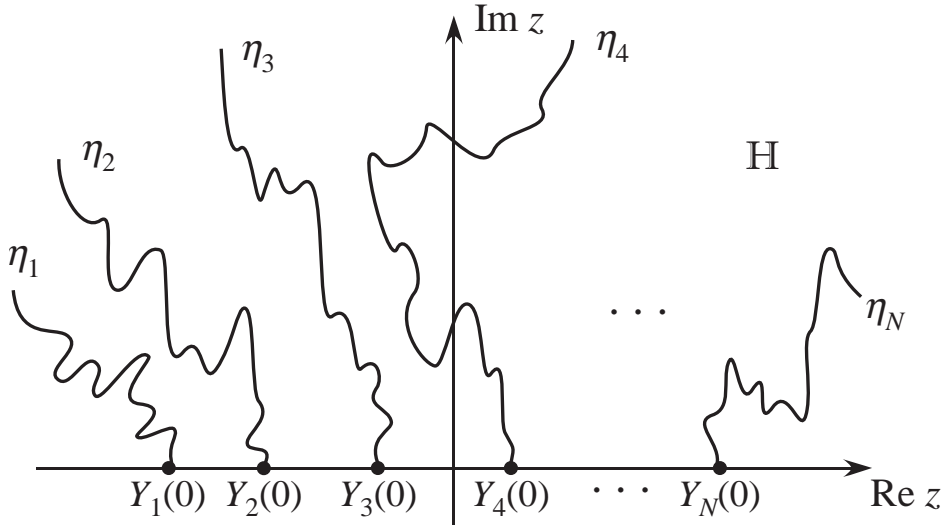}
\caption{Schematic picture of multiple SLE curves
$\eta_i, i=1, \dots, N$ in $\overline{\HH}$ driven by
$\Y(t)=(Y_1(t), \dots, Y_N(t)), t \geq 0$.
Here we assume that they are all simple and avoiding from
each other in $\HH$.
}
\end{center}
\label{fig:mSLE}       
\end{figure}

\subsection{Gaussian free field (GFF) 
with Dirichlet boundary condition } \label{sec:GFF}
\subsubsection{Bochner--Minlos Theorem} 
\label{sec:GFF_Dirichlet}

Here we start with the classical \textbf{Bochner theorem},
which states that
a probability measure on a finite dimensional Euclidean space
is determined by a \textbf{characteristic function}
which is a Fourier transform of
the probability measure. 
Note that we have identified the Fourier transform
of the transition probability density of BM
(\ref{eqn:charaBM0}) with the
characteristic function of BM (\ref{eqn:charaBM1})
in Section \ref{sec:1dimBM}, 
and the multitime Laplace transforms 
of the joint probability distribution function 
(\ref{eqn:GF1}) with the generating function 
of correlation functions for DYS$_2$ in 
Section \ref{sec:DMR_DP}.
First we define a functional of positive type.

\begin{df}
\label{thm:positive_type}
Let $\cV$ be a finite or infinite dimensional vector space.
A function $\psi: \cV \to \C$ is said to be 
a functional of positive type if for arbitrary $N \in \N$,
$\xi_1, \dots, \xi_N \in \cV$, and $z_1, \dots, z_N \in \C$, we have
\[
\sum_{n=1}^N \sum_{m=1}^N
\psi(\xi_n - \xi_m) z_n \overline{z_m} \geq 0.
\]
\end{df}
Then the following is proved.
\begin{lem}
\label{thm:positive_type2}
Let $\psi: \cV \to \C$ be a functional of positive type
on a vector space $\cV$. Then it follows that
{\rm (i)} 
$\psi(0) \geq 0$, 
{\rm (ii)} 
$\psi(\xi)=\overline{\psi(-\xi)}$ for all $\xi \in \cV$,
and
{\rm (iii)} 
$|\psi(\xi)| \leq \psi(0)$ for all $\xi \in \cV$.
\end{lem}
For $x, y \in \R^N$, the standard inner product is denoted by
$x \cdot y$ and we write $|x| := \sqrt{x \cdot x}$.
Let $\cB^N$ be the family of Borel sets in $\R^N$.
Then the following is known as the Bochner theorem.

\begin{thm}[Bochner theorem]
\label{thm:Bochner}
Let $\psi: \R^N \to \C$ be a continuous functional
of positive type such that $\psi(0)=1$. 
Then there exists a unique probability measure
$\rP$ on $(\R^N, \cB^N)$ such that
\[
\psi(\xi)=\int_{\R^N} e^{\sqrt{-1} x \cdot \xi} \, \rP(d x)
\quad \mbox{for \, $\xi \in \R^N$}.
\]
\end{thm}

If we consider the case that $\psi(\xi)$ 
is given by 
$\Psi(\xi):= e^{-|\xi|^2/2}, \xi \in \R^N$,
then 
the probability measure $\rP$ given by the Bochner theorem
is the 
\textbf{$N$-dimensional standard Gaussian measure},
\begin{align*}
\rP(d x) &=\frac{1}{(2\pi)^{N/2}} e^{-|x|^2/2} d x
=\prod_{i=1}^N \lambda_{\rN(0,1)}(dx_i),
\quad x=(x_1, \dots, x_N) \in \R^N.
\end{align*}
Hence we can say that
the finite-dimensional standard Gaussian measure 
$\rP$ is determined by
the characteristic function $\Psi(\xi)$ as
\begin{align*}
\Psi(\xi) &=\int_{\R^N} e^{\sqrt{-1} x \cdot \xi} \rP(d x)
= e^{-|\xi|^2/2} \quad \mbox{for \, $\xi \in \R^N$}.
\end{align*}

Now consider the case that $\cH$ is an infinite dimensional Hilbert space with inner product 
$\bra \cdot, \cdot \ket=\bra \cdot, \cdot \ket_{\cH}$
with $\| x \|=\| x \|_{\cH}=\sqrt{\bra x,x \ket_{\cH}}$,
$x \in \cH$.
The dual space of $\cH$ will be denoted by
$\cH^{\ast}$. 
Suppose that there were a probability measure $\rP$
on $\cH$ with a suitable $\sigma$-algebra such that
\[
\psi(\xi) 
= \int_{\cH} e^{\sqrt{-1} \bra x, \xi \ket} \rP(d x)
= e^{-\|\xi\|^2/2}
\quad \mbox{for $\xi \in \cH$}.
\]
Let $\{e_n\}_{n=1}^{\infty}$ be a complete orthonormal system (CONS)
of $\cH$.
If we set $\xi=t e_n, t \in \R$
for an arbitrary $n \in \N$, then
\[
\int_{\cH} e^{\sqrt{-1} t \bra x, e_n \ket} \rP(d x) = e^{-t^2/2},
\quad t \in \R.
\]
Since $x \in \cH$, we have $\bra x, e_n \ket \to 0$ as $n \to \infty$.
Therefore in the limit $n \to \infty$, 
the above equality gives $e^{-t^2/2}=1$, which is a contradiction.
This observation suggests that the application of the Bochner theorem
to an infinite dimensional space requires more consideration. 
The following arguments are base on \cite{Ara10}
and a note given by Koshida \cite{Koshida_note19}.

Let $D \subsetneq \C$ be a simply connected 
proper domain of $\C$
that is bounded.
We consider the case $\cH=L^2(D, \mu(dz))$ with 
$\bra f, g \ket := \int_D f(z) g (z) d \mu(z)$, $f, g \in L^2(D, \mu(dz))$,
where $\mu(dz)$ is the Lebesgue measure on $D \subset \C$;
$\mu(dz)=d \Re z d \Im z= \sqrt{-1}dz d\overline{z}/2$.
Let $\Delta$ be the 
\textbf{Dirichlet Laplacian} acting on $L^2(D, \mu(dz))$.
Then $-\Delta$ has positive discrete eigenvalues so that
\begin{equation}
-\Delta e_n = \lambda_n e_n, \quad
e_n \in L^2(D, \mu(dz)), \quad n \in \N.
\label{eqn:lambdaB}
\end{equation}
We assume that the eigenvalues are labeled 
in a non-decreasing order;
$0 < \lambda_1 \leq \lambda_2 \leq \cdots$.
The system of eigenvalue functions $\{e_n\}_{n \in \N}$ forms
a \textbf{complete orthonormal system} (CONS) of $L^2(D)$. 
The following is known as the
\textbf{Weyl formula}
\begin{lem}
\label{thm:Weyl}
Let $D \subsetneq \C$ be a simply connected finite domain.
The eigenvalues $\{\lambda_n\}_{n \in \N}$ 
of the operator $-\Delta$ on $D$
exhibit the
following asymptotic behavior, 
\[
\lim_{n \to \infty} \frac{\lambda_n}{n} =\rO(1).
\]
\end{lem}

For two functions $f, g \in \cC_{\rm c}^{\infty}(D)$,
their 
\textbf{Dirichlet inner product} is defined as
\begin{equation}
\bra f, g \ket_{\nabla} 
:= \frac{1}{2 \pi} \int_D (\nabla f)(z) \cdot (\nabla g)(z)
\mu(d z).
\label{eqn:Dirichlet_IP}
\end{equation}
The Hilbert space completion of $\cC_{\rm c}^{\infty}(D)$ with respect
to this Dirichlet inner product
will be denoted by $W(D)$.
We write 
$\|f\|_{\nabla}=\sqrt{\bra f, f \ket_{\nabla}}, f \in W(D)$.
If we set $u_n=\sqrt{2\pi/\lambda_n} \, e_n, n \in \N$,
then by integration by parts, we have
\[
\bra u_n, u_n \ket_{\nabla}=
\frac{1}{2 \pi} \bra u_n, (-\Delta) u_m \ket=\delta_{nm},
\quad n, m \in \N.
\]
Therefore $\{u_n\}_{n \in \N}$ forms a CONS of $W(D)$.

Let $\widehat{\cH}(D)$ be the space of formal real infinite
series in $\{u_n\}_{n \in \N}$.
This is obviously isomorphic to $\R^{\N}$
by setting 
$\widehat{\cH}(D) \ni \sum_{n \in \N} f_n u_n 
\mapsto (f_n)_{n \in \N} \in \R^{\N}$.
As a subspace of $\widehat{\cH}(D)$, $W(D)$ is isomorphic to
$\ell^2(\N) \subset \R^{\N}$.
For two formal series 
$f=\sum_{n \in \N} f_n u_n$,
$g=\sum_{n \in \N} g_n u_n \in \widehat{\cH}(D)$
such that
$\sum_{n \in \N} |f_n g_n| < \infty$,
we define their \textbf{pairing} as
\[
\bra f, g \ket_{\nabla}:=\sum_{n \in \N} f_n g_n.
\]
In the case when $f, g \in W(D)$, their pairing of course
coincides with the Dirichlet inner product 
(\ref{eqn:Dirichlet_IP}).

Notice that, for any $a \in \R$, the operator
$(-\Delta)^a$ acts on $\widehat{\cH}(D)$ as
\[
(-\Delta)^a \sum_{n \in \N} f_n u_n
:= \sum_{n \in \N} \lambda_n^a f_n u_n,
\quad (f_n)_{n \in \N} \in \R^{\N}.
\]
Using this fact, we define
$\cH_a(D) := (-\Delta)^a W(D)$, $a \in \R$,
each of which is a Hilbert space with inner product
\[
\langle f, g \rangle_a
:=\bra(-\Delta)^{-a}f, (-\Delta)^{-a}g \ket_{\nabla},
\quad f, g \in \cH_a(D).
\]
We write $\|\cdot\|_a:=\sqrt{\langle \cdot, \cdot \rangle_a}, a \in \R$.

\vskip 0.3cm
\begin{example}
\label{thm:a=1/2}
When $a=1/2$, we have
\[
\langle f, g \rangle_{1/2}
=\Big\bra (-\Delta)^{-1/2} f, (-\Delta)^{-1/2}g \Big\ket_{\nabla}
=\frac{1}{2 \pi} \bra f, g \ket, \quad
f, g \in \cH_{1/2}(D).
\]
Therefore $\cH_{1/2}(D) = L^2(D, \mu(dz))$.
\end{example}
We can prove the following two lemmas.
\begin{lem}
\label{thm:Ha<Hb}
Assume $a< b$. Then
$\cH_a(D) \subset \cH_b(D)$.
\end{lem}
\begin{lem}
\label{thm:dual}
Let $a \in \R$ and fix $h \in \cH_a(D)$.
Then the assignment
\[
\bra h, \cdot \ket_{\nabla} : \cH_{-a}(D) \to \R
\quad \mbox{such that} \quad
\cH_{-a}(D) \ni f \mapsto \bra h, f \ket_{\nabla} \in \R
\]
is well-defined and continuous. 
In particular, $\cH_a(D)$ and $\cH_{-a}(D)$
makes a dual pair of Hilbert spaces with respect to 
the Dirichlet inner product $\bra \cdot, \cdot \ket_{\nabla}$. 
\end{lem}

\vskip 0.3cm
\begin{rem}
\label{thm:Remark3_1}
Since $\cH_{1/2}(D)=L^2(D, \nu(dz))$ as mentioned
in Example \ref{thm:a=1/2}, 
the members of $\cH_a(D)$
with $a > 1/2$ cannot be functions, but are 
\textbf{distributions}.
\end{rem}
\vskip 0.3cm

Define 
\begin{equation}
\cE(D) := \bigcup_{a > 1/2} \cH_a(D).
\label{eqn:cE}
\end{equation}
Then its dual Hilbert space is identified with
$\cE(D)^{\ast} :=\bigcap_{a < -1/2} \cH_a(D)$ 
by Lemma \ref{thm:dual}, and 
\[
\cE(D)^{\ast} \subset W(D) \subset \cE(D)
\]
is established (by definition and Lemma \ref{thm:Ha<Hb}).
Here $(\cE(D)^{\ast}, W(D), \cE(D))$ is called a
\textbf{Gel'fand triple}.
We set 
$\Sigma_{\cE(D)} 
=\sigma(\{\bra \cdot, f \ket_{\nabla} : f \in \cE(D)^{\ast} \})$.
On such a setting, the following is
obtained.
This theorem is the extension of the Bochner theorem 
(Theorem \ref{thm:Bochner})
and is called 
the \textbf{Bochner--Minlos theorem}
(see, for instance, \cite{Hid80,She07,Ara10}). 

\begin{thm}[Bochner--Minlos theorem]
\label{thm:BM}
Let $\psi$ be a continuous function of positive type
on $W(D)$ such that $\psi(0)=1$.
Then there exists a unique probability measure 
$\bP$ on $(\cE(D), \Sigma_{\cE(D)})$
such that
\begin{equation}
\psi(f)=\int_{\cE(D)} e^{\sqrt{-1} \bra h, f \ket_{\nabla}} \bP(d h)
\quad \mbox{for $f \in \cE(D)^{\ast}$}.
\label{eqn:BMeq1}
\end{equation}
\end{thm}

Under certain conditions on $\psi$,
the domain of function $f$ for (\ref{eqn:BMeq1})
can be extended from $\cE(D)^{\ast}$ to $W(D)$.
It is easy to verify that the functional
\begin{equation}
\Psi(f):= e^{-\|f\|_{\nabla}^2/2}
\label{eqn:PsiA}
\end{equation}
satisfies the conditions.
Then the following holds
with a probability measure $\bP$ on $(\cE(D), \Sigma_{\cE(D)})$,
\begin{align}
\Psi(f) &=\int_{\cE(D)} e^{\sqrt{-1} \bra h, f \ket_{\nabla}} \bP(d h)
=e^{-\|f\|_{\nabla}^2/2}
\quad \mbox{for $f \in W(D)$}.
\label{eqn:GFF}
\end{align}

\begin{df}[Dirichlet boundary GFF]
\label{thm:GFF_Dirichlet}
A \textbf{Gaussian free field} (GFF) \textbf{with 
Dirichlet boundary condition} is defined 
as a pair $((\Omega^{\GFF}, \cF^{\GFF}, \P^{\GFF}),H)$ 
of a probability space $(\Omega^{\GFF}, \cF^{\GFF}, \P^{\GFF})$ 
and an isometry
\[
H : W(D) \to L^{2}(\Omega^{\GFF}, \cF^{\GFF}, \P^{\GFF})
\]
such that each $H(f)$, $f\in W(D)$ is a Gaussian random variable.
\end{df}

For each $f \in W(D)$, (\ref{eqn:GFF}) gives
a Gaussian random variable
$H(f):=\bra H, f \ket_{\nabla}$ 
for the variable $h \in \cE(D)$ by 
\[
\cE(D) \ni h \mapsto \bra h, f \ket_{\nabla} 
\in L^{2}(\cE(D), \Sigma_{\cE(D)}, \bP).
\]
In this way, (\ref{eqn:GFF}) ensures that the pair of 
of the probability space
\[
(\Omega^{\GFF}, \cF^{\GFF}, \P^{\GFF})
:= (\cE(D), \Sigma_{\cE(D)}, \bP)
\]
with the isometry $H$, 
$((\cE(D), \Sigma_{\cE(D)}, \bP), H)$ gives a GFF with Dirichlet boundary condition. 
We often just call $H$ a 
\textbf{Dirichlet boundary GFF} without referring 
to the probability space $ (\cE(D), \Sigma_{\cE(D)}, \bP)$.
By this definition, (\ref{eqn:GFF}) is written as
\begin{equation}
\E^{\GFF}[ e^{\sqrt{-1} \bra H, f \ket_{\nabla}} ]
=e^{-\|f\|_{\nabla}^2/2}, \quad f \in W(D).
\label{eqn:exp_av}
\end{equation}
This determines all the moments of 
the family of Gaussian random variables 
$\{\bra H, f \ket_{\nabla} : f \in W(D)\}$.
For examples, the \textbf{covariance} is given by
\begin{equation}
\Cov[\bra H, f \ket_{\nabla},
\bra H, g \ket_{\nabla} ]
:= \E^{\GFF} \Big[
\bra H, f \ket_{\nabla} \bra H, g \ket_{\nabla} \Big]
=\bra f, g \ket_{\nabla}, \quad
\mbox{$f, g \in W(D)$}.
\label{eqn:cov1}
\end{equation}
In particular, the variance is written as
\[
\Var[\bra H, f \ket_{\nabla} ]
:=\E^{\GFF} \Big[\bra H, f \ket_{\nabla}^2 \Big]
=\|f\|_{\nabla}^2, \quad f \in W(D).
\]
By introducing a parameter 
$\theta \in \R$, and by replacing $f \to \theta f$
in (\ref{eqn:exp_av}), the above result implies the following
useful formula.

\begin{lem}
\label{thm:chara1}
For the family of \textbf{centered (mean zero) 
Gaussian random variables}  
$\{\bra H, f \ket_{\nabla} : f \in W(D)\}$
given by the Dirichlet boundary GFF, 
\begin{equation}
\E^{\GFF}[ e^{\sqrt{-1} \theta \bra H, f \ket_{\nabla}} ]
=\exp \left(- \frac{\theta^2}{2} \Var\Big[ \bra H, f \ket_{\nabla} \Big]
\right), 
\quad \theta \in \R,
\quad f \in W(D).
\label{eqn:GFF_chara}
\end{equation}
\end{lem}
We call (\ref{eqn:GFF_chara}) 
the \textbf{characteristic function of
Dirichlet boundary GFF}. 

\subsubsection{Conformal invariance of GFF} 
\label{sec:conformal_invariance}

Assume that $D, D' \subsetneq \C$ are simply connected domains
and let $\varphi: D' \to D$ be a conformal map.

\begin{lem}
\label{thm:conformal_inv}
The Dirichlet inner product (\ref{eqn:Dirichlet_IP}) is 
conformally invariant, 
that is,
\[
\int_{D} (\nabla f)(z) \cdot (\nabla g)(z) \mu(d z)
=\int_{D'} (\nabla (f \circ \varphi))(z) \cdot 
(\nabla (g \circ \varphi))(z) \mu(d z)
\quad \mbox{for $f, g \in \cC_{\rm c}^{\infty}(D)$}.
\]
\end{lem}

From the above lemma, we see that the
\textbf{pull-back} 
$\varphi^{\ast} : W(D) \ni f \mapsto f \circ \varphi \in W(D')$
is an isomorphism. This allows one to consider
a GFF on an unbounded domain.
Namely, if $D'$ is bounded on which a Dirichlet GFF $H$ is defined, 
but $D$ is unbounded, 
we can define a family 
$\{\bra \varphi_{\ast} H, f \ket_{\nabla} : f \in W(D)\}$ by
\[
\bra \varphi_{\ast} H, f\ket_{\nabla} 
:= \bra H, \varphi^{\ast} f \ket_{\nabla}, \quad f \in W(D).
\]
The \textbf{transformation rule of covariance} is given by
\begin{align}
\Cov[\bra \varphi_{\ast} H, f \ket_{\nabla}, 
\bra \varphi_{\ast} H, g \ket_{\nabla}]
&=
\E^{\GFF} \Big[
\bra \varphi_{\ast} H, f \ket_{\nabla} 
\bra \varphi_{\ast} H, g \ket_{\nabla} \Big]
=\bra \varphi^{\ast} f, \varphi^{\ast} g \ket_{\nabla}
=\bra f, g \ket_{\nabla}
\nonumber\\
&
=\E^{\GFF}\Big[
\bra H, f \ket_{\nabla} 
\bra H, g \ket_{\nabla} \Big]
=\Cov[\bra H, f \ket_{\nabla}, 
\bra H, g \ket_{\nabla} ].
\label{eqn:cov_str}
\end{align}
Relying on the following formal computation
\begin{align*}
\bra \varphi_{\ast} H, f \ket_{\nabla}
=\bra H, \varphi^{\ast} f \ket_{\nabla}
&=\frac{1}{2 \pi} \int_{D'} (\nabla H)(z) 
\cdot (\nabla f \circ \varphi)(z) \mu(d z)
\nonumber\\
&=\frac{1}{2 \pi} \int_D (\nabla H \circ \varphi^{-1})(z) 
\cdot (\nabla f)(z) \mu(d z)
\end{align*}
we understand the equality $\varphi_{\ast} H = H \circ \varphi^{-1}$.
By the fact (\ref{eqn:cov_str}) such that
the covariance structure does not change under 
a conformal map $\phi$, we say that
\textbf{the GFF is conformally invariant}. 

\subsubsection{The Green's function of GFF} 
\label{sec:Green}

By a formal integration by parts, we see that
\begin{align*}
\bra H, f\ket_{\nabla} 
&= \frac{1}{2\pi} \int_{D} (\nabla H)(z) \cdot (\nabla f)(z) \mu(d z)
=\frac{1}{2\pi} \int_{D} H(z) (-\Delta f)(z) \mu(d z)
\nonumber\\
&=\frac{1}{2 \pi} \bra H, (-\Delta) f \ket.
\end{align*}
Motivated by this observation, we define
\begin{equation}
\bra H, f \ket:=2 \pi \bra H, (-\Delta)^{-1} f \ket_{\nabla} 
\quad \mbox{for $f \in \sD((-\Delta)^{-1})$}, 
\label{eqn:inner_product}
\end{equation}
where
$\sD((-\Delta)^{-1})$ denotes the domain of
$(-\Delta)^{-1}$ in $W(D)$.
Note that if $D$ is bounded, then $(-\Delta)^{-1}$ is
a bounded operator, but if $D$ is unbounded,
then $(-\Delta)^{-1}$ is not defined on $W(D)$.
The action of $(-\Delta)^{-1}$ is expressed as
an integral operator and the integral kernel
is known as \textbf{the Green's function}. Namely,
\[
((-\Delta)^{-1} f)(z)
=\frac{1}{2 \pi} \int_{D} G_D(z,w) f(w) \mu(d w),
\quad \mbox{a.e.} \, z \in D, \quad
f \in \sD((-\Delta)^{-1}),
\]
where $G_D(z,w)$ denotes the Green's function of $D$
under the Dirichlet boundary condition.
Hence the covariance of $\bra H, f \ket$ and $\bra H, g \ket$ with
$f, g \in \sD((-\Delta)^{-1})$ is written as
\begin{equation}
\E^{\GFF}[ \bra H, f \ket \bra H, g \ket] 
=\int_{D \times D} f(z) G_D(z, w) g(w) \mu(d z) \mu(d w).
\label{eqn:Green}
\end{equation}
When we symbolically write
\[
\bra H, f \ket=\int_{D} H(z) f(z) \mu(d z),
\quad f \in \sD((-\Delta)^{-1}),
\]
the covariance structure can be expressed as
\[
\E^{\GFF}[H(z) H(w)]= G_D(z, w), \quad
z, w \in D, \quad n \not=w.
\]
From the formula (\ref{eqn:Green}), we see that
$\cC_{\rm c}^{\infty}(D) \subset \sD((-\Delta)^{-1})$.
In the following, we will consider the
family of random variables 
$\{\bra H, f \ket : f \in \cC_{\rm c}^{\infty}(D) \}$
to characterize the GFF, $H$.

\vskip 0.3cm
\begin{example}
\label{thm:G_H}
When $D$ is the upper half plane $\HH$,
\begin{align*}
G_{\HH}(z,w) &= \log \left|
\frac{z-\overline{w}}{z-w} \right|
=\log|z-\overline{w}| -\log|z-w|
\nonumber\\
&= \Re \log(z-\overline{w}) - \Re \log(z-w), 
\end{align*}
$z, w \in \HH, z \not= w$.
\end{example}
\vskip 0.3cm
\begin{example}
\label{thm:G_O}
When $D$ is the first orthant 
$\O:=\{z : \Re z >0, \Im z >0\}$,
\begin{align*}
G_{\O}(z,w) &= \log \left|
\frac{(z-\overline{w})(z+\overline{w})}{(z-w)(z+w)} \right|
\nonumber\\
&=\log|z-\overline{w}| + \log|z+\overline{w}|
-\log|z-w| - \log|z+w|, 
\nonumber\\
&=\Re \log (z-\overline{w}) + \Re \log(z+\overline{w})
- \Re \log (z-w) - \Re \log(z+w), 
\end{align*}
$z, w \in \O, z \not= w$.
\end{example}
\vskip 0.3cm

The conformal invariance of GFF implies that
for a conformal map
$\varphi: D' \to D$, we have the equality,
\begin{equation}
G_{D'}(z, w)=G_D(\varphi(z), \varphi(w)), \quad z, w \in D'.
\label{eqn:G_conformal}
\end{equation}
In the following, we will regard the upper half plane $\HH$
as the representative of the simply connected
domain $D \subsetneq \C$.
Since each $D' \subsetneq \C$ is specified by 
the conformal transformation 
$\varphi: D' \to \HH$,
we put this in the superscript and write
\begin{equation}
G^{\varphi}(z, w) := G_{\varphi^{-1}(\HH)}(z, w)
=G_{\HH}(\varphi(z), \varphi(w)).
\label{eqn:G_notation1a}
\end{equation}
By Example \ref{thm:G_H}, 
\begin{align}
G^{\varphi}(z, w)
&= \log|\varphi(z)-\overline{\varphi(w)}| 
-\log|\varphi(z)-\varphi(w)|
\nonumber\\
&=\Re \log(\varphi(z)-\overline{\varphi(w)}) 
- \Re \log(\varphi(z)-\varphi(w)).
\label{eqn:G_notation1b}
\end{align}

\subsection{GFF coupled with stochastic log-gase} 
\label{sec:GFF_BPs}
\subsubsection{GFF transformed by multiple SLE}
\label{sec:GFF_SLE}

We write the probability space of the multiple SLE as
$(\Omega^{\SLE}, \cF^{\SLE}, (\cF^{\SLE})_{\geq 0}, \P^{\SLE})$.
Now we define the direct product of this space and
that for the Dirichlet boundary GFF, 
\begin{equation}
(\Omega, \cF, \P)=(\Omega^{\SLE} \times \Omega^{\GFF},
\cF^{\SLE} \times \cF^{\GFF}, 
\P^{\SLE} \times \P^{\GFF}).
\label{eqn:direct}
\end{equation}
We consider the multiple SLE as well as GFF in this
extended probability space, and 
the multiple SLE is assumed to be adapted to
the filtration, 
\[
\cF_{t}=\cF^{\SLE}_{t} \times \{\emptyset, \Omega^{\GFF} \},
\quad t \geq 0.
\]

We assume that the multiple SLE (\ref{eqn:mSLE1}) 
driven by $\Y(t)=(Y_1(t), \dots, Y_N(t))$ 
under the initial condition $g_0(z)=z \in \HH$ 
has a unique solution $(g_t)_{t \geq 0}$.
Moreover, we assume that this solution can be extended to
$\R$ and determines 
\[
\eta_i := \{\eta_i(t) : t \geq 0\}, \quad
i=1, \dots, N
\]
by the equations
\begin{equation}
Y_i(t)=g_t(\eta_i(t)), \quad i=1, \dots, N.
\label{eqn:m_slit1}
\end{equation}
We consider the situation such that
each $\eta_i, i=1, \dots, N$ make a continuous curve in 
$\overline{\HH}$
\footnote{
This is proved by \cite{KK21b}. 
See Theorem \ref{thm:3phases_multi} shown
in Section \ref{sec:solution}.
}, 
and for each $t \geq 0$, we define
\begin{equation}
\mbox{
$\HH^{\eta}_t :=$ the unbounded component
of $\displaystyle{\HH \backslash \bigcup_{i=1}^N \eta_i(0, t], 
\quad t \geq 0}$}. 
\label{eqn:HetaB2}
\end{equation}
Then the solution of (\ref{eqn:mSLE1}) 
gives the conformal transformation
from $\HH^{\eta}_t$ to $\HH$ and hence it is denoted by
$g_{\HH^{\eta}_t}$.

In summary, we assume the following;
\begin{align}
\frac{d g_{\HH^{\eta}_t}(z)}{dt}
&=\sum_{i=1}^N \frac{2}{g_{\HH^{\eta}_t}(z)-Y_i(t)}, \quad t \geq 0,
\quad g_{\HH^{\eta}_0}(z)=g_{\HH}(z)=z \in \HH,
\nonumber\\
Y_i(t) &=
g_{\HH^{\eta}_t}(\eta_i(t)) :=
\lim_{z \to \eta_i(t), z \in \HH^{\eta}_t} g_{\HH^{\eta}_t}(z),
\quad i=1, \dots, N, \quad t \geq 0, 
\nonumber\\
g_{\HH^{\eta}_t} &: \mbox{
conformal : $\HH^{\eta}_t \to \HH$}.
\label{eqn:mSLE2}
\end{align}

By the notation 
(\ref{eqn:G_notation1a}), we put
\begin{equation}
G^{g_{\HH^{\eta}_t}}(z, w)
:=G_{\HH}(g_{\HH^{\eta}_t}(z), g_{\HH^{\eta}_t}(w)),
\quad t \geq 0
\label{eqn:GreenB}
\end{equation}
and define the Dirichlet boundary GFF,
${g_{\HH^{\eta}_t}}_* H$ on $\HH^{\eta}_t, t \geq 0$, 
so that its Green's function is given by (\ref{eqn:GreenB}).
We can prove the following.
\begin{prop}
\label{thm:dG}
For $t \geq 0$, 
\begin{align}
\frac{dG^{g_{\HH^{\eta}_t}}(z, w)}{dt}
&= -4 \sum_{i=1}^N \Im \left(\frac{1}{g_{\HH^{\eta}_t}(z)-Y_i(t)} \right)
\Im \left( \frac{1}{g_{\HH^{\eta}_t}(w)-Y_i(t)} \right), 
\quad z, w \in \HH^{\eta}_t.
\label{eqn:dG}
\end{align}
\end{prop}
\noindent{\it Proof.} 
By (\ref{eqn:G_notation1b}), 
\begin{align*}
\frac{dG^{g_{\HH^{\eta}_t}}(z, w)}{dt}
&= \frac{d}{dt}
\Big[ \Re \log(g_{\HH^{\eta}_t}(z)-\overline{g_{\HH^{\eta}_t}(w)})
- \Re \log(g_{\HH^{\eta}_t}(z)-g_{\HH^{\eta}_t}(w)) \Big]
\nonumber\\
&= \Re \frac{1}
{g_{\HH^{\eta}_t}(z)-\overline{g_{\HH^{\eta}_t}(w)}}
\left(\frac{d g_{\HH^{\eta}_t}(z)}{dt}-
\frac{d \overline{g_{\HH^{\eta}_t}(w)}}{dt} \right)
\nonumber\\
&\quad - \Re \frac{1}
{g_{\HH^{\eta}_t}(z)-g_{\HH^{\eta}_t}(w)}
\left(\frac{d g_{\HH^{\eta}_t}(z)}{dt}-
\frac{d g_{\HH^{\eta}_t}(w)}{dt} \right), 
\quad t \geq 0.
\end{align*}
Since the multiple SLE (\ref{eqn:mSLE2}) gives
\begin{align*}
\frac{d g_{\HH^{\eta}_t}(z)}{dt}-
\frac{d \overline{g_{\HH^{\eta}_t}(w)}}{dt}
& = \sum_{i=1}^N \frac{2}{g_{\HH^{\eta}_t}(z)-Y_i(t) }
- \sum_{i=1}^N \frac{2}{\overline{g_{\HH^{\eta}_t}(w)}-Y_i(t)}
\nonumber\\
&= - 2(g_{\HH^{\eta}_t}(z)-\overline{g_{\HH^{\eta}_t}(w)})
\sum_{i=1}^N \frac{1}{(g_{\HH^{\eta}_t}(z)-Y_i(t))
(\overline{g_{\HH^{\eta}_t}(w)}-Y_i(t) )}, 
\end{align*}
the above is written as 
\begin{align*}
\frac{dG^{g_{\HH^{\eta}_t}}(z, w)}{dt}
&=-2 \sum_{i=1}^N \Re \frac{1}{(g_{\HH^{\eta}_t}(z)-Y_i(t))
(\overline{g_{\HH^{\eta}_t}(w)}-Y_i(t) )}
\nonumber\\
& \quad
+2 \sum_{i=1}^N \Re \frac{1}{(g_{\HH^{\eta}_t}(z)-Y_i(t))
(g_{\HH^{\eta}_t}(w)-Y_i(t) )}, \quad t \geq 0.
\end{align*}
In general, the equality
\[
\Re \zeta \overline{w} - \Re \zeta w =
2 \Im \zeta \, \Im w
\]
holds for $\zeta, w \in \C$, and hence 
(\ref{eqn:dG}) is proved.
\qed

For a given domain $A \subset \HH$, we assume that
$\supp(f)$ of
$f \in \cC^{\infty}_{\rm c}(\HH)$ satisfies 
$\supp(f) \subset A$. 
Then we define
\begin{align}
E_A^{g_{\HH^{\eta}_t}}(f) &:= \int_{A \times A} f(z) G^{g_{\HH^{\eta}_t}}(z, w)
f(w) \mu(dz) \mu (dw).
\label{eqn:energy}
\end{align}
By Proposition \ref{thm:dG}, we have
\begin{align}
\frac{d E_A^{g_{\HH^{\eta}_t}}(f)}{dt}
&= - \sum_{i=1}^N \left( \int_A 
\Im \frac{2}{g_{\HH^{\eta}_t}(z) - Y_i(t)} f(z) \mu(dz) \right)^2,
\label{eqn:energy2}
\end{align}
It implies that $E^{g_{\HH^{\eta}_t}}_A(f)$ 
is non-increasing function of $t$. 

\subsubsection{Imaginary surface and extended GFF}
\label{sec:IS}

Consider a simply connected domain $D \subsetneq \C$
and write $\cC_{\rm c}^{\infty}(D)$ for the space of 
real smooth functions
on $D$ with compact support.
Assume $h \in \cC_{\rm c}^{\infty}(D)$ and 
consider a smooth vector field
$e^{\sqrt{-1} (h/\chi + \theta)}$ with parameters $\chi, \theta \in \R$.
Then a \textbf{flow line} along this vector field, 
$\eta : (0, \infty) \ni t \mapsto \eta(t) \in D$,
starting from 
$\lim_{t \to 0} \eta(t) =: \eta(0)=x \in \partial D$ 
is defined (if exists) as the 
solution of the ordinary differential equation (ODE)
\cite{She16,MS16a}
\begin{equation}
\frac{d \eta(t)}{dt} = e^{\sqrt{-1}\{h(\eta(t))/\chi + \theta\}},
\quad t \geq 0, \quad \eta(0) = x.
\label{eqn:flow1}
\end{equation}
Let $\widetilde{D} \subsetneq \C$ be another 
simply connected domain and
consider a conformal map $\varphi: \widetilde{D} \to D$.
Then we define the pull-back of the flow line $\eta$ by $\varphi$ as
$\widetilde{\eta}(t)=(\varphi^{-1} \circ \eta)(t)$.
That is, 
$\varphi(\widetilde{\eta}(t))=\eta(t)$, 
and the derivatives with respect to $t$ of the both sides of this equation 
gives 
\[
\varphi'(\widetilde{\eta}(t)) 
\frac{d \widetilde{\eta}(t)}{dt} = \frac{d \eta(t)}{dt} 
\]
with the notation
\[
\varphi'(z) := \frac{d \varphi(z)}{dz}.
\]
We use the polar coordinate
$\varphi'(\cdot)=|\varphi'(\cdot)| e^{\sqrt{-1} \arg \varphi'(\cdot)}$,
where $\arg \zeta$ of $\zeta \in \C$ is a priori 
defined up to additive multiples of $2 \pi$, and hence
we have 
$d \widetilde{\eta}(t)/dt=
e^{\sqrt{-1}\{ (h \circ \varphi 
- \chi \arg \varphi')(\widetilde{\eta}(t))/\chi + \theta\}}
/|\varphi'(\widetilde{\eta}(t))|, t \geq 0$.
If we perform a time change $t \to \tau=\tau(t)$ by putting
$t=\int_0^{\tau} ds/|\varphi'(\widetilde{\eta}(s))|$ and 
$\widehat{\eta}(t):=\widetilde{\eta}(\tau(t))$, then 
the above equation becomes
\[
\frac{d \widehat{\eta}(t)}{dt}
= e^{\sqrt{-1}\{ (h \circ \varphi 
- \chi \arg \varphi')(\widehat{\eta}(t))/\chi + \theta\}},
\quad t \geq 0.
\]
Since a time change preserves the image of a flow line,
we can identify $h$ on $D$ and 
$h \circ \varphi - \chi \arg \varphi'$ 
on $\widetilde{D}=\varphi^{-1}(D)$.
In \cite{She16,MS16a,MS16b,MS16c,MS17}, such a flow line
is considered also in the case that $h$ is given by 
an instance of a GFF defined as follows.
\begin{df}
\label{thm:def_GFF}
Let $D\subsetneq \mathbb{C}$ be a simply connected 
domain and $H$ be 
a Dirichlet boundary GFF on $D$. 
A GFF on $D$ is a random distribution $h$ 
of the form 
\[
h=H+u, 
\]
where $u$ is a 
deterministic harmonic function on $D$.
\end{df} 
\noindent
Since a GFF is not function-valued, 
but it is a \textbf{distribution-valued random field}
(see Remark 4.1 in Section \ref{sec:GFF}),
the ODE in the form (\ref{eqn:flow1}) 
no longer makes sense mathematically 
in the classical sense. 
Using the theory of SLE, however, 
the notion of flow lines was generalized as follows
\cite{SS09,SS13}.

Consider the collection
\[
\sS:=\left\{(D,h) \Bigg| 
\substack{\, \ D\subsetneq \mathbb{C}:\ 
\mbox{\small simply connected} \\ 
h:\ \mbox{\small GFF on $D$}} \right\}.
\]
Fixing a parameter $\chi \in \R$, 
we define the following equivalence relation in $\sS$.
\begin{df}
\label{thm:def_IS}
Two pairs $(D,h)$ and
$(\widetilde{D}, \widetilde{h}) \in \sS$ are 
\textbf{$\chi$-equivalent} if there exists a conformal map 
$\varphi: \widetilde{D} \to D$ and
\[
\widetilde{h} 
\law= \varphi_* h - \chi \arg \varphi^{\prime}. 
\]
In this case, we write 
$(D, h) \sim (\widetilde{D}, \widetilde{h})$.
\end{df}
\noindent
We call each element belonging to 
the equivalence class $\sS/\sim$ 
an \textbf{imaginary surface} \cite{MS16a} 
(or an \textbf{AC surface} \cite{She16};  
(AC means a combination of an altimeter and a compass.) 
That is, in this equivalence class, 
a conformal map $\varphi$ causes not only 
a coordinate change of a GFF as $h \mapsto \varphi_* h$
associated with changing the domain of definition of the field as
$D \mapsto \varphi^{-1}(D)$,
but also an addition of a deterministic 
harmonic function $-\chi \arg \varphi'$ to the field.
Notice that this definition depends on
one parameter $\chi \in \R$.

\subsubsection{Complex-valued logarithmic potentials and
local martingales}
\label{sec:complex_log}

Assume that the driving process $(\Y(t))_{t \geq 0}$
of the multiple SLE is given by (\ref{eqn:SDE_A1}),
that is 
\begin{equation}
d Y_i(t)=\sqrt{\kappa} dB_i(t)+
F_i(\Y(t)) dt, \quad t \geq 0,
\quad i=1, \dots, N, 
\label{eqn:SDE_A1b}
\end{equation}

Now for $z \in \C, y_i \in \R, i=1, \dots, N$,
we consider a \textit{sum} of
\textbf{complex-valued logarithmic potentials}
\cite{KK20,KK21a,KK21b,Kos21a}, 
\begin{equation}
\Phi(z, \y):= \sum_{i=1}^N \log(z-y_i),
\label{eqn:C_log}
\end{equation}
and consider a stochastic process
$(\Phi(g_{\HH^{\eta}_t}(z), \Y(t)))_{t \geq 0}$.
It should be noted that introduction of 
a complex-valued logarithmic potential 
in the form
$\log(z-y)$
has been reported in the previous papers
\cite{Dub09,SS13,KM13} in order to
establish the coupling between the original SLE$_{\kappa}$ 
and the extended GFFs.
In this sense, the present generalization for
the multiple SLE/GFF coupling seems to be very straightforward.
But this generalization is nontrivial.
Moreover, it was shown by \cite{Kos21a} that
this choice is the only possibility to obtain the
multiple SLE/GFF coupling when $\kappa \not=4$. 
See also \cite{KM13} for the consideration from the
viewpoint of conformal field theory.

Since 
\begin{align*}
& \frac{\partial}{\partial z} \log(z-x)=\frac{1}{z-x},
\quad
\frac{\partial}{\partial x} \log(z-x)=- \frac{1}{z-x},
\quad
\frac{\partial^2}{\partial x^2} \log(z-x)= -\frac{1}{(z-x)^2},
\end{align*}
It\^o's formula (\ref{eqn:Ito1}) with 
the multiple SLE$_{\kappa}$ (\ref{eqn:mSLE1}) gives
\begin{align*}
d \log(g_{\HH^{\eta}_t}(z)-Y_i(t))
&= \frac{1}{g_{\HH^{\eta}_t}(z)-Y_i(t)}
\Big[ d g_{\HH^{\eta}_t}(z)-d Y_i(t) \Big]
\nonumber\\
&= \frac{1}{g_{\HH^{\eta}_t}(z)-Y_i(t)}
\left[
2 \sum_{j=1}^N \frac{1}{g_{\HH^{\eta}_t}(z)-Y_i(t)}dt
-\sqrt{\kappa} d B_i(t) - F_i(\Y(t)) dt \right]
\nonumber\\
& \, 
- \frac{1}{2} \frac{1}{(g_{\HH^{\eta}_t}(z)-Y_i(t))^2} \kappa dt, 
\end{align*}
and hence we have
\begin{align}
d \Phi(g_{\HH^{\eta}_t}(z), \Y(t))
&= - \sqrt{\kappa} \sum_{i=1}^N \frac{1}{g_{\HH^{\eta}_t}(z)-Y_i(t)} 
d B_i(t)
+ 2  
\left(\sum_{i=1}^N \frac{1}{g_{\HH^{\eta}_t}(z)-Y_i(t)} \right)^2 dt
\nonumber\\
& \quad
- \sum_{i=1}^N \frac{1}{g_{\HH^{\eta}_t}(z)-Y_i(t)} F_i(\bxi) dt
- \frac{\kappa}{2} \sum_{i=1}^N 
\frac{1}{(\varphi_t(z)-Y_i(t))^2} dt.
\label{eqn:dPhi1}
\end{align}
The factor appearing in the second term of RHS
is written as 
\begin{align}
\left(\sum_{i=1}^N \frac{1}{g_{\HH^{\eta}_t}(z)-Y_i(t)} \right)^2 
&= \sum_{i=1}^N \sum_{1 \leq j \leq N, j \not=i}
\frac{1}{(g_{\HH^{\eta}_t}(z)-Y_i(t))(g_{\HH^{\eta}_t}(z)-Y_j(t))} 
+ \sum_{i=1}^N \frac{1}{(g_{\HH^{\eta}_t}(z)-Y_j(t))^2} 
\nonumber\\
&= 2
\sum_{i=1}^N 
\frac{1}{g_{\HH^{\eta}_t}(z)-Y_i(t)}
\sum_{\substack{1 \leq j \leq N,  \cr j \not=i}}
\frac{1}{Y_i(t)-Y_j(t)} 
+ \sum_{i=1}^N \frac{1}{(g_{\HH^{\eta}_t}(z)-Y_j(t))^2}.
\label{eqn:key}
\end{align}
By differentiate the multiple SLE (\ref{eqn:mSLE1})
with respect to $z$, we have
\[
d g_{\HH^{\eta}_t}'(z)=-2 g_{\HH^{\eta}_t}'(z) \sum_{i=1}^N 
\frac{1}{(g_{\HH^{\eta}_t}(z)-Y_i(t))^2} dt,
\]
and hence the equality
\[
d \log g_{\HH^{\eta}_t}'(z)
=-2 \sum_{i=1}^N \frac{1}{(g_{\HH^{\eta}_t}(z)-Y_i(t))^2} dt
\]
is established. 
Therefore, if (\ref{eqn:dPhi1}) is multiplied by $2/\sqrt{\kappa}$,
then we obtain 
\begin{align}
\frac{2}{\sqrt{\kappa}}
d \Phi(g_{\HH^{\eta}_t}(z), \Y(t))
&= - 2 \sum_{i=1}^N \frac{1}{g_{\HH^{\eta}_t}(z)-Y_i(t)} dB_i(t)
\nonumber\\
& \quad 
+\frac{2}{\sqrt{\kappa}} 
\sum_{i=1}^N \frac{1}{g_{\HH^{\eta}_t}(z)-Y_i(t)}
\left[ 4 \sum_{\substack{1 \leq j \leq N,  \cr j \not=i}}
\frac{1}{Y_i(t)-Y_j(t)} - F_i(\Y(t)) 
\right] dt
\nonumber\\
& \quad
- \frac{1}{\sqrt{\kappa}}
\left( 2- \frac{\kappa}{2} \right)
d \log g_{\HH^{\eta}_t}'(z).
\label{eqn:C_log2}
\end{align}

Comparing the above calculation result with
the $\chi$-equivalence 
\[
H_{D_1} \law=
\varphi_* H_{D_2} 
- \chi \arg \varphi' 
\]
defined by Definition \ref{thm:def_IS}, we find the following:
\begin{align}
\mbox{\rm (i)} \quad & 
\mbox{if} \quad
\chi=\frac{2}{\sqrt{\kappa}}- \frac{\sqrt{\kappa}}{2},
\nonumber\\
& \qquad 
\, \mbox{then 
the imaginary part of the last term in RHS of (\ref{eqn:C_log2})}
=- d \Big( \chi \arg g_{\HH^{\eta}_t}'(z) \Big),
\nonumber\\
\mbox{\rm (ii)} \quad & 
\mbox{if} \quad 
F_i(\Y)=4 \sum_{\substack{1 \leq j \leq N, \cr j \not=i}}
\frac{1}{Y_i(t)-Y_j(t)}, \quad i=1, \dots, N,
\label{eqn:DysonC1}
\\
& \qquad 
\, \mbox{then 
the second term in RHS of (\ref{eqn:C_log2})}=0.
\end{align}

\begin{rem}
\label{thm:Remark3_2}
If we choose the drift term as (\ref{eqn:Dyson3}), 
the system of SDEs of the driving process of the multiple SLE
(\ref{eqn:SDE_A1b}) is determined as
\begin{equation}
d Y_i(t)=\sqrt{\kappa} dB_i(t)+
4 \sum_{\substack{1 \leq j \leq N, \cr j \not=i}}
\frac{1}{Y_i(t)-Y_j(t)}dt, \quad t \geq 0,
\quad i=1, \dots, N, 
\label{eqn:DysonC2}
\end{equation}
As we have already claimed that
$(\sqrt{\kappa} B(t))_{t \geq 0} \law= (B(\kappa t))_{t \geq 0}$,
we perform a time change $\kappa t \to t$ and define
$\X(t)=\Y(\kappa t), t \geq 0$. Then 
we have the following set of SDEs for $(\X(t))_{t \geq 0}$,
\begin{align}
& d X_i(t)=dB_i(t)+ \frac{4}{\kappa}
\sum_{\substack{1 \leq j \leq N, \cr j \not=i}}
\frac{1}{X_i(t)-X_j(t)}dt
\nonumber\\
\iff \quad & d X_i(t)=dB_i(t)+ \frac{\beta}{2}
\sum_{\substack{1 \leq j \leq N, \cr j \not=i}}
\frac{1}{X_i(t)-X_j(t)}dt, \quad t \geq 0,
\quad i=1, \dots, N, 
\label{eqn:DysonC3}
\end{align}
with
\begin{equation}
\beta=\frac{8}{\kappa}.
\label{eqn:beta_kappa}
\end{equation}
Hence, we can say that the $N$-particle system $(\Y(t))_{t \geq 0}$
satisfying (\ref{eqn:DysonC2}) is a time change
of the Dyson model with
$\beta=8/\kappa$; DYS$_{8/\kappa}$.
\end{rem}

We introduce the following notation,
\begin{align}
\mfh_t(\cdot)
&:=-\frac{2}{\sqrt{\kappa}} 
\Im \Phi(g_{\HH^{\eta}_t}(\cdot), \Y(t)) 
-\chi \Im \log {g_{\HH^{\eta}_t}}'(\cdot)
\nonumber\\
&=-\frac{2}{\sqrt{\kappa}} 
\sum_{i=1}^N \arg (g_{\HH^{\eta}_t}(\cdot) - Y_i(t))
-\chi \arg {g_{\HH^{\eta}_t}}'(\cdot). 
\label{eqn:mfh}
\end{align}
The above observation implies the following.

\begin{prop}
\label{thm:martingale}
Assume (\ref{eqn:DysonC2}), that is,
the driving process of the multiple SLE is
given by DYS$_{8/\kappa}$. 
If 
$\displaystyle{
\chi=\frac{2}{\sqrt{\kappa}}-\frac{\sqrt{\kappa}}{2}}$,
then $(\mfh_t(\cdot))_{t \geq 0}$ is
a continuous local martingale and satisfies
\begin{equation}
d \mfh_t(z)
= 2 \sum_{i=1}^N
\Im \frac{1}{g_{\HH^{\eta}_t}(z)-Y_i(t)} d B_i(t),
\quad z \in \HH_t^{\eta}, 
\quad t \in [0, \infty).
\label{eqn:martingale2b}
\end{equation}
\end{prop}

\begin{rem}
\label{thm:chi_c}
Between the central charge $c=c_{\kappa}$ given by
(\ref{eqn:CFT}) and 
$\displaystyle{\chi=
\chi_{\kappa} := \frac{2}{\sqrt{\kappa}} - \frac{\sqrt{\kappa}}{2}
}$, the following simple relation is established,
\begin{equation}
c_{\kappa}=1- 6 \chi_{\kappa}^2.
\label{eqn:c_kappa}
\end{equation}
By this formula, for the SLE parameter $\kappa >0$, 
$c_{\kappa} \leq 1$,
and the maximum central charge
$c_{\kappa}=1$ is realized if and only if 
$\chi_{\kappa}=0 \Longleftrightarrow 
\kappa=4$.
\end{rem}
\vskip 0.3cm

Comparing Proposition \ref{thm:martingale}
with (\ref{eqn:dG}) in Proposition \ref{thm:dG},
we see that 
\begin{align*}
d \langle \mfh_{\cdot}(z), 
\mfh_{\cdot}(w) \rangle_t
&=-d G^{g_{\HH^{\eta}_t}}(z, w),
\quad z, w \in \HH^{\eta}_t, \quad t \in [0, \infty).
\end{align*}
Moreover, with 
(\ref{eqn:energy2}) we have the equality, 
\begin{align}
d \langle 
\langle \mfh_{\cdot}, f \rangle, \langle \mfh_{\cdot}, f \rangle 
\rangle_t
&=-d E^{g_{\HH^{\eta}_t}}_A(f)
\quad \mbox{for
$f \in \cC^{\infty}_{\rm c}(\HH)$ with
$\supp(f) \subset A \subset \HH$}.
\label{eqn:energy_B}
\end{align}

\subsubsection{Extended GFF-valued stochastic process
and its stationarity}
\label{sec:GFF_process}

Put 
\[
\chi=\frac{2}{\sqrt{\kappa}}-\frac{\sqrt{\kappa}}{2}, 
\]
and consider the following time-evolution of
extended GFF, 
\begin{align}
H_t
&: = {g_{\HH^{\eta}_t}}_* H
+\mfh_t, 
\quad t \geq 0.
\label{eqn:GFF_process1}
\end{align}

\begin{figure}[t]
\begin{center}
\includegraphics[scale=.40]{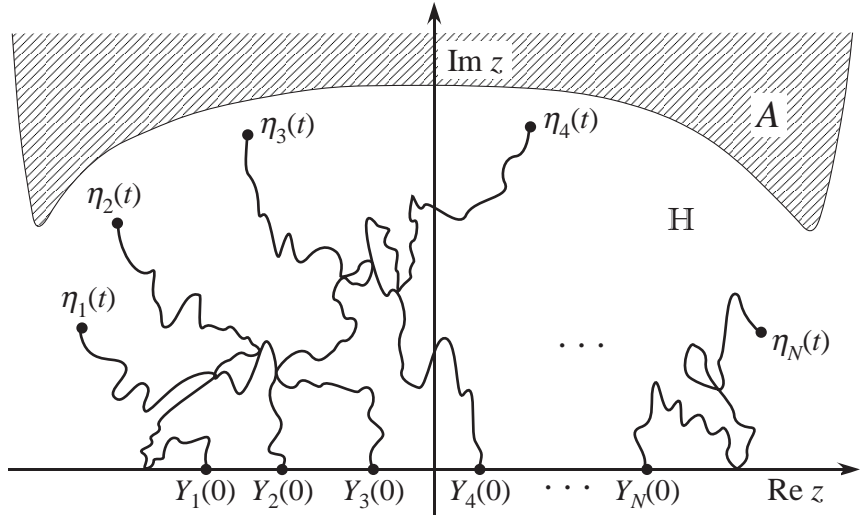}
\caption{
Consider a domain $A \subset \HH$ such that
there is a gap between $A$ and the real axis $\R$.
For such $A$, we define 
$\tau_A := \sup \Big\{ t \geq 0 : 
A \subset \HH^{\eta}_t \Big\} >0$ 
and consider $H_t$ 
during the time period 
$t \in [0, \tau_A]$.
}
\end{center}
\label{fig:mSLEA}       
\end{figure}

\begin{rem}
\label{thm:A}
By Definition \ref{thm:def_GFF}, 
a GFF added by a harmonic function has been called 
an extended GFF.
The additional harmonic function in
(\ref{eqn:GFF_process1}) is
\[
\arg(g_{\HH^{\eta}_t}(\cdot)-Y_i(t))
=\Im \log (g_{\HH^{\eta}_t}(\cdot)-Y_i(t)),
\]
which is not well-defined on the SLE curves $\eta_{i}$,
$i=1, \dots, N$.
Hence, $H_t$ can be regarded as an extended GFF
not in the whole $\HH$, but only on  $\HH^{\eta}_t$,
$t \geq 0$, where we remember that 
$\HH^{\eta}_t, t > 0$ is defined by (\ref{eqn:HetaB2}), 
that is 
\[
\mbox{
$\HH^{\eta}_t :=$ the unbounded component
of $\displaystyle{\HH \backslash \bigcup_{i=1}^N \eta_i(0, t], 
\quad t \geq 0}$}. 
\]
We consider a domain $A \subset \HH$ such that
there is a gap between $A$ and the real axis $\R$.
In this case, if we set $v_A:=\inf\{\Im z: z \in A\}$,
then $v_A \geq {^{\exists}\delta} > 0$. 
For such $A$, we define an $(\cF_t)_{t \geq 0}$--stopping time, 
\begin{equation}
\tau_A := \sup \Big\{ t \geq 0 : 
A \subset \HH^{\eta}_t \Big\}.
\label{eqn:tauA}
\end{equation}
Then, $\tau_A >0$ and 
during the time period 
$t \in [0, \tau_A]$,
$A$ is separated from any SLE curve and
the multiple SLE $g_{\HH^{\eta}_t}$ is well defined 
in $A$, and hence 
(\ref{eqn:GFF_process1}) 
can be regarded as an extended GFF-valued process.
See Fig.\ref{fig:mSLEA}.
Note that test functions $f \in \cC^{\infty}_{\rm c}(\HH)$ 
should be considered with the condition
$\supp(f) \subset A$.
\end{rem}

Under the restrictions mentioned in
Remark \ref{thm:A}, 
the following Proposition asserts that
a kind of stationarity is established if
we couple the multiple SLE and GFF properly.

\begin{prop}
\label{thm:GFF_mSLE1}
Assume the following.
\item{\bf (A)} \,
$(\mfh_t)_{t \geq 0}$ is a continuous local martingale and
its quadratic covariation satisfies the equality, 
\begin{equation}
d \langle \mfh_{\cdot}(z), 
\mfh_{\cdot}(w) \rangle_t
=-d G^{g_{\HH^{\eta}_t}}(z, w),
\quad z, w \in \HH^{\eta}_t, \quad t \in [0, \infty).
\label{eqn:A2}
\end{equation}
Let $A$ be an arbitrary domain in $\HH$ such that 
$v_A \geq {^{\exists} \delta} > 0$.
With the $(\cF_t)_{t \geq 0}$--stopping time $\tau_A$ 
given by (\ref{eqn:tauA}), 
put $0 < T < \tau_A$.
Then for any $f \in \cC^{\infty}_{\rm c}(\HH)$
with $\supp(f) \subset A$, 
\begin{equation}
\bra H_0, f \ket \law= \bra H_t, f \ket,
\quad t \in [0, T].
\label{eqn:stationary1}
\end{equation}
\end{prop}
\noindent{\it Proof.} 
Introduce a real parameter $\theta \in \R$ 
and consider the characteristic function for 
$\bra H_t, f \ket$, 
$\E[ e^{\sqrt{-1} \theta \bra H_t, f \ket}]$.
Here $\E$ denotes the expectation with respect to
the joint probability law $\P$ of
the multiple SLE and GFF.
$\mfh_{t}(\cdot)$ is $\cF_{t}$--measurable, 
\[
\E \Big[ e^{\sqrt{-1} \theta \bra H_t, f \ket} \Big]
=\E \Bigg[ 
\E 
\Big[
\exp( \sqrt{-1} \theta
\bra {g_{\HH^{\eta}_t}}_* H, f \ket ) \Big|
\cF_t \Big]
e^{\sqrt{-1} \theta \bra \mfh_{t}, f \ket}
\Bigg].
\]
Since
$\supp(f) \subset A$,
by the definition (\ref{eqn:energy}) of
$E^{g_{\HH^{\eta}_t}}_A(f)$ 
$\Var[\bra {g_{\HH^{\eta}_t}}_* H, f \ket]
=E^{g_{\HH^{\eta}_t}}_A(f)$.
Hence the formula (\ref{eqn:GFF_chara})
of Lemma \ref{thm:chara1},
which was concluded from the 
Bochner--Minlos theorem, gives
\[
\E \Big[
\exp( \sqrt{-1} \theta
\bra {g_{\HH^{\eta}_t}}_* H, f \ket ) \Big|
\cF_t \Big]
=\exp\left(
-\frac{\theta^2}{2} E^{g_{\HH^{\eta}_t}}_A(f) \right)
\]
Then we obtain 
\begin{equation}
\E \Big[ e^{\sqrt{-1} \theta \bra H_t, f \ket} \Big]
=\E \left[
\exp \left(
-\frac{\theta^2}{2} E^{g_{\HH^{\eta}_t}}_A(f)
+\sqrt{-1} \theta \bra \mfh_{t}, f \ket
\right) \right].
\label{eqn:E_mart}
\end{equation}
By It\^o's formula, 
\begin{align*}
& d \exp \left(
-\frac{\theta^2}{2} E^{g_{\HH^{\eta}_t}}_A(f)
+\sqrt{-1} \theta \bra \mfh_{t}, f \ket
\right)
\nonumber\\
& \,
= \left\{
\sqrt{-1} \theta d \bra \mfh_{t}, f \ket
-\frac{\theta^2}{2}
\left(d E^{g_{\HH^{\eta}_t}}_A(f) 
+ d \bra \bra \mfh_{\cdot}, f \ket, \bra \mfh_{\cdot}, f \ket 
\ket_t \right)
\right\}
\exp \left(
-\frac{\theta^2}{2} E^{g_{\HH^{\eta}_t}}_A(f)
+\sqrt{-1} \theta \bra \mfh_{t}, f \ket
\right).
\end{align*}
By the assumption {\bf (A)}, (\ref{eqn:energy_B}) holds
and we can conclude that
$\exp \left(
-\frac{\theta^2}{2} E^{g_{\HH^{\eta}_t}}_A(f)
+\sqrt{-1} \theta \bra \mfh_{t}, f \ket
\right)$ 
is an $\cF_t$-local martingale
and that the characteristic function 
is given by its expectation.
Therefore, the characteristic function is
time independent and given by its initial value, 
\[
\E \left[
\exp \left(
-\frac{\theta^2}{2} E^{g_{\HH}}_A(f)
+\sqrt{-1} \theta \bra \mfh_{0}, f \ket
\right) \right]
=
\E \Big[ e^{\sqrt{-1} \theta \bra H_0, f \ket} \Big].
\]
The proof of the equality (\ref{eqn:stationary1}) is
complete. 
\qed

\subsubsection{Multiple SLE/GFF coupling realized by 
DYS$_{8/\kappa}$}
\label{sec:solution}
\begin{df}[multiple SLE/GFF coupling]
\label{eqn:coupling}
For any domain $A \subset \HH$ such that 
$v_A \geq {^{\exists} \delta} > 0$,
define the $(\cF_t)_{t \geq 0}$--stopping time 
$\tau_A$ by (\ref{eqn:tauA}) and assume
$0< T < \tau_A$. 
Given the driving process
$(\Y(t))_{t \geq 0}$ for the multiple SLE,
assume the equality 
\[
\bra H_0, f \ket \law= \bra H_t, f \ket,
\quad t \in [0, T]
\]
for any $f \in \cC^{\infty}_{\rm c}(\HH)$ with
$\supp(f) \subset A$.
Then we say that the 
\textbf{multiple SLE/GFF coupling} is established.  
\end{df}
The following is the main result in this
section \cite{KK20,KK21a,KK21b}.
\begin{thm}
\label{thm:main}
Assume
\begin{equation}
\chi=\frac{2}{\sqrt{\kappa}}-\frac{\sqrt{\kappa}}{2}.
\label{eqn:chiA}
\end{equation}
Then, if and only if the driving process 
of multiple SLE $(\Y(t))_{t \geq 0}$ 
is given by the solution of the SDEs
\begin{equation}
dY_i(t) = \sqrt{\kappa} d B_i(t)
+ 4 \sum_{\substack{1 \leq j \leq N, \cr j \not= i}}
\frac{1}{Y_i(t)-Y_j(t) } dt, 
\quad t \geq 0, \quad
i=1, \dots, N,
\label{eqn:DysonA1_B}
\end{equation}
the multiple SLE/GFF coupling is established. 
In other words, under (\ref{eqn:chiA}),
the driving process of multiple SLE
is uniquely determined to be 
the Dyson model with
\begin{equation}
\beta=\frac{8}{\kappa},
\label{eqn:beta_kappaB}
\end{equation}
DYS$_{8/\kappa}$, 
if the multiple SLE/GFF coupling holds.
\end{thm}
\noindent{\it Proof.} 
Assume that $(\mfh_t)_{t \geq 0}$ given by 
(\ref{eqn:mfh}) satisfies the assumption {\bf (A)} in 
Proposition \ref{thm:GFF_mSLE1}.
Then by Propositions \ref{thm:martingale} and \ref{thm:GFF_mSLE1},
the present theorem will be proved if 
we can derive DYS$_{8/\kappa}$ 
as $(\Y(t))_{t \geq 0}$ under the conditions.
Owing to the implicit function theorem, 
we can say that each $Y_i(t), t \geq 0,  i=1, \dots, N$ 
is a $\cC^{\infty}$ function of
a continuous local martingale $\mfh_t(z_j)$, 
multiple SLE $g_t(z_j)$ solving
(\ref{eqn:mSLE1}), and 
$\log g'_t(z_j)$, $j=1, \dots, N$,
which satisfy
\[
\frac{d}{dt} \log g_t'(z)=
-\sum_{i=1}^N \frac{2}{(g_t(z)-Y_i(t))^2},
\quad z \in \HH^{\eta}_t, \quad t \geq 0.
\]
Then by It\^o's formula they are semi-martingales
and given in the forms,
\begin{equation}
Y_i(t)=M_i(t)+F_i(\Y_t), \quad t \geq 0, \quad
i=1, \dots, N, 
\label{eqn:semi_martingale}
\end{equation}
where $M_i$, $i=1, \dots,N$
denote the local martingale parts
and $F_i$, $i=1, \dots, N$ do
the finite-variation parts.
Using these representations, as already seen in Section
\ref{sec:complex_log},
It\^o's formula gives
\begin{align}
d \mfh_t(z)
&= \Im \sum_{i=1}^N \frac{1}{(g_t(z)-Y_i(t))^2}
\left\{ \left( -\frac{4}{\sqrt{\kappa}}+2 \chi \right) dt
+\frac{1}{\sqrt{\kappa}} 
d \langle M_i(\cdot), M_i(\cdot) \rangle_t \right\}
\nonumber\\
& \,
+\frac{2}{\sqrt{\kappa}} \Im \frac{1}{g_t(z)-Y_i(t)}
\left( d F_i(\Y(t)) - 
\sum_{\substack{1 \leq j \leq N, \cr j \not=i}} 
\frac{4}{Y_i(t)-Y_j(t)} dt \right)
\nonumber\\
& \,
+ \frac{2}{\sqrt{\kappa}} \Im
\sum_{i=1}^N \frac{1}{g_t(z)-Y_i(t)} d M_i(t),
\quad z \in \HH^{\eta}_t, \quad t \geq 0. 
\label{eqn:key2}
\end{align}
By the assumption {\bf (A)}, the first and the second lines in RHS
should vanish.
The each term in the first line
represents the contribution from the pole of
second order at $g_t(z)=Y_i(t)$, $i=1, \dots, N$.
Under (\ref{eqn:chiA}), vanishing of all these contribution
implies 
\[
d \langle M_i(\cdot), M_i(\cdot) \rangle_t
= \kappa dt, \quad i=1, \dots, N.
\]
The each term in the second line
represents the contribution from the pole of
first order at $g_t(z)=Y_i(t)$, $i=1, \dots, N$,
and vanishing of all these contribution
implies 
\[
d F_i(\Y_t) =
\sum_{\substack{1 \leq j \leq N, \cr j \not=i}} 
\frac{4}{Y_i(t)-Y_j(t)} dt,
\quad t \geq 0, \quad i=1, \dots, N. 
\]
Moreover, the quadratic covariations of 
$(\mfh_t(\cdot))_{t \geq 0}$ is calculated as
\begin{align*}
d \langle \mfh_{\cdot}(z), \mfh_{\cdot}(w) \rangle_t
&= -d G^{g_{\HH^{\eta}_t}}(z, w)
\nonumber\\
& \, 
+\frac{4}{\kappa} \sum_{\substack{1 \leq i, j \leq N, \cr
i\not=j}}
\Im \left(\frac{1}{g_t(z)-Y_i(t)} \right)
\Im \left( \frac{1}{g_t(w)- Y_j(t)} \right)
d \langle M_i(\cdot), M_j(\cdot) \rangle_t,
\end{align*}
$z, w \in \HH^{\eta}_t, t \geq 0$.
By the assumption {\bf (A)}, 
\[
d \langle M_i(\cdot), M_j(\cdot) \rangle_t=0,
\quad t \geq 0, \quad 1 \leq i \not=j \leq N.
\]
Hence the proof is complete.
\qed

\begin{rem}
\label{thm:3phases}
Remember the relation (\ref{eqn:DbetaB}) 
in Section \ref{sec:multi}, 
\[
\beta=D-1,
\]
and the relation (\ref{eqn:D_kappa}) in Section
\ref{sec:SLE},
\[
\kappa=\frac{4}{D-1}.
\]
If we simply combine these two relations,
we will have
\[
\beta=\frac{4}{\kappa}.
\]
The relation (\ref{eqn:beta_kappaB}) which establishes 
the multiple SLE/GFF coupling is 
\textit{different} from such a
simple-minded result.
The multiple SLE driven by 
DYS$_{\beta}$ with the true relation (\ref{eqn:beta_kappaB}),
$\beta=8/\kappa$, inherits many properties from the original
SLE$_{\kappa}$ with a single SLE curve.
Actually, we have proved that
our \textbf{multiple SLE$_{\kappa}$} also shows 
\textbf{phase transitions} 
at $\kappa_{\rm c}=4$ and
$\overline{\kappa}_{\rm c}=8$ 
\cite{KK21b}.
\begin{thm}
\label{thm:3phases_multi}
For each $i=1, \dots, N$,
the limit $\eta_i(t)=\lim_{\varepsilon \downarrow 0}
g_{\HH^{\eta}_t}(Y_i(t)+\sqrt{-1} \varepsilon)$ exists
for all $t \geq 0$ and
$\lim_{t \to \infty} \eta_i(t)=\infty$ with probability one.
Moreover, the following 
\textbf{three phases} are observed.
\begin{description}
\item{\rm (a)} \,
If $0 < \kappa \leq \kappa_{\rm c}=4$,
$\eta_i(0, \infty), i=1, \dots, N$ are simple disjoint curves
such that $\eta_i \subset \HH, i=1, \dots, N$ with
probability one.

\item{\rm (b)} \,
If $\kappa_{\rm c}=4 < \kappa < \overline{\kappa}_{\rm c}=8$,
$\eta_i, i=1, \dots, N$ are continuous curves with probability one,
and they intersect themselves and $\R$ with positive probability.

\item{\rm (c)} \,
If $\kappa=8$,
$\eta_i, i=1, \dots, N$ are space filling continuous 
curves with probability one.
\end{description}
\end{thm}
\end{rem}

\subsection{Exercises 3}
\label{sec:exercises3}
\subsubsection{Exercise 3.1}
\label{sec:ex3_1}
Let $\alpha \in (0, 1)$ and
\begin{equation}
\kappa=\kappa(\alpha)=
\frac{4(1-2 \alpha)^2}{\alpha(1-\alpha) }.
\label{eqn:ex3_V1}
\end{equation}
Then consider the case such that
\begin{equation}
V(t)=\begin{cases}
\sqrt{\kappa t}, & \mbox{if $\alpha \leq 1/2$}, \cr
-\sqrt{\kappa t}, & \mbox{if $\alpha > 1/2$}.
\end{cases}
\label{eqn:ex3_V2}
\end{equation}
\begin{description}
\item{(1)} \,
Show that
the inverse of $g_t$ is solved as
\begin{equation}
g_{\HH^{\eta}_t}^{-1}(z) 
= \left( z + 2 \sqrt{\frac{\alpha}{1-\alpha}} \sqrt{t} \right)^{1-\alpha}
\left(z - 2 \sqrt{\frac{1-\alpha}{\alpha}} \sqrt{t} \right)^{\alpha},
\quad z \in \HH.
\label{eqn:ex3_V3}
\end{equation}

\item{(2)} \,
Derive that the following formula for the slit,
\begin{equation}
\eta(t)=g_{\HH^{\eta}_t}^{-1}(V(t))
=2 \left( \frac{1-\alpha}{\alpha} \right)^{1/2-\alpha}
e^{\sqrt{-1} \alpha \pi} t^{1/2},
\quad t \geq 0.
\label{eqn:ex3_V4}
\end{equation}
\end{description}
The slit grows from the origin along a straight line in $\HH$
which makes an angle $\alpha \pi$ with respect to the positive
direction of the real axis.
When $\alpha=1/2$, this is reduced to Example \ref{thm:SLE0}.

\subsubsection{Exercise 3.2}
\label{sec:ex3_2}

Assume that $D, D' \subsetneq \C$ are simply connected domains
and let 
\begin{equation}
\label{eqb:ex3_CI1}
\varphi: D' \to D \quad \mbox{conformal map}.
\end{equation}
\begin{description}
\item{(1)} \,
We write the conformal map $\varphi$ using the
real functions (harmonic functions) $u, v$ as
\begin{equation}
\varphi(z)=u(x, y)+\sqrt{-1} v(x, y).
\label{eqn:ex3_CI2}
\end{equation}
Using the \textbf{Cauchy--Riemann equations},
\begin{equation}
\frac{\partial u}{\partial x}=\frac{\partial v}{\partial y}, \quad
\frac{\partial u}{\partial y}=-\frac{\partial v}{\partial x},
\label{eqn:ex3_CI3}
\end{equation}
show that the Jacobian of $\varphi$ is given by
\begin{equation}
\frac{\partial (u,v)}{\partial (x,y)}
=\frac{\partial u}{\partial x} \frac{\partial v}{\partial y}
-\frac{\partial u}{\partial y} \frac{\partial v}{\partial x}
=\left(\frac{\partial u}{\partial x}\right)^2
+\left( \frac{\partial u}{\partial y}\right)^2
\label{eqn:ex3_CI4}
\end{equation}

\item{(2)} \,
Assume that $f, g \in \cC_{\rm c}^{\infty}(D)$
Using the chain rule of differentials and the
Cauchy--Riemann equations (\ref{eqn:ex3_CI3}),
show the equality,
\begin{equation}
(\nabla f \circ \varphi)(z) \cdot (\nabla g \circ \varphi)(z)
=\left( \frac{\partial f}{\partial u} \frac{\partial g}{\partial u}
+\frac{\partial f}{\partial v} \frac{\partial g}{\partial v} \right)
\left\{ \left(\frac{\partial u}{\partial x}\right)^2
+\left( \frac{\partial u}{\partial y}\right)^2 \right\}.
\label{eqn:ex3_CI5}
\end{equation}

\item{(3)} \,
Prove the equality,
\begin{equation}
\int_{D} (\nabla f)(z) \cdot (\nabla g)(z) \mu(d z)
=\int_{D'} (\nabla (f \circ \varphi))(z) \cdot 
(\nabla (g \circ \varphi))(z) \mu(d z)
\label{eqn:ex3_CI6}
\end{equation}
for $f, g \in \cC_{\rm c}^{\infty}(D)$.
\end{description}

\subsubsection{Exercise 3.3}
\label{sec:ex3_3}

Assume that
\begin{equation}
\E^{\GFF}[ e^{\sqrt{-1} \theta \bra H, f \ket_{\nabla}} ]
=e^{-\theta^2 \|f\|_{\nabla}^2/2}, 
\quad \theta \in \R, \quad f \in W(D).
\label{eqn:ex3_GFF1}
\end{equation}
\begin{description}
\item{(1)} \,
From (\ref{eqn:ex3_GFF1}), derive the following
equalities,
\begin{align}
\Var[\bra H, f \ket_{\nabla} ]
&:=\E^{\GFF} \Big[\bra H, f \ket_{\nabla}^2 \Big]
=\|f\|_{\nabla}^2, \quad f \in W(D),
\label{eqn:ex3_GFF2}
\\
\Cov[\bra H, f \ket_{\nabla},
\bra H, g \ket_{\nabla} ]
&:= \E^{\GFF} \Big[
\bra H, f \ket_{\nabla} \bra H, g \ket_{\nabla} \Big]
=\bra f, g \ket_{\nabla}, \quad
f, g \in W(D). 
\label{eqn:ex3_GFF3}
\end{align}

\item{(2)} \,
Let $\varphi: D' \to D$ be a conformal map,
and denote its pull-back as
\begin{equation}
\varphi^{\ast} : W(D) \ni f \mapsto f \circ \varphi \in W(D').
\label{eqn:ex3_GFF4}
\end{equation}
Then the equality
\begin{equation}
\int_{D} (\nabla f)(z) \cdot (\nabla g)(z) \mu(d z)
=\int_{D'} (\nabla (f \circ \varphi))(z) \cdot 
(\nabla (g \circ \varphi))(z) \mu(d z)
\label{eqn:ex3_GFF5}
\end{equation}
in Lemma \ref{thm:conformal_inv}
is simply written as
\begin{equation}
\bra \varphi^{\ast} f, \varphi^{\ast} g \ket_{\nabla}
=\bra f, g \ket_{\nabla}.
\label{eqn:ex3_GFF6}
\end{equation}
Prove the \textbf{conformal invariance of GFF}, 
\begin{equation}
\Cov[\bra \varphi_{\ast} H, f \ket_{\nabla}, 
\bra \varphi_{\ast} H, g \ket_{\nabla}]
=\Cov[\bra H, f \ket_{\nabla}, 
\bra H, g \ket_{\nabla} ].
\label{eqn:ex3_GFF7}
\end{equation}
\end{description}

\vskip 1cm
\noindent{\bf Acknowledgements} \,
This manuscript was prepared for the 
lectures given in
the 4th ZiF Summer School
`Randomness in Physics and Mathematics: 
From Integrable Probability to Disordered Systems'
held at 
ZiF--Center for Interdisciplinary Research, Bielefeld University,
Germany, from 1st to 13th August 2022,
which was organized by Gernot Akemann
and Friedrich G\"{o}tze.
The present author would like to thank
Gernot and Friedrich for kind hospitality 
and well-organization of the summer school.
The lectures were based on the joint work with
Hideki Tanemura (Keio University),
Tomoyuki Shirai (Kyushu University)
and Shinji Koshida (Aalto University). 
The author is grateful to them for fruitful collaborations.
The present study has been supported by
the Grant-in-Aid for Scientific Research (C) (No.19K03674),
(B) (No.18H01124), (S) (No.16H06338),
and (A) (No.21H04432)
of Japan Society for the Promotion of Science (JSPS).
This work is also supported by JSPS Grant-in-Aid 
for Transformative
Research Areas (A) JP22H05105. 


\end{document}